\documentclass[12pt, reqno, english]{amsart}  
\usepackage[utf8]{inputenc}
\usepackage[T1]{fontenc}
\usepackage{amsmath,amsthm}
\usepackage{amsfonts,amssymb}
\usepackage{mathtools}  
\usepackage[colorlinks=true,urlcolor=blue,linkcolor=red,citecolor=magenta]{hyperref}
\usepackage[margin=1.25in]{geometry}
\usepackage{xcolor}
\usepackage[capitalise]{cleveref}
\usepackage{pgfplots}
\pgfplotsset{compat=1.18}
\usepackage{tikz}
\usepackage{tikz-cd}

\theoremstyle{plain}
\newtheorem{theorem}{Theorem}[section]
\newtheorem*{theorem*}{Theorem}
\newtheorem{lemma}[theorem]{Lemma}
\newtheorem*{lemma*}{Lemma}
\newtheorem{corollary}[theorem]{Corollary}

\theoremstyle{definition}
\newtheorem{definition}[theorem]{Definition}

\newtheorem{goal}[theorem]{Goal}
\newtheorem*{goal*}{Goal}
\theoremstyle{remark}
\newtheorem{remark}[theorem]{Remark}
\newtheorem{example}[theorem]{Example}

\numberwithin{equation}{section}

\newcommand{\kmax}{\text{kmax}}

\newcommand{\N}{\mathbb{N}}

\newcommand{\R}{\mathbb{R}}

\newcommand{\Z}{\mathbb{Z}}

\newcommand{\im}{\mathrm{im}}
\newcommand\abs[1]{\left|#1\right|}
\newcommand\norm[1]{\lVert#1\rVert}
\newcommand{\supp}{\mathrm{supp}}

\newcommand{\eps}{\varepsilon}

\DeclareMathOperator{\vol}{vol}
\DeclareMathOperator{\Id}{Id}

\newif\ifdraft
\drafttrue 

\begin{document}

\title[Explainable TDA using persistence heatmaps]{Explainable topological data analysis using persistence heatmaps}

\author{Peter Bubenik}
\address[PB]{Department of Mathematics, University of Florida, Gainesville, FL 32611, USA}
\email{peter.bubenik@ufl.edu}

\author{Alexander Wagner}
\email{alex.y.wagner@gmail.com}

\author{Himanshu Yadav}
\address[HY]{Carlson School of Management, University of Minnesota, Minneapolis, MN 55455, USA}
\email{hyadav@umn.edu}


\begin{abstract}
  Topological data analysis (TDA) leverages tools from algebraic topology to aid in various machine learning tasks. 
  Numerous TDA constructions are provably stable in the sense that changes in the input produce linearly bounded changes in the output, with a specified bound.
  We use representative cycles, which are unstable TDA constructions, to produce stable visualizations to aid in explaining TDA. 
  For example, we produce stable heatmaps on images containing the data such that summing the values of the pixels gives the value of a learned regression function.
\end{abstract}

\maketitle

\ifdraft
\fi

\section{Introduction}
\label{sec:introduction}

Topological data analysis (TDA) uses topological methods to summarize and learn from the ``shape'' of data~\cite{10.3389/frai.2021.667963}. 
For example, TDA may be used for classification~\cite{MR3464008} and regression~\cite{MR4063201}.
Data is encoded in a simplicial complex together with real values associated to its simplices. 
Then one computes the persistence diagram, which describes how the topology of the simplicial complex changes if one builds it over time according to the weights of the simplices.
Next, one applies a feature map, mapping persistence diagrams to elements of a Hilbert space or Banach space. 
With labeled training data, one may then construct a classifier (e.g. using SVM). 
For testing data, one obtains class predictions by applying the classifier to the feature map of their persistence diagrams. 
We summarize these steps in the following, which we call the standard TDA pipeline.

\begin{equation} \label{eq:pipeline}
\framebox[1.2\width]{Data} 
\to
\fbox{%
 \begin{minipage}{8em}
    \centering
    Filtered\\
    simplicial complex
  \end{minipage}%
}
\to
\fbox{%
 \begin{minipage}{5em}
    \centering
    Persistence\\
    diagram
  \end{minipage}%
}
\to
\fbox{%
 \begin{minipage}{3.5em}
    \centering
    Banach\\
    space
  \end{minipage}%
}
\to 
\framebox[1.4\width]{$\R$} 
\end{equation}

\vspace{0.5ex}

Assuming that our data is parametrized by a vector of length $m$, we represent 
the pipeline \eqref{eq:pipeline} with a function $f: \R^m \to \R$,
and we let $w \in \R^m$ represent the data.
Thus, the result of applying our pipeline to our data is the number $f(w)$.
A key property of many TDA constructions is that \eqref{eq:pipeline} is \emph{stable}, meaning that the function $f$ is Lipschitz continuous. 
That is, there is a constant $L$ such that 
for all $w,\eps \in \R^m$, 
$\abs{f(w+\eps) - f(w)} \leq L \norm{\eps}$,
where $\norm{-}$ is some norm on $\R^m$.

For examples, see the database DONUT~\cite{MR4658089}, which contains over 500 applications of TDA.
However, persistence diagrams and their feature maps are challenging for non-experts to interpret.
End users (e.g. scientists and engineers) would like an explanation of the results of these learning algorithms. 

\begin{figure}[!htb]
  \centering
  \begin{minipage}[c]{0.22\textwidth}
    \centering
    \includegraphics[width=\linewidth]{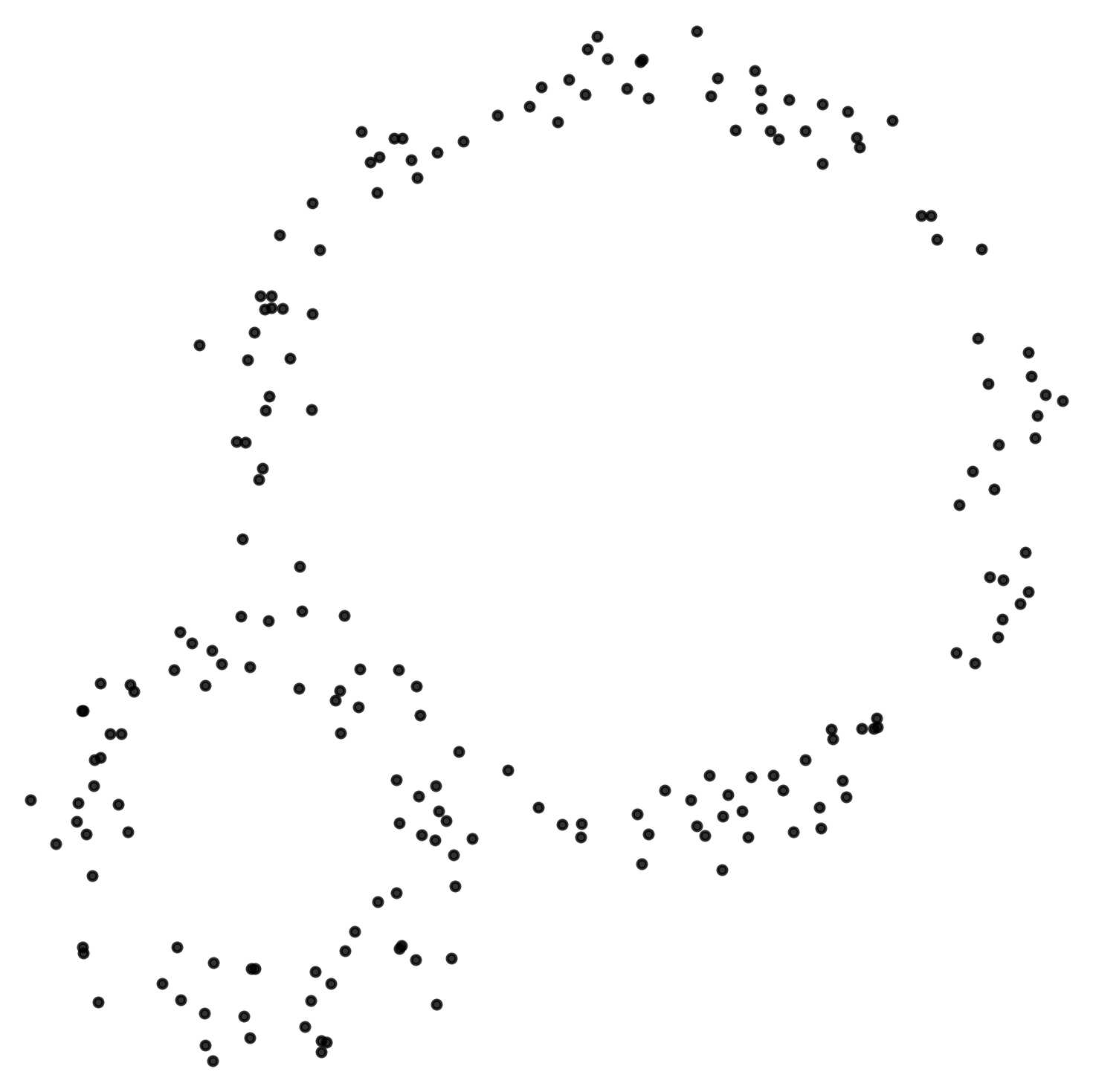}
  \end{minipage}
  \begin{minipage}[c]{0.22\textwidth}
    \centering
    \includegraphics[width=\linewidth]{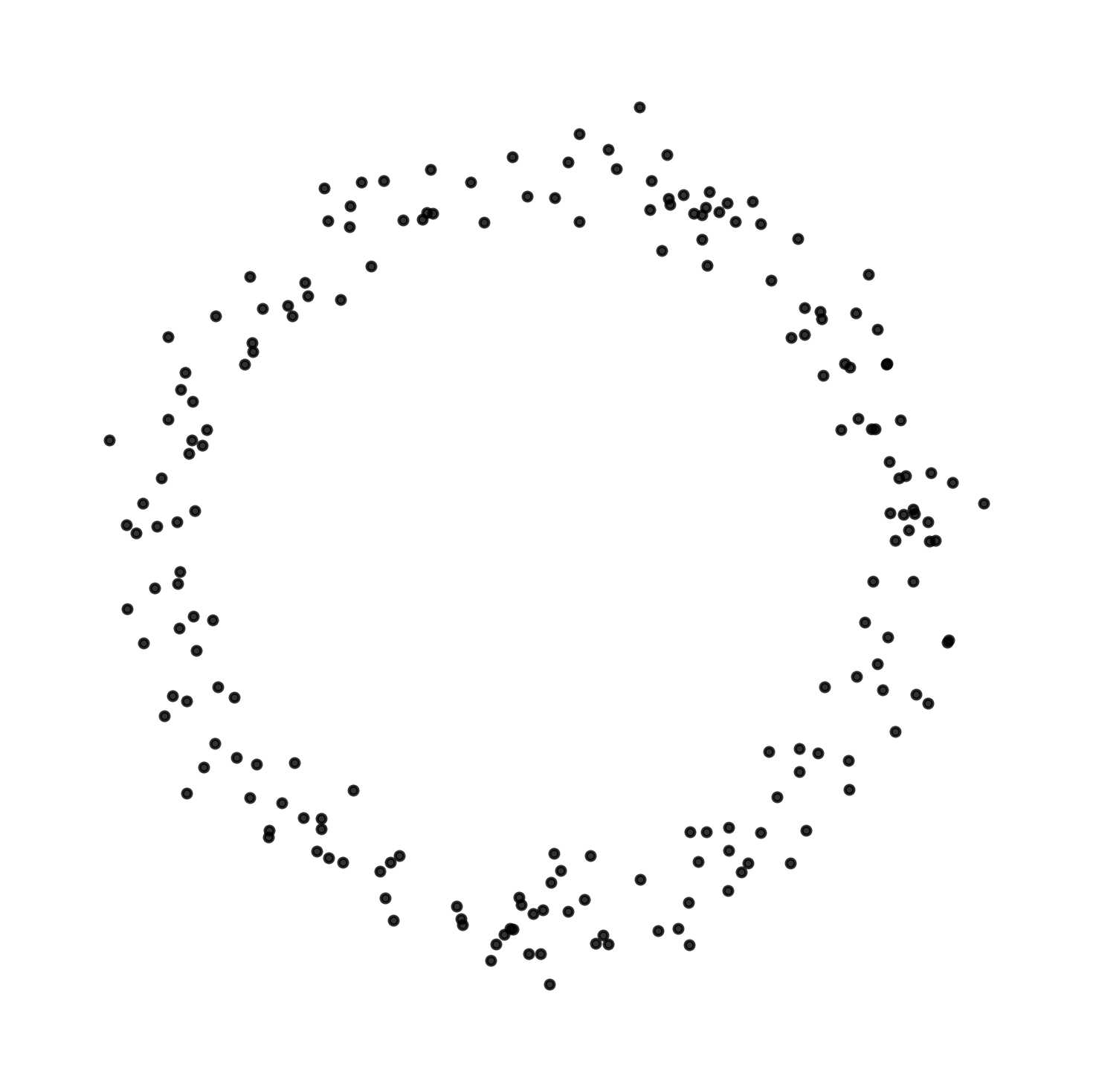}  
  \end{minipage}
  \hfill
  \begin{minipage}[c]{0.26\textwidth}
    \centering
    \includegraphics[width=\linewidth]{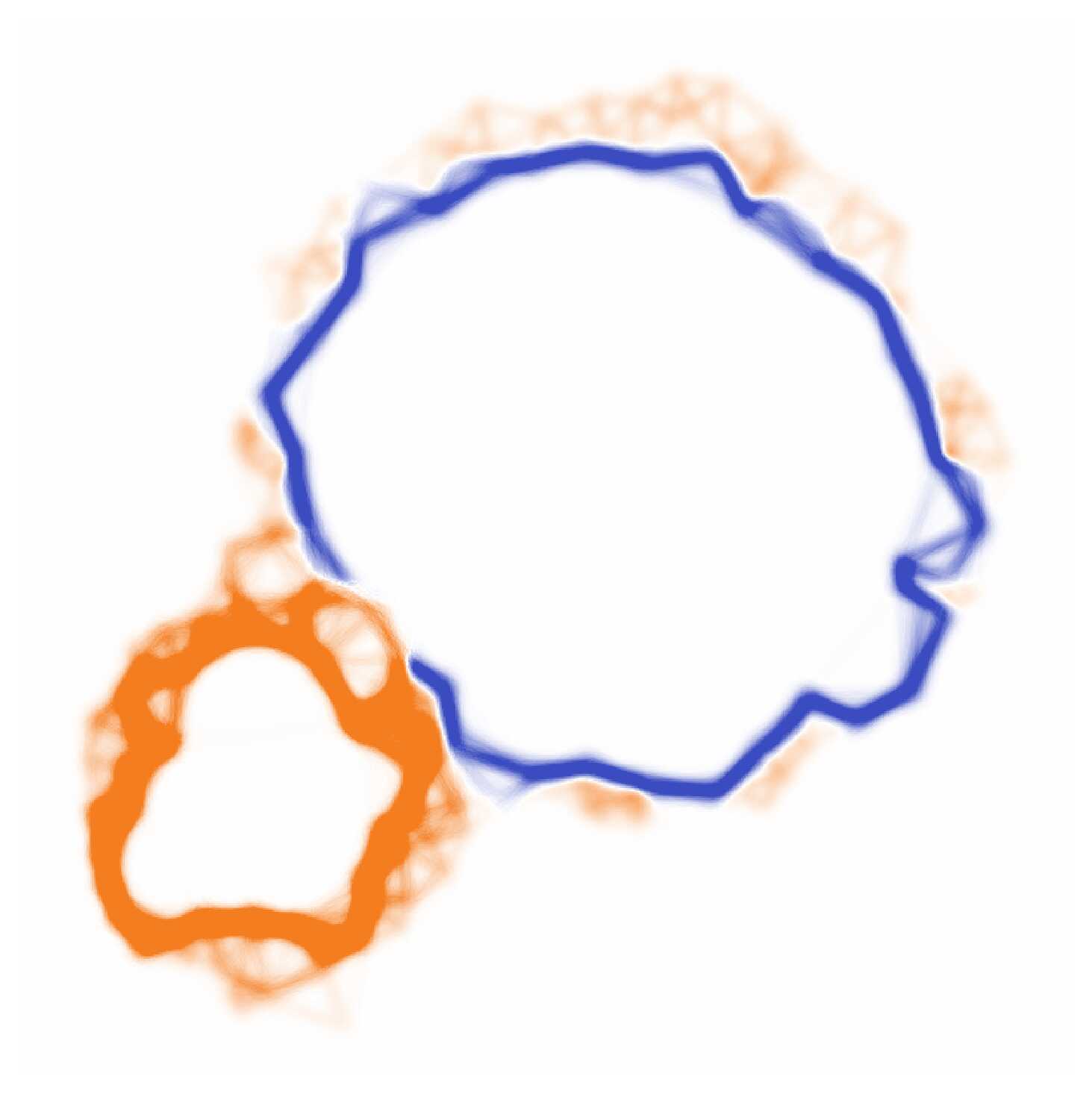}  
  \end{minipage}
  \begin{minipage}[c]{0.26\textwidth}
    \centering
    \includegraphics[width=\linewidth]{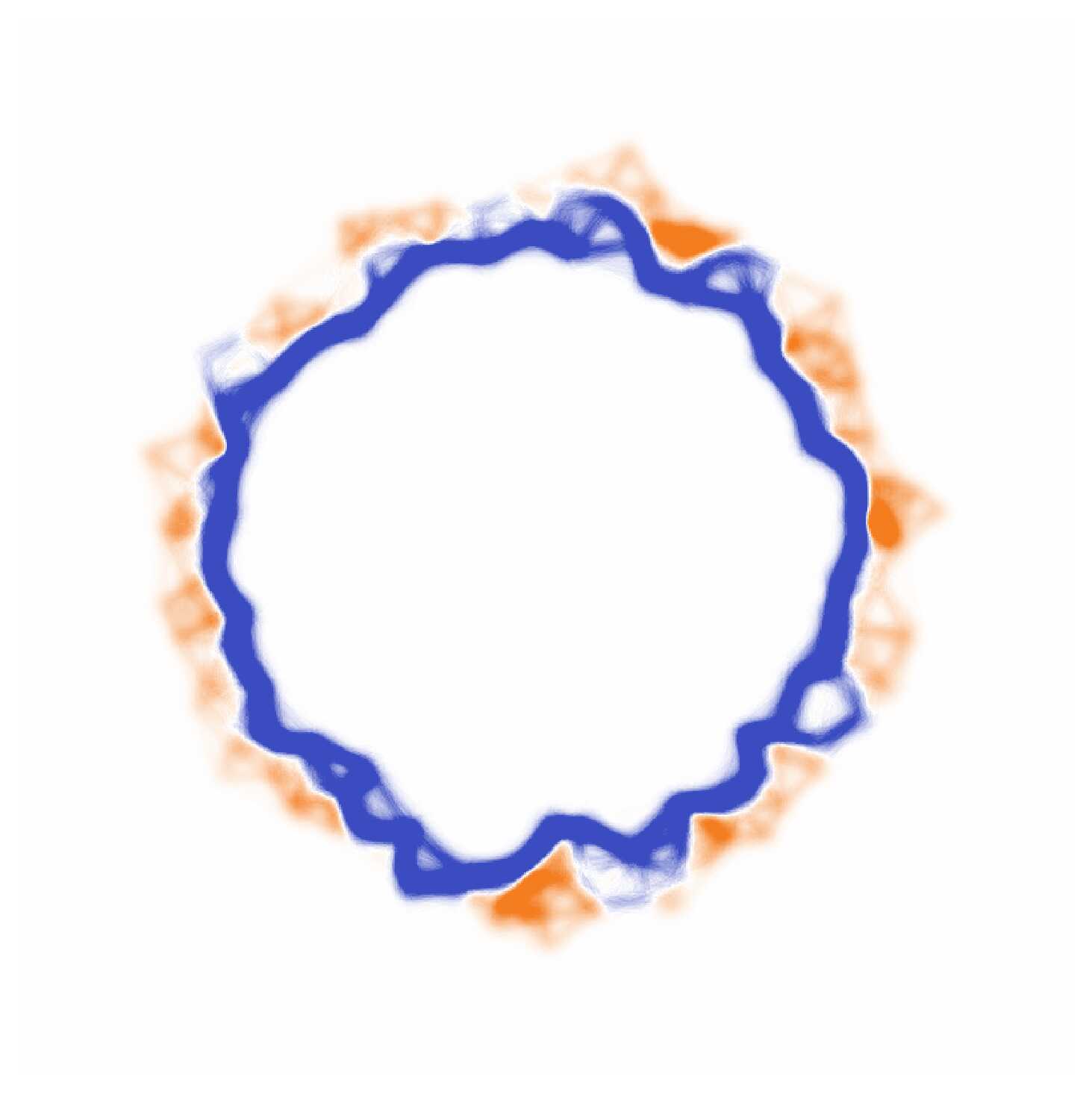}  
  \end{minipage}
    \caption{An explanation of a TDA classifier. 
    For points sampled from a double annulus (far left) and annulus (middle left),
    a classifier returns a negative expected value for the first and a positive expected value for the second, correctly identifying the class.
    We explain this result with visualizations using representative cycles for persistent homology classes (middle right and far right).
    We use an orange-to-blue color gradient in which orange corresponds to the double annulus class 
    and blue corresponds to the annulus class. 
    The sum of the expected values of the pixels equals the expected value of the classifier.
    Furthermore, the pixel values are stable with respect to perturbations of the sampled points.
    }
    \label{fig:annulus-visualization-0}
\end{figure}

\begin{figure}[!htb]
    \centering
    \includegraphics[height=35mm]{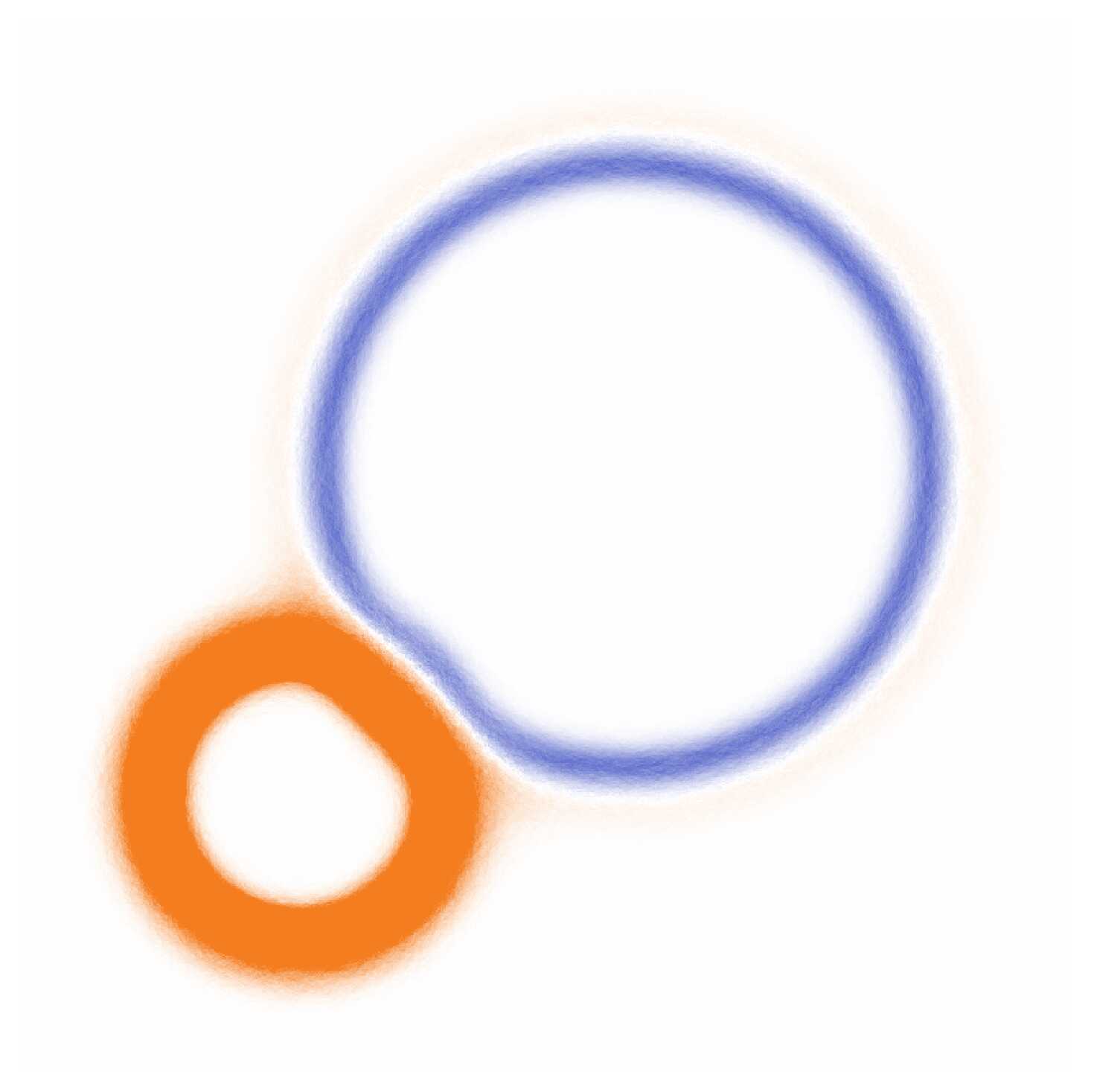}
    \quad \quad
    \includegraphics[height=35mm]{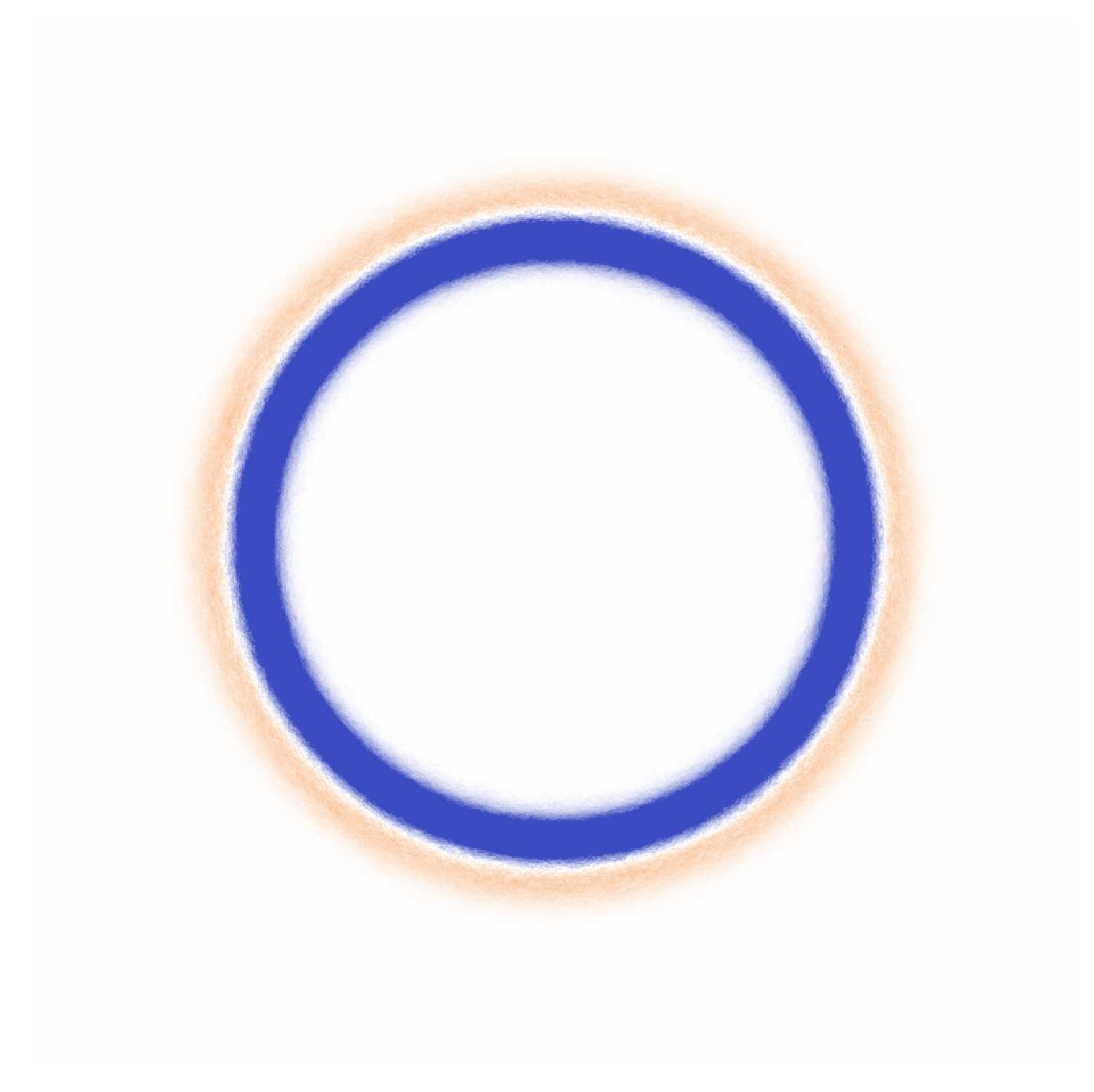}
    \caption{Instead of considering perturbations of a single sample from each class as in \cref{fig:annulus-visualization-0}, here we consider all samples from the double annulus (left) and the annulus (right).
    }
    \label{fig:annulus-visualization-family}
\end{figure}

In the computation that produces the persistence diagram, there is additional data that has potential for aiding interpretation, namely the representative cycles of persistent homology classes~\cite{csem:vineyards}. 
However, unlike the persistence diagram, representative cycles may change discontinuously with respect to perturbations of the data~\cite{BBW:2020}. 
In the present work, we develop a TDA and ML framework that produces stable visualizations using representative cycles. 
See \cref{fig:annulus-visualization-0,fig:annulus-visualization-family}.

We are motivated by the following.

\begin{goal} \label{goal:a}
    We would like to \emph{explain} the value $f(w)$ for the data $w$ by \emph{distributing} it to simplices or to pixels/voxels.
\end{goal}

If $K$ is our filtered simplicial complex then we aim to distribute $f(w)$ to the simplices of $K$. 
That is, we want for each simplex $\sigma$ in $K$ a function 
$\bar{f}_\sigma: \R^m \to \R$ such that for all $w \in \R^m$, $\sum_{\sigma \in K} \bar{f}_\sigma(w) = f(w)$.
Or, given a set $\mathcal{C}$ of pixels or voxels, we aim to distribute $f(w)$ to the pixels or voxels of $\mathcal{C}$.
That is, we want for each pixel or voxel $c$ in $\mathcal{C}$ a function 
$\tilde{f}_c: \R^m \to \R$ such that for all $w \in \R^m$, $\sum_{c \in \mathcal{C}} \tilde{f}_c(w) = f(w)$.

As a step towards \cref{goal:a},
we introduce the following variation of 
\eqref{eq:pipeline}.

\begin{equation} \label{eq:pipeline-labeled}
\framebox[1.2\width]{Data} 
\to
\fbox{%
 \begin{minipage}{4.5em}
    \centering
    Filtered\\
    simplicial\\
    complex
  \end{minipage}%
}%
\to
\fbox{%
 \begin{minipage}{5em}
    \centering
    Labeled\\
    persistence\\ 
    diagram
  \end{minipage}%
}%
\to
\fbox{%
 \begin{minipage}{4em}
    \centering
    Vector in\\
    Banach\\
    space
  \end{minipage}%
}%
\to 
\fbox{%
 \begin{minipage}{4em}
    \centering
    Vector in\\
    $\R$
  \end{minipage}%
}%
\end{equation}

\begin{remark}
    As we will explain in detail in Section~\ref{sec:annotation}, in \eqref{eq:pipeline-labeled}, each of the points in the persistence diagram is labeled by a representative cycle and we obtain vectors indexed by the simplicial chains of the filtered simplicial complex.
\end{remark}

We represent \eqref{eq:pipeline-labeled} with a function $\hat{f}: \R^m \to \R^{\hat{n}}$.
In \cref{sec:annotation}, we will construct the pipeline \eqref{eq:pipeline-labeled} so that
$\hat{f}$ has the property that for each $w \in \R^m$, the sum of the entries
of the vector $\hat{f}(w)$ equals $f(w)$.

In \cref{sec:localization}, we will use $\hat{f}$ to define a function $\bar{f}:\R^m \to \R^{\bar{n}}$ whose image vector is indexed by the simplices of the filtered simplicial complex and
which also has the property that 
for each $w \in \R^m$ the sum of the entries of the vector $\bar{f}(w)$ equals $f(w)$.
In addition,
given a set of pixels or voxels containing the data, 
we will define a function $\tilde{f}: \R^m \to \R^{\tilde{n}}$
whose image vector is indexed by the set of  pixels or voxels and
which also has the property that 
for each $w \in \R^m$ the sum of the entries of the vector $\tilde{f}(w)$ equals $f(w)$.

At this point, it seems that we have succeeded in attaining \cref{goal:a}.
However, our maps $\bar{f}_\sigma$ and $\tilde{f}_c$ are not continuous. 
Hence our explanation of $f(w)$ is not stable: arbitrarily small changes in the data $w$ can lead to large changes in the way that $f(w)$ is distributed.

Perhaps surprisingly, we will resolve this problem by \emph{adding noise}. 
Let $\eps$ be a Gaussian random variable in $\R^d$ with mean $0$.
Then $f(w+\eps)$ is a random variable in $\R$.
We will shift our focus from $f(w)$ to $E(f(w+\eps)$, where $E$ denotes expectation.
Note that in practice it will often be the case that $E(f(w+\eps))$ exists and is close to $f(w)$.

\begin{goal}
    Explain the value $E(f(w+\eps))$ for the data $w$ by distributing it to simplices or pixels/voxels in a stable way.
\end{goal}

With the following results, we attain this goal.
Note that since the expectation of a multivariate function is given coordinatewise, 
we have that for all $\sigma \in K$ and for all $c \in \mathcal{C}$, 
$E(\bar{f})_\sigma = E(\bar{f}_\sigma)$ and $E(\tilde{f})_c = E(\tilde{f}_c)$.

\begin{theorem*}[\cref{thm:distribution}]
    $ E(f(w+\eps)) = \sum_{\sigma \in K} E(\bar{f}_\sigma(w+\eps)) = \sum_{c \in \mathcal{C}} E(\tilde{f}_c(w+\eps))$.
\end{theorem*}

\begin{theorem*}[\cref{thm:slln,thm:clt}]
    For all $\sigma \in K$ and for all $c \in \mathcal{C}$, 
    $\bar{f}_\sigma(w+\eps)$ and $\tilde{f}_c(w+\eps)$ satisfy 
    the strong law of large numbers and central limit theorem.
\end{theorem*}

Let $\kappa$ denote the distribution of the Gaussian random variable $\eps$.

\begin{lemma*}[\cref{sec:convolution}]
    For all $\sigma \in K$ and for all $c \in \mathcal{C}$, 
    $E(\bar{f}_\sigma(w+\eps)) = (\bar{f}_\sigma * \kappa)(w)$ and
    $E(\tilde{f}_c(w+\eps)) = (\tilde{f}_c * \kappa)(w)$, 
    where $*$ denotes convolution.
\end{lemma*}

\begin{theorem*}[\cref{thm:stability-gaussian}]
    For all $\sigma \in K$ and for all $c \in \mathcal{C}$, 
    $E(\bar{f}_\sigma(w+\eps))$ and
    $E(\tilde{f}_c(w+\eps))$ are Lipschitz
    as functions in $w$ from $\R^m$ to $\R$.
\end{theorem*}

\begin{remark}
    We also have stability for other kernels (\cref{sec:stability}), such as the triangular kernel (\cref{thm:stability-triangular}) and the Epanechnikov kernel (\cref{thm:stability-epanechnikov}).
\end{remark}

\begin{theorem*}[\cref{thm:empirical-mean}]
    Let $\eps_1,\eps_2,\ldots,\eps_n$ be sampled according to the distribution $\kappa$.
    Then 
    $\frac{1}{n} \sum_{i=1}^n \bar{f}(w + \eps_i)$ converges to 
    $E(\bar{f}(w+\eps))$ and
    $\frac{1}{n} \sum_{i=1}^n \tilde{f}(w + \eps_i)$ converges to 
    $E(\tilde{f}(w+\eps))$.
\end{theorem*}

In \cref{fig:annulus-visualization-0}, we have taken a single sample, which we denote $w$, from each class. 
We add Gaussian noise independently to each coordinate of each of the points. 
Then we apply our function $\tilde{f}$.
For each pixel, $c$, this results in a random variable, $\tilde{f}_c(w+\eps)$.
It satisfies the strong law of large numbers and the central limit theorem. 
Its expectation is given by convolution with $\kappa$.
This expectation is Lipschitz continuous: if we perturb the point cloud by some small amount, then this expectation will change by at most a fixed scalar multiple of this small amount.
Since it is not feasible to determine the expectation analytically, we approximate it using the empirical mean. We plot the empirical mean for each pixel using a color gradient, resulting in the images in the figure.
The sum of the values of the pixels equals the empirical mean value of the classifier.

Let $W$ be a random variable with values in $\R^d$.
For example, $W$ could be given by a scientific experiment.

\begin{theorem*}[\cref{thm:empirical-mean-family}]
    Let $w_1,w_2,\ldots,w_n$ be samples from $W$.
    Then $\frac{1}{n} \sum_{i=1}^n \bar{f}(w_i)$ converges to $E(\bar{f}(W))$ and
    $\frac{1}{n} \sum_{i=1}^n \tilde{f}(w_i)$ converges to $E(\tilde{f}(W))$.
\end{theorem*}

See \cref{fig:annulus-visualization-family}, where $W$ is given by sampling points from each class.

\begin{remark}
    In addition to using representative cycles, our methods apply equally well to birth simplices, death simplices, and bounding chains.
    For example, see \cref{fig:annulus-visualization-single,fig:annulus-visualization-mean}.
\end{remark}

\subsection{Outline}

In \cref{sec:background} we give background material.
In \cref{sec:annotation} we show how representative cycles may be used to decompose feature maps on persistence diagrams.
If the feature map is additive, such as persistence images~\cite{Adams:2017}, then we will see that this trivial. 
However, for the persistence landscape, which is not additive, this is nontrivial.
We define the \emph{labeled persistence landscape}, which may be of independent interest.
In \cref{sec:localization} we show how to use the decomposition of a feature map to distribute it to simplices and to pixels or voxels. 
In \cref{sec:randomization} we discuss random feature maps. 
In \cref{sec:stabilization} we give our main theoretical results.
In Section~\ref{sec:examples} we apply our methods to a number of computational examples.

\subsection{Related work}

There is a considerable literature on computing representative cycles.
See for example,
\cite{MR4141870,10.3389/frai.2021.681117,MR4655024,MR4672572,MR4842013,Giunti_2026}, and the references therein. 
While the focus of these papers is somewhat different than ours, we expect that our methods could be profitably combined with these results. 

There has also been much previous work on interpreting the results of TDA. For example,
\cite{math9151723,EDWARDS2021100367,Hartsock:2025,Mishra_2025,toproitopologyinformednetworkapproach},
apply TDA to biological data with an emphasis on interpretation.

\section{Background}
\label{sec:background}

We start with background material from topological data analysis.

\subsection{Ordered simplicial complexes}

Let $K$ be a finite (abstract) simplicial complex. 
That is, $K$ is a collection of nonempty subsets of a finite set $V$ such that for all $v \in V$, $\{v\} \in K$ and if $\sigma \in K$ and $\tau \subseteq \sigma$ and $\tau \neq \emptyset$ then $\tau \in K$.
Then $K$ has a partial order $\subseteq$ given by inclusion.
Let $K_i$ denote the set of $i$-simplices in $K$, i.e., the subset of elements of $K$ of cardinality $i+1$.
If $\sigma \subseteq \tau$ and $\abs{\sigma} = \abs{\tau}-1$, then $\sigma$ is called a \emph{facet} of $\tau$.

Assume that $(K,\subseteq)$ has a linear extension $\leq$, i.e., a total order such that $\sigma \subseteq \tau$ implies $\sigma \leq \tau$.
Call $(K,\leq)$ an \emph{ordered simplicial complex}.

\subsection{Simplicial chains and homology}

Unless specified otherwise, our vector spaces will be over the field $\Z / 2\Z$.
Let $(K,\subseteq)$ be a finite simplicial complex.
Let $S(K)$ be the vector space with basis the set $K$. 
If $\abs{K} = n$ then $\abs{S(K)} = 2^n$.
Elements of $S(K)$ are called \emph{simplicial chains} or just \emph{chains}, and $S(K)$ is called the \emph{simplicial chain vector space}.
An advantage of using the field $\Z / 2\Z$ is that each simplicial chain corresponds to a subset of simplices.
For $\gamma \in S(K)$, the norm $\norm{\gamma}_1$ equals the number of simplices in the simplicial chain.

The \emph{boundary} of a simplex is the simplicial chain given by the sum of its facets.
This map extends by linearity to a linear operator $\partial$ on $S(K)$, called the \emph{boundary operator}.
An element of the kernel of $\partial$ is called a \emph{cycle} and an element of the image of $\partial$ is called \emph{boundary}.
It is a crucial elementary exercise that $\partial \circ \partial = 0$.
It follows that the image of $\partial$ is contained in the kernel of $\partial$.
The \emph{homology vector space} of $K$ is defined to be quotient $H(K) = \ker(\partial)/\im(\partial)$.

We may also compute homology one \emph{degree} at a time as follows.
For $j \in \Z$, let $S_j(K)$ be the vector space with basis $K_j$, where $S_j(K) = 0$ if $K_j = \emptyset$.
Elements of $S_j(K)$ are called simplicial chains in degree $j$.
Then $\partial$ restricts to a linear operator $\partial_j: S_j(K) \to S_{j-1}(K)$.
Let $H_j(K) = \ker(\partial_j) / \im(\partial_{j+1})$, called the homology vector space of $K$ in degree $j$.
Then $S_j(K)$ and $H_j(K)$ are zero if $j<0$ or if  $j$ is greater than or equal to $\max\{\abs{\sigma} \ | \ \sigma \in K\}$.
We remark that 
$S(K) = \bigoplus_{j \in \Z} S_j(K)$ and
$H(K) = \bigoplus_{j \in \Z} H_j(K)$.

Note that $S(K)$ has a canonical basis but not a canonical ordered basis.
However, if we are given an ordered simplicial complex $(K,\leq)$ then $S(K)$ has a canonical ordered basis.

\subsection{Weighted simplicial complexes and persistent homology}

Let $(K,\subseteq)$ be a finite simplicial complex.
A \emph{weight} (also called a \emph{filtration value}) is a order-preserving map $w:(K,\subseteq) \to (\R,\leq)$,
i.e. $\sigma \subseteq \tau$ implies $w(\sigma) \leq w(\tau)$.
Call $(K,w)$ a \emph{weighted simplicial complex}.

For each $a \in \R$, let $K_a = \{\sigma \in K \ | \ w(\sigma) \leq a\}$.
Then with the partial order given by inclusion, $K_a$ is a simplicial complex.
The family of simplicial complexes $\{K_a\}_{a \in \R}$ together the inclusion maps $i_{ab}: K_a \hookrightarrow K_b$ for $a \leq b$ is called a \emph{filtered simplicial complex}. 
It is a crucial elementary fact that for all $a \leq b$, the inclusion maps $i_{ab}: K_a \hookrightarrow K_b$ induce linear maps $H(i_{ab}): H(K_a) \to H(K_b)$ such that $H(i_{aa})$ is the identity map and $H(i_{bc}) \circ H(i_{ab}) = H(i_{ac})$.
The collection of vector spaces $\{H(K_a)\}_{a \in \R}$ together with the linear maps $H(i_{ab})$ for $a \leq b$ is called a \emph{persistence module}, denoted $H(K_\bullet)$, and the images of the maps $H(i_{ab})$ are called \emph{persistent homology vector spaces}.

\subsection{Computational persistent homology}

Let $(K,\leq)$ be a ordered simplicial complex with cardinality $n$.
We now give a computational description of its \emph{persistent homology}~\cite{elz:tPaS,csem:vineyards,edelsbrunnerHarer:survey}.
Recall that the vector space $S(K)$ has a canonical ordered basis. 
Using this basis, we can represent the boundary operator $\delta$ by a $n \times n$ matrix $D$, called the \emph{boundary matrix}.
The rows and columns of $D$ are indexed by the elements of $K$ in the order given by $\leq$ and $D_{\sigma \tau}$ equals $1$ if and only if 
$\sigma$ is a facet of $\tau$.

A \emph{persistence algorithm} is an algorithm that reduces the boundary matrix $D$ using left-to-right column additions to a matrix $R$ that is \emph{reduced} as follows.
For any pair of nonzero columns $R_{\cdot \tau}$ and $R_{\cdot \tau'}$ 
the largest $\sigma$ such that $R_{\sigma \tau}$ is nonzero
does not equal
the largest $\sigma'$ such that $R_{\sigma' \tau'}$ is nonzero.
Equivalently, a persistence algorithm gives, for each boundary matrix $D$ a decomposition $R=DV$ where $R$ is a reduced matrix and $V$ is a upper triangular matrix with ones on the diagonal.
Note that the nonzero entries of $V$ above the diagonal correspond on left-to-right column additions.

Given a reduced matrix $R$, consider a pair $(\sigma,\tau)$ 
such that $R_{\sigma \tau} = 1$ and $\sigma$ is the largest $\rho$ such that $R_{\rho \tau}$ is nonzero.
In this case, we call $\sigma$ a \emph{birth simplex}, $\tau$ a \emph{death simplex}, and $(\sigma,\tau)$ a \emph{birth-death pair}.
We call the set of all such birth-death pairs a \emph{birth-death matching}.
This birth-death matching depends only on the boundary matrix $D$ and is independent of the choice of persistence algorithm~\cite[Pairing Uniqueness Lemma]{csem:vineyards}.

Given an $R=DV$ decomposition and a birth-death pair $(\sigma,\tau)$, there are two
corresponding simplicial chains:
$\beta$ given by the sum of the simplices $\kappa$ such that $V_{\kappa \tau}$ is nonzero;
and $\gamma$ given by the sum of the simplices $\rho$ such that $R_{\rho \tau}$ is nonzero~\cite{edelsbrunnerHarer:survey}.
These chains satisfy $\partial(\beta) = \gamma$
and hence $\partial(\gamma) = 0$.
We call $\gamma$ a \emph{representative cycle} for $(\sigma, \tau)$ and $\beta$ a \emph{bounding chain} for $(\sigma,\tau)$.
Furthermore, the collection all such representative cycles and their images under the inclusion maps $i_{ab}$ where $a \leq b$ gives a minimal set of representatives for bases of the homology vector spaces $H(K_b)$ of the persistence module $H(K_\bullet)$~\cite{edelsbrunnerHarer:survey}. 
Unlike the birth-death matching, the representative cycles and bounding chains depend on the choice of persistence algorithm. 

\begin{figure}[!htb]
    \centering
    \begin{minipage}{.3\textwidth}
        \centering
        \includegraphics[height=25mm]
        {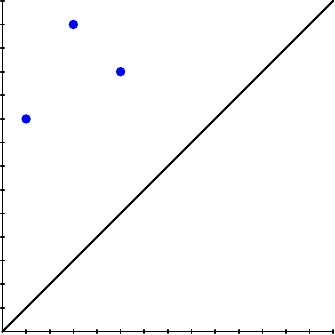}
    \end{minipage}%
    \quad \quad
    \begin{minipage}{0.5\textwidth}
        \centering
        \includegraphics[height=20mm]
        {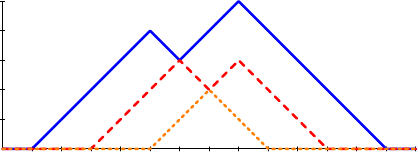}
    \end{minipage}
    \caption{Left: a persistence diagram consisting of three birth-death pairs. Right: the corresponding persistence landscape.}
    \label{fig:pd-pl}
\end{figure}

Consider an order-preserving map $w:(K,\leq) \to (\R,\leq)$,
i.e. $\sigma \subseteq \tau$ implies $w(\sigma) \leq w(\tau)$.
Call such a map a \emph{weight} on $(K,\leq)$ and call $(K,\leq,w)$ a \emph{weighted ordered simplicial complex}.
Now consider the formal sum $\alpha = \sum (w(\sigma),w(\tau))$, where the sum is taken over all birth-death pairs $(\sigma,\tau)$ in the birth-death matching of $(K,\leq)$. 
This formal sum is independent of the choice of the linear extension $\leq$ of $\subseteq$~\cite{csem:vineyards}.
It is called the \emph{persistence diagram} of $(K,w)$.
See the left side of \cref{fig:pd-pl}.
Just as with homology, we may also compute persistence diagrams one degree at a time.

\begin{remark}
    A formal sum is also called a multiset. A \emph{formal sum} on a set $S$ is a nonnegative function $f:S \to \Z$ which takes a nonzero value for only finitely many elements of $S$.
    Thus, a formal sum is element of the free commutative monoid on $S$, which consists of such functions together with elementwise addition and the zero map.
\end{remark}

\subsection{Feature maps} \label{sec:feature-maps}

A \emph{feature map} is a function from a set to a separable real Banach space.
We will assume throughout that $\mathcal{B}$ is a separable real Banach space with norm $\norm{-}$.
Let $\mathcal{D}$ denote the set of all persistence diagrams.
Then a function $F: \mathcal{D} \to \mathcal{B}$ is a feature map.
If $S$ is a finite set, then $\mathcal{B}^S$ is a separable real Banach space with the norm of $x = (x_s)_{s \in S} \in \mathcal{B}^S$ given by $\norm{x} = \sum_{s \in X} \norm{x_s}$.

\begin{example} \label{ex:total-persistence}
    A simple example of a feature map is given the \emph{total persistence} of the persistence diagram. 
    It is the function $p:\mathcal{D} \to \R$ given by $p(\sum_{i=1}^n(b_i,d_i)) = \sum_{i=1}^n d_i - b_i$.
    This map is also called the \emph{total lifetime}.
\end{example}

\begin{example} \label{ex:persistence-landscape}
    A more elaborate example is given by the \emph{persistence landscape} of the the persistence diagram~\cite{bubenik:landscapes}. 
    See the right side of \cref{fig:pd-pl}.
    First, given $b \leq d$, let us define the \emph{triangle function} to be the piecewise linear nonnegative function $\Lambda_{b,d}:\R \to \R$ determined by the points $(b,0)$, $(\frac{b+d}{2},\frac{d-b}{2})$, and $(d,0)$. 
    That is, $\Lambda_{b,d}(t)$ equals $t-b$ if $b \leq t \leq \frac{b+d}{2}$, $d-t$ if $\frac{b+d}{2} < t \leq d$ and $0$ otherwise.
    Second, let $\kmax$ be the function that takes a list of $n$ nonnegative real numbers and if $k \leq n$ it returns the $k$th largest number and otherwise it returns $0$.
    Then,
    the persistence landscape of a persistence diagram $\sum_{i=1}^n(b_i,d_i)$ is the nonnegative function $\lambda: \N \times \R \to \R$ given by
    \begin{equation}
        \lambda(k,t) = \kmax (\Lambda_{b_1,d_1}(t),\ldots,\Lambda_{b_n,d_n}(t)).
    \end{equation}
    Finally, the persistence landscape feature map is the function $\Lambda: \mathcal{D} \to L^p(\N \times \R)$, where $1 \leq p \leq \infty$, given by 
    mapping a persistence diagram to its persistence landscape.
    This feature map is injective~\cite{bubenik:landscapes}. 
    In fact, viewed as a kernel it is characteristic for certain generic empirical measures~\cite{Bubenik:2020b}.
\end{example}

\begin{remark}
    The definition of the persistence landscape may seem somewhat ad-hoc, but it can be defined more intrinsically using persistence modules~\cite{bubenik:landscapes}, and it may also be defined more directly using graded persistence diagrams~\cite{MR4356248}.
\end{remark}

\begin{example} \label{ex:death-vector}
    For persistence diagrams of the form $\alpha = \sum_{i=1}^n (0,d_i)$, we have the \emph{death vector} given by the sequence $(d_{(1)},d_{(2)},d_{(3)},\ldots)$, where $d_{(k)} = \kmax(d_1,d_2,\ldots,d_n)$,
    and where $\kmax$ is defined in \cref{ex:persistence-landscape}.
    Let $\mathcal{D}_0$ denote the subset of such persistence diagrams. 
    Then the death vector is a map $\delta: \mathcal{D}_0 \to \ell^p$, where $1 \leq p \leq \infty$.
\end{example}

\subsection{Kernels}

A \emph{kernel} on $\R^d$ is nonnegative, integrable function $\kappa: \R^d \to \R$ such that $\int_{\R^d} \kappa(x) dx = 1$,
$\int_{\R^d} x \kappa(x) dx = 0$, and
$\int_{\R^d} x^2 \kappa(x) dx < \infty$.

\begin{example} \label{ex:gaussian-kernel}
    Let $a > 0$.
    The (isotropic) \emph{Gaussian kernel} on $\R^d$ is given by
    \begin{equation} \label{eq:gaussian-kernel}
        \kappa_a(x) = \frac{1}{a^d (2\pi)^{d/2}} \exp\left(-\frac{1}{2a^2} \norm{x}^2\right).
    \end{equation}
\end{example}

Let $V_d$ denote the volume of the unit ball in $\R^d$.

\begin{example} \label{ex:triangular-kernel}
    Let $a > 0$.
    The \emph{triangular kernel} on $\R^d$ is given by
    \begin{equation} \label{eq:triangular-kernel}
        \kappa_a(x) = \frac{d+1}{a^d V_d} \left(1 - \frac{\norm{x}}{a}\right)_+,
    \end{equation}
    where $x_+$ equals $x$ if $x \geq 0$ and otherwise equals $0$.
\end{example}

\begin{example} \label{ex:epanchknikov-kernel}
    Let $a > 0$.
    The \emph{Epanechnikov kernel} on $\R^d$ is given by
    \begin{equation} \label{eq:epanechnikov-kernel}
        \kappa_a(x) = \frac{d+2}{2a^d V_d} \left(1 - \frac{\norm{x}^2}{a^2}\right)_+.
    \end{equation}
\end{example}

\section{Using representative cycles}
\label{sec:annotation}

In this section we discuss representative cycles of persistent homology and how they may be incorporated in topological data analysis.

\subsection{The four simplicial chains of a birth-death pair}

Let $(K,w)$ be a finite weighted simplicial complex.
Let $\leq$ be a linear extension of $\subseteq$ such that $w:{(K,\leq)} \to (\R,\leq)$ is an order-preserving function.
Then $(K,\leq,w)$ is a weighted ordered simplicial complex.

Let $\mathcal{A}$ be a choice of persistence algorithm.
Applying $\mathcal{A}$ to $(K,\leq)$, we obtain a birth-death matching, 
and for each birth-death pair we have a corresponding birth simplex, death simplex, representative cycle, and bounding chain.
Each of the latter four is a simplicial chain, i.e. a formal sum of simplices. 
In what follows, we will fix a choice of one of these four simplicial chains. 
For concreteness, we may refer to this choice the representative cycle, but
our results apply equally to all four cases. 


\subsection{Labeling persistence diagrams with representative cycles}
\label{sec:labeled-pd}

Recall that a persistence diagram 
is a formal sum $\alpha = \sum_{i=1}^n (b_i,d_i)$ where $b_i \leq d_i$ and $n \geq 0$.

\begin{figure}[!htb]
    \centering
    \begin{minipage}{.3\textwidth}
        \centering
        \includegraphics[height=25mm]
        {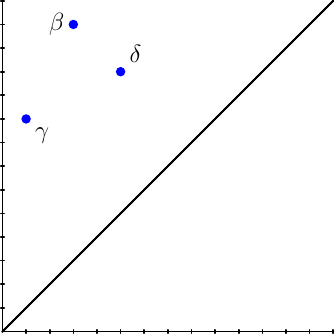}
    \end{minipage}%
    \quad \quad
    \begin{minipage}{0.5\textwidth}
        \centering
        \includegraphics[height=20mm]
        {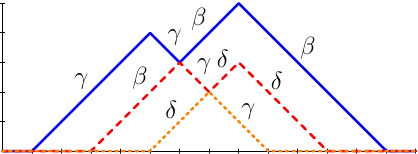}
    \end{minipage}
    \caption{Left: a labeled persistence diagram consisting of three labeled birth-death pairs. Right: the corresponding labeled persistence landscape.}
    \label{fig:lpd-lpl}
\end{figure}

\begin{definition} \label{def:lpd}
    A \emph{labeled persistence diagram} 
    is a formal sum $\hat{\alpha} = \sum_{i=1}^n (b_i,d_i,\ell_i)$, where $b_i \leq d_i$, 
    $\ell_i \in \N$, and $n \geq 0$.
    See the left side of \cref{fig:lpd-lpl}.
    Let $\hat{\mathcal{D}}$ denote the set of labeled persistence diagrams.
    Let $\varphi: \hat{\mathcal{D}} \to \mathcal{D}$ be the map that forgets the labels, i.e. $\varphi(\sum_{i=1}^n (b_i,d_i,\ell_i)) = \sum_{i=1}^n(b_i,d_i)$.
\end{definition}

Let $(K,\leq,w)$ be a weighted ordered finite simplicial complex. 
Let $\mathcal{A}$ be a choice of persistence algorithm.
From this data, we have a birth-death matching $\{(\sigma_i,\tau_i)\}_i$ and
persistence diagram $\alpha = \sum_i(w(\sigma_i),w(\tau_i))$, where the birth-death pair $(\sigma_i,\tau_i)$ has the corresponding representative cycle $\gamma_i$.
The set $S(K)$ is finite and we may choose an injective map $S(K) \to \N$.
Therefore, we have the labeled persistence diagram 
$\hat{\alpha} = \sum_i (w(\sigma_i),w(\tau_i),\gamma_i)$.
Note that the persistence diagram is independent of the choice of $\leq$ and $\mathcal{A}$, but the labeled persistence diagram is not.

\subsection{Decomposing feature maps using representative cycles} \label{sec:decomposing}

Recall that a feature map is a function 
$F: \mathcal{D} \to \mathcal{B}$ 
from the set of persistence diagrams 
to a Banach space $\mathcal{B}$, 
and that $\hat{\mathcal{D}}$ is the set of labeled persistence diagrams.
Given a labeled persistence diagram $\hat{\alpha} = \sum_{i=1}^n (b_i,d_i,\ell_i)$, let $\alpha$ denote the corresponding persistence diagram $\sum_{i=1}^n (b_i,d_i)$ that forgets the labels.
Given a map $G: \hat{\mathcal{D}} \to \mathcal{B}^{\N}$, and $\ell \in \N$,
let $G_{\ell}(\hat{\alpha}) = G(\hat{\alpha})(\ell) \in \mathcal{B}$.
Then the function $G(\hat{\alpha})$ may be alternatively viewed as the sequence $G_1(\hat{\alpha}),G_2(\hat{\alpha}),\ldots$.

\begin{definition} \label{def:labeled-decomposition}
    A \emph{labeled decomposition} of a feature map $F: \mathcal{D} \to \mathcal{B}$ 
    is a function $\hat{F}: \hat{\mathcal{D}} \to \mathcal{B}^{\N}$ such that 
    for all $\hat{\alpha} \in \hat{\mathcal{D}}$, the sequence
    $\hat{F}_1(\hat{\alpha}),\hat{F}_2(\hat{\alpha}),\ldots$ is eventually zero and
    $\sum_\ell \hat{F}_\ell(\hat{\alpha}) = F(\alpha)$.
    That is, we have the following commutative diagram,
\begin{equation} \label{eq:labeled-decomposition}
    \begin{tikzcd}
        \mathcal{D} \ar[r,"F"] & \mathcal{B}\\
        \hat{\mathcal{D}} \ar[r,"\hat{F}"] \ar[u,"\varphi"] & \mathcal{B}^\N \ar[u,"\Sigma_\N"']
    \end{tikzcd}
\end{equation}
    where 
    $\Sigma_\N$ is the function given by the (finite) sum over $\N$.
    We call $\hat{F}_\ell(\hat{\alpha})$ the part of $F(\alpha)$ labeled by $\ell$.
\end{definition}

\begin{example}[Total persistence] \label{ex:labled-total-persistence}
    Consider the feature map $p: \mathcal{D} \to \R$ given by total persistence (\cref{ex:total-persistence}).
    Then the labeled decomposition of total persistence is given by the sequence $\hat{p}_1(\hat{\alpha}),\hat{p}_2(\hat{\alpha}),\ldots$ where $\hat{p}_\ell(\sum_{i=1}^n(b_i,d_i,\ell_i))$ is the sum of all $d_i-b_i$ for which $\ell_i=\ell$.
    We have that $\sum_\ell \hat{p}_\ell(\hat{\alpha}) = p(\alpha)$.
\end{example}

\begin{example}[Linear feature maps] 
\label{ex:linear-feature-map}
Consider a feature map $p:\mathcal{D} \to \R$ such that $p(\sum_{i=1}^n (b_i,d_i)) = \sum_{i=1}^n p((b_i,d_i))$.
That is, $p$ is additive. 
For example, total persistence (\cref{ex:total-persistence}) is such a linear feature map.
Other examples include 
the persistence scale-space feature map~\cite{Reininghaus:2015}, 
the persistence-weighted Gaussian feature map~\cite{Kusano:2016}, 
the persistence surface~\cite{Adams:2017}, and
persistence images~\cite{Adams:2017}.
 
Any such feature map has a labeled decomposition given by $\hat{p}_{\ell}(\sum_{i=1}^n(b_i,d_i,\ell_i)) = \sum p((b_i,d_i))$, where the sum is taken over all $i$ such that $\ell_i=\ell$.
\end{example}

\begin{example}[Persistence landscape] \label{ex:label-persistence-landscape}
    Consider the persistence landscape feature map $\Lambda: \mathcal{D} \to {L^p(\N \times \R)}$, where $1 \leq p \leq \infty$ (\cref{ex:persistence-landscape}).
    First, we consider elements $(t,i) \in \R \times \N$ which consist of a value $t \in \R$ and a label $i \in \N$.
    Given $b \leq d$ and $\ell \in \N$, 
    define the \emph{labeled triangle function} to be the function $\Lambda_{b,d,\ell}:\R \to \R \times \N$ given by $\Lambda_{b,d,\ell}(t) = (\Lambda_{b,d}(t),\ell)$.
    If we let $\pi_1: \R \times \N \to \R$ and $\pi_2: \R \times \N \to \N$ be the coordinate projection maps,
    then $\pi_1(\Lambda_{b,d,\ell}(t)) = \Lambda_{b,d}(t)$ is the value and
    $\pi_2(\Lambda_{b,d,\ell}(t)) \in \N$ is the label.
    
    Second, let $\R \times \N$ have the total order given by lexicographic order, i.e. $(s,i) \leq (t,j)$ if and only if either $s < t$ or $s=t$ and $i \leq j$.
    Then, let $\kmax$ be the function that takes a list of $n$ elements of $\R \times \N$ whose values are nonnegative, and if $k \leq n$ it returns the $k$th largest element and otherwise it returns $(0,0)$.

    Third,
    the \emph{labeled persistence landscape} of a labeled persistence diagram, $\hat{\alpha} = \sum_{i=1}^n(b_i,d_i,\ell_i)$, is the function $\hat{\lambda}: \N \times \R \to \R \times (\N \cup \{0\})$ given by
    \begin{equation}
        \hat{\lambda}(k,t) = \kmax (\Lambda_{b_1,d_1,\ell_1}(t),\ldots,\Lambda_{b_n,d_n,\ell_n}(t)).
    \end{equation}
    See the right side of \cref{fig:lpd-lpl}.
    Note that $\pi_1(\hat{\lambda}(k,t)) = \lambda(k,t)$, where $\lambda$ is the persistence landscape of $\alpha = \sum_{i=1}^n(b_i,d_i)$.
    
    Finally, for $\hat{\alpha} \in \hat{\mathcal{D}}$, let $\hat{\lambda}$ be its labeled persistence landscape and let $\lambda = \pi_1 \hat{\lambda}$ be its persistence landscape.
    Then the labeled decomposition of the persistence landscape feature map
    is given by the sequence $\hat{\Lambda}_1(\hat{\alpha}), \hat{\Lambda}_2(\hat{\alpha}),\ldots$,
    where 
    \begin{equation*}
        \hat{\Lambda}_\ell(\hat{\alpha})(k,t) = \begin{cases}
            \lambda(k,t) & \text{if } \pi_2 \hat{\lambda}(k,t) = \ell\\
            0 & \text{otherwise}.
        \end{cases}
    \end{equation*}
    For $\ell \in \N$, $\hat{\Lambda}_\ell(\hat{\alpha})$ is the part of the persistence landscape of $\alpha$ labeled by $\ell$.
    We have that $\sum_\ell \hat{\Lambda}_\ell(\hat{\alpha}) = \Lambda(\alpha)$.
\end{example}

\begin{example}[Death vector] 
\label{ex:labeled-death-vector}
    Recall the death vector feature map $\delta: \mathcal{D}_0 \to \ell^p$, where $1 \leq p \leq \infty$ (\cref{ex:death-vector}).
    Consider a labeled persistence diagram $\hat{\alpha} = \sum_{i=1}^n (0,d_i,\ell_i)$.
    Let $(d\ell_{(1)},d\ell_{(2)},d\ell_{(3)},\ldots)$ be the sequence in $\R \times \N$ given by $d\ell_{(k)} = \kmax ((d_1,\ell_1),\ldots,(d_n,\ell_n))$, where $\kmax$ is defined in \cref{ex:label-persistence-landscape}.
    Let $\pi_1:\R \times \N \to \R$ and $\pi_2:\R \times \N \to \N$ denote the coordinate projection maps. 
    For $k \geq 1$, let $d_{(k)} = \pi_1(d\ell_{(k)})$ and let $\ell_{(k)} = \pi_2(d\ell_{(k)})$.
    Then the \emph{labeled death vector} of $\hat{\alpha}$ is given by 
    the sequences $\hat{d}_\ell$, where
    $(\hat{d}_\ell)_i$ equals $d_{(i)}$ if $\ell_{(i)} = \ell$ and equals $0$ otherwise.
    This defines the labeled death vector feature map $\hat{\delta}: \hat{\mathcal{D}}_0 \to (\ell^p)^\N$.
\end{example}

\subsection{Feature maps and linear representations}

Assume that we have a linear representation $g:\mathcal{B} \to \R$.
Then \eqref{eq:labeled-decomposition} extends to the following commutative diagram,
\begin{equation*} 
    \begin{tikzcd}
        \mathcal{D} \ar[r,"F"] & \mathcal{B} \ar[r,"g"] & \R \\
        \hat{\mathcal{D}} \ar[r,"\hat{F}"] \ar[u,"\varphi"] & \mathcal{B}^\N \ar[u,"\Sigma_\N"'] \ar[r,"g^\N"] & \R^\N \ar[u,"\Sigma_\N"']
    \end{tikzcd}
\end{equation*}
where $g^\N$ is the map given by applying $g$ to each component.
Therefore, $g \circ F: \mathcal{D} \to \R$ is a feature map, and $g^\N \circ \hat{F}$ is a labeled decomposition of $g \circ F$. 
An important special case is the following.

\begin{example} \label{ex:svm}
    Assume that we have a feature map $F: \mathcal{D} \to \mathcal{H}$, where $\mathcal{H}$ is a Hilbert space, and that we have a labeled decomposition $\hat{F}: \hat{\mathcal{D}} \to \mathcal{H}$ of $F$.
    Now assume that we have training data for a classification or regression problem consisting of elements $\mathcal{H}$ obtained from applying the feature map to persistence diagrams, together with corresponding class or scalar values.
    We apply support vector machines to learn a normal vector $a \in \mathcal{H}$.
    Let $g: \mathcal{B} \to \R$ be given by the mapping $v \mapsto \langle v, a \rangle$.

    For testing data, we have scalar estimates $\langle F(\alpha), a \rangle$. 
    For classification, we obtain the estimated class by seeing whether this scalar is at least a number $b$, also learned by SVM, and for 
    regression, we obtained the estimated scalar by subtracting $b$.
    In either case we would like to know what parts of the simplicial complex $K$ are responsible for the value $\langle F(\alpha), a \rangle$.
    
    Since the inner product is linear, we have 
    \begin{equation*} 
        \langle F(\alpha),a \rangle = \sum_{\ell} \langle \hat{F}_\ell(\hat{\alpha}), a \rangle.
    \end{equation*}
    Thus, using 
    this equation, we are able to determine the contribution of each of the simplicial chains to the estimator.
\end{example}

In the next section, we will furthermore distribute this contribution to simplices, and then use the positions of these simplices to produce a visualization on the ambient Euclidean space.

\section{Distributing feature maps to simplices and pixels/voxels}
\label{sec:localization}

In the previous section we decomposed feature maps according to the labels of the persistence diagram. 
In this section we show how to distribute feature maps to the simplices.
If the simplicial complex has an image in a Euclidean space, we furthermore distribute the feature maps to a choice of pixels/voxels of this space.

\subsection{Distribution to simplices}

Let $F: \mathcal{D} \to \mathcal{B}$ be a feature map with a labeled decomposition $\hat{F}: \mathcal{D} \to \mathcal{B}^\N$.
Recall that $\hat{F}$ may be represented as a sequence of functions $(\hat{F}_\ell:\hat{\mathcal{D}} \to \mathcal{B})_{\ell \in \N}$.

Let $K$ be a finite simplicial complex. 
In this section, we will restrict $\hat{\mathcal{D}}$ to the subset $\hat{\mathcal{D}}_K$ of labeled persistence diagrams arising from $K$. 
So, all labels will be elements of $S(K)$.
In particular, for each such labeled persistence diagram $\hat{\alpha}$, we have $\sum_{\gamma \in S(K)} \hat{F}_\gamma(\hat{\alpha}) = F(\alpha)$.
We will define a function $\bar{F}: \hat{\mathcal{D}}_K \to \mathcal{B}^K$, which is given by a sequence of functions $(\bar{F}_\sigma: \hat{\mathcal{D}}_K \to \mathcal{B})_{\sigma \in K}$.

\begin{definition} \label{def:f-bar}
    Let $\hat{f}:S(K) \to \mathcal{B}$, where $\mathcal{B}$ is a Banach space. 
    Define $\bar{f}:K \to \mathcal{B}$ by
    \begin{equation*}
        \bar{f}(\sigma) = \sum_{\gamma \ni \sigma} \frac{\hat{f}(\gamma)}{\norm{\gamma}_1},
    \end{equation*}
    where sum is taken over all simplicial chains $\gamma$ that contain the simplex $\sigma$.
\end{definition}

\begin{lemma} \label{lem:f-bar}
    $\sum_{\sigma \in K} \bar{f}(\sigma) = \sum_{\gamma \in S(K)} \hat{f}(\gamma)$.
\end{lemma}

\begin{proof}
    $\sum_{\sigma \in K} \bar{f}(\sigma) = \sum_{\sigma \in K} \sum_{\gamma \ni \sigma} \frac{\hat{f}(\gamma)}{\norm{\gamma}_1} = \sum_{\gamma \in S(K)} \sum_{\sigma \in \gamma} \frac{\hat{f}(\gamma)}{\norm{\gamma}_1} = \sum_{\gamma \in S(K)} \hat{f}(\gamma)$.
\end{proof}

As a special case of \cref{def:f-bar}, we have the following.

\begin{definition} \label{def:simplex-feature-map}
    The \emph{simplex distribution} of $F$ is the function $\bar{F}: \hat{\mathcal{D}}_K \to \mathcal{B}^K$ given by 
    \begin{equation*}
        \bar{F}_\sigma(\hat{\alpha}) = \sum_{\gamma \ni \sigma} \frac{\hat{F}_\gamma(\hat{\alpha})}{\norm{\gamma}_1},
    \end{equation*}
    where the sum is taken over all 
    simplicial chains $\gamma$ that contain the
    simplex $\sigma$.
\end{definition}

Combining \cref{lem:f-bar} and \cref{def:labeled-decomposition}, we have the following.

\begin{corollary} \label{cor:f-bar}
    $\sum_{\sigma \in K} \bar{F}_\sigma(\hat{\alpha}) = F(\alpha)$.
\end{corollary}

\begin{example}[\cref{ex:svm} continued] \label{ex:svm2}
    Using \cref{cor:f-bar}, we have $\langle F(\alpha), a \rangle = \sum_{\sigma \in K} \langle \bar{F}_\sigma(\hat{\alpha}), a \rangle$.
    We can also distribute the bias among the simplices as follows.
    We will assume that $F(\alpha)$ only depends on the persistence diagram in degree $j$. Then,
    \begin{equation*}
        \langle F(\alpha), a) - b = \sum_{\sigma \in K_j} \left( \langle \overline{F}_\sigma(\hat{\alpha}), a \rangle - \frac{b}{\abs{K_j}} \right).
    \end{equation*}
\end{example}

\subsection{Distribution to pixels/voxels}

\begin{figure}[!htb]
    \centering
    \includegraphics[height=20mm]{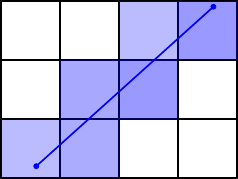}
    \caption{The value of the simplex is distributed to the pixels containing its image.}
    \label{fig:pixel-localization}
\end{figure}

In this section, we assume that our simplicial complex lies in a Euclidean space that has a decomposition into a grid of cubes (i.e. pixels/voxels).
Given values on the simplices of $K$, we distribute them onto the cubes according their intersections of the simplices and the cubes.
See \cref{fig:pixel-localization}.

Let $K$ be a finite simplicial complex.
Assume that we have a function $v: K_0 \to \R^d$ for some $d \geq 1$, which we call \emph{position}.
For each $\sigma \in K$, let $\abs{\sigma}$ denote the convex hull of the finite set $v(\sigma)$.
Let $d_\sigma$ denote the dimension of $\abs{\sigma}$.
For any $k \geq 0$ and any polytope $p$ of dimension at most $k$,
let $\vol_k(p)$ denote its $k$-dimensional volume.
Let $\abs{K} = \bigcup_{\sigma \in K} \abs{\sigma}$.

Let $\eps > 0$.
Consider a partition 
of $\R^d$ into translates of $[0,\eps)^d$.
Since $K$ is finite, $\abs{K}$ intersects only finitely many of these translates.
Let $\mathcal{C}$ be the finite set of these translates that intersect $\abs{K}$.

\begin{definition} \label{def:f-tilde}
    Let $f:K \to \mathcal{B}$, where $\mathcal{B}$ is a Banach space.
    Define $\tilde{f}: \mathcal{C} \to \mathcal{B}$ by
    \begin{equation*}
        \tilde{f}(c) = \sum_{\sigma \in K} \frac{\vol_{d_\sigma}(\abs{\sigma} \cap c)}{\vol_{d_\sigma}(\abs{\sigma})} f(\sigma).
    \end{equation*}
\end{definition}

\begin{lemma} \label{lem:f-tilde}
    $\sum_{c \in \mathcal{C}} \tilde{f}(c) = \sum_{\sigma \in K}f(\sigma)$.
\end{lemma}

\begin{proof}
    $\sum_{c \in \mathcal{C}} \tilde{f}(c)  
    = \sum_{c \in \mathcal{C}} \sum_{\sigma \in K} \frac{\vol_{d_\sigma}(\abs{\sigma} \cap c)}{\vol_{d_\sigma}(\abs{\sigma})} f(\sigma) 
    = \sum_{\sigma \in K} f(\sigma)
    \sum_{c \in \mathcal{C}} 
    \frac{\vol_{d_\sigma}(\abs{\sigma} \cap c)}{\vol_{d_\sigma}(\abs{\sigma})}  
    = \sum_{\sigma \in K} f(\sigma)$.
\end{proof}

As a special case of \cref{def:f-tilde}, we have the following.

\begin{definition} \label{def:spatial-localization}
    The 
    \emph{pixel/voxel distribution} of $F$ is the function $\tilde{F}: \hat{\mathcal{D}}_K \to \mathcal{B}^\mathcal{C}$ given by
    \begin{equation*}
        \tilde{F}_c(\hat{\alpha},v) = \sum_{\sigma \in K} \frac{\vol_{d_\sigma}(\abs{\sigma} \cap c)}{\vol_{d_\sigma}(\abs{\sigma})} \bar{F}_\sigma(\hat{\alpha}).
    \end{equation*}
\end{definition}

Combining \cref{lem:f-tilde} and \cref{cor:f-bar} we have the following.

\begin{corollary} \label{cor:f-tilde}
     $\sum_{c \in \mathcal{C}} \tilde{F}_c(\hat{\alpha},v) = F(\alpha)$.
\end{corollary}


\begin{example}[\cref{ex:svm2} continued] \label{ex:svm3}
    Using \cref{cor:f-tilde}, we have $\langle F(\alpha), a \rangle = \sum_{c \in \mathcal{C}} \tilde{F}_c(\hat{\alpha}), a \rangle$.
    Assuming that $F(\alpha)$ only depends on the persistence diagram in degree $j$, we furthermore distribute the bias term to get
    \begin{equation*}
        \langle F(\alpha),a \rangle = \sum_{c \in \mathcal{C}} \left( \langle \tilde{F}_c(\hat{\alpha}), a \rangle - \frac{b}{\abs{K_j}} \sum_{\sigma \in K_j} \frac{\vol_{d_\sigma}(\abs{\sigma} \cap c)}{\vol_{d_\sigma}(\abs{\sigma})} \right).
    \end{equation*}
\end{example}

\section{Randomizing feature maps}
\label{sec:randomization}

In this section, we consider our feature maps to be functions of the weight on the simplicial complex.
We observe that they are discontinuous with respect to perturbations of the weight.
As a step toward removing this discontinuous behavior,
we introduce random weights, which produce random feature maps.
In Section~\ref{sec:stabilization}, we will obtain stable versions of our feature maps by averaging over random perturbations.

\subsection{Feature maps as functions on a space of weights}

Let $K$ be a finite simplicial complex. 
Recall that a weight on $K$ is an order-preserving function $w:{(K,\subseteq)} \to (\R,\leq)$.
Let $W = \{w \in \R^K \ | \ w \text{ is a weight on } K\}$.
The assignment of $w$ to the persistence diagram of $(K,w)$ defines a map $\alpha: W \to \mathcal{D}$.
For appropriate choices of metrics on $\mathcal{D}$, this map is continuous and piecewise linear~\cite{csem:vineyards,Skraba:2020,Bubenik:2023a}.


Next, we will use labeled persistence diagrams to define a map $\hat{\alpha}: \R^K \to \hat{\mathcal{D}}$.
Recall that the labeled persistence diagram depends on a choice of persistence algorithm and linear extension $\leq$ of $\subseteq$.
For the former, choose a persistence algorithm $\mathcal{A}$.
The latter requires more care, since we need to choose a linear extension $\leq$ such that $w$ is a weight on $(K,\leq)$, i.e. $w:(K,\leq) \to (\R,\leq)$ is order-preserving.
To do this, we follow Bubenik and Bush~\cite{Bubenik:2023a}.
First, choose a linear extension $\trianglelefteq$ of $\subseteq$.
Second, define $\sigma \leq_{(w,\trianglelefteq)} \tau$ if $w(\sigma) < w(\tau)$ or $w(\sigma) = w(\tau)$ and $\sigma \trianglelefteq \tau$.
Then $\leq_{(w,\trianglelefteq)}$ is a linear extension of $\subseteq$ and $w$ is a weight on $(K,\leq_{(w,\trianglelefteq)})$~\cite{Bubenik:2023a}.
Now we define $\hat{\alpha}(w)$ to be the labeled persistence diagram of $(K,\leq_{(w,\trianglelefteq)},w)$, where we use the persistence algorithm $\mathcal{A}$.
For appropriate choices of metric on $\hat{\mathcal{D}}$, this function is discontinuous but piecewise linear~\cite{Bubenik:2023a}.
Recall that $\varphi: \hat{\mathcal{D}} \to \mathcal{D}$ is the map the forgets the labels.


\begin{lemma} \label{lem:weight}
    For every weight $w$ on $K$, $\varphi(\hat{\alpha}(w)) = \alpha(w)$.
\end{lemma}

\begin{proof}
    We have that $\hat{\alpha}(w)$ is the labeled persistence diagram of $(K,\leq_{w,\trianglelefteq)},w)$ using the persistence algorithm $\mathcal{A}$.
    Also, $\alpha(w)$ is the persistence diagram of $(K,w)$, which is independent of the choice of linear extension of $\subseteq$ and the choice of persistence algorithm.
    Therefore $\varphi(\hat{\alpha}(w)) = \alpha(w)$.
\end{proof}

Let $F: \mathcal{D} \to \mathcal{B}$ be a feature map with labeled decomposition $\hat{F}: \hat{\mathcal{D}} \to \mathcal{B}^\N$.
Recall that $\hat{\mathcal{D}}_K$ consists of those labeled persistence diagrams arising from $K$.
Combining \cref{lem:weight} with \eqref{eq:labeled-decomposition}, we have the following commutative diagram.
\begin{equation*}
    \begin{tikzcd}
        W \ar[r,"\alpha"] & \mathcal{D} \ar[r,"F"] & \mathcal{B}\\
        W \ar[u,"\Id"] \ar[r,"\hat{\alpha}"] & \hat{\mathcal{D}}_K \ar[u,"\varphi"] \ar[r,"\hat{F}"] & \mathcal{B}^{S(K)} \ar[u,"\Sigma_{S(K)}"']
    \end{tikzcd}
\end{equation*}

Next, we assume that we have a parametrized family of weights.
That is, we have a set $A \subseteq \R^m$ for some $m \geq 1$, and a map $\psi:A \to W$.

\begin{example}
    For example, if $K$ has $M$ vertices then $A \subset (\R^d)^M$ may be a set of allowed positions, and $\psi$ may assign a set of vertices the radius of their smallest enclosing ball.
\end{example}

We define $\alpha': \R^m \to \mathcal{D}$ by defining $\alpha'(a) = \alpha(\psi(a))$ for $a \in A$, and defining $\alpha'(b) = 0$ for $b \in \R^m \setminus A$.
We similarly define $\hat{\alpha}': \R^m \to \hat{\mathcal{D}}_K$.
For convenience, we henceforth denote $\alpha'$ and $\hat{\alpha}'$ by $\alpha$ and $\hat{\alpha}$.
Hence, we have the following commutative diagram.
\begin{equation*} 
    \begin{tikzcd}
        \R^m \ar[r,"\alpha"] & \mathcal{D} \ar[r,"F"] & \mathcal{B}\\
        \R^m \ar[u,"\Id"] \ar[r,"\hat{\alpha}"] & \hat{\mathcal{D}}_K \ar[u,"\varphi"] \ar[r,"\hat{F}"] & \mathcal{B}^{S(K)} \ar[u,"\Sigma_{S(K)}"']
    \end{tikzcd}
\end{equation*}
Elements of $\R^m$ are parameters for weights, but, for brevity, we may refer to them as weights.
Let $\eta_{S(K)}: \R^m \to \mathcal{B}^{S(K)}$ be the map given by the composition $\hat{F} \circ \hat{\alpha}$.
We call $\eta_{S(K)}$ the \emph{chain feature map}.

Let $\bar{F}: \hat{\mathcal{D}}_K \to \mathcal{B}^K$ be the induced distribution to simplices (\cref{def:simplex-feature-map}).
Combining \cref{lem:weight} with \cref{cor:f-bar}, 
we have the following commutative diagram.
\begin{equation} \label{cd:heatmap-simplices}
    \begin{tikzcd}
        \R^m \ar[r,"\alpha"] & \mathcal{D} \ar[r,"F"] & \mathcal{B}\\
        \R^m \ar[u,"\Id"] \ar[r,"\hat{\alpha}"] & \hat{\mathcal{D}}_K \ar[u,"\varphi"] \ar[r,"\bar{F}"] & \mathcal{B}^K \ar[u,"\Sigma_K"']
    \end{tikzcd}
\end{equation}
Let $\eta_K: \R^m \to \mathcal{B}^K$ be the map given by the composition $\bar{F} \circ \hat{\alpha}$.
We call $\eta_K$ the \emph{simplex feature map}.

Now assume that $K$ has a position function $v:K_0 \to \R^d$.
Let $\tilde{F}: \hat{\mathcal{D}}_K \to \mathcal{B}^\mathcal{C}$ be the induced distribution to pixels/voxels (\cref{def:spatial-localization}).
Combining \cref{lem:weight} with \cref{cor:f-tilde}, 
we have the following commutative diagram.
\begin{equation} \label{cd:heatmap-pixels}
    \begin{tikzcd}
        \R^m \ar[r,"\alpha"] & \mathcal{D} \ar[r,"F"] & \mathcal{B}\\
        \R^m \ar[u,"\Id"] \ar[r,"\hat{\alpha}"] & \hat{\mathcal{D}}_K \ar[u,"\varphi"] \ar[r,"\tilde{F}"] & \mathcal{B}^{\mathcal{C}} \ar[u,"\Sigma_{\mathcal{C}}"']
    \end{tikzcd}
\end{equation}
Let $\eta_{\mathcal{C}}: \R^m \to \mathcal{B}^\mathcal{C}$ be the map given by the composition $\tilde{F} \circ \hat{\alpha}$.
We call $\eta_{\mathcal{C}}$ the 
\emph{pixel/voxel feature map}.

The chain feature map, simplex feature map, and pixel/voxel feature map are discontinuous, but are piecewise continuous and hence Borel measurable.

\subsection{Random chain, simplex, and pixel/voxel feature maps}

Let $K$ be a finite simplicial complex with position function $v: K_0 \to \R^d$.
Let $F: \mathcal{D} \to \mathcal{B}$ with a labeled decomposition $\hat{F}: \hat{\mathcal{D}}_K \to \mathcal{B}^{S(K)}$.
Choose a persistence algorithm $\mathcal{A}$ and a linear extension $\trianglelefteq$ of the partial order $\subseteq$ on $K$.
From the previous section, we have 
the chain feature map $\eta_{S(K)}: \R^m \to \mathcal{B}^{S(K)}$,
the simplex feature map $\eta_K: \R^m \to \mathcal{B}^K$, and 
the pixel/voxel feature map $\eta_\mathcal{C} : \R^m \to \mathcal{B}^\mathcal{C}$.

Let $(\Omega,\mathcal{F},P)$ be a probability space.
Let $X:(\Omega,\mathcal{F},P) \to \R^m$ be a Borel random variable, i.e. a random weight.
Then 
$\eta_{S(K)} \circ X: (\Omega,\mathcal{F},P) \to \mathcal{B}^{S(K)}$ is a \emph{random chain feature map},
$\eta_K \circ X: (\Omega,\mathcal{F},P) \to \mathcal{B}^K$ is a \emph{random simplex feature map}, and
$\eta_\mathcal{C} \circ X: (\Omega,\mathcal{F},P) \to \mathcal{B}^\mathcal{C}$ is a \emph{random pixel/voxel feature map}.

\begin{example} \label{ex:gaussian-noise}
    Let $w \in \R^m$ be a (parameter for a) weight on $K$.
    Let $\eps$ be a $m$-dimensional Gaussian random variable with mean $0$ and $m \times m$ covariance (positive definite) matrix $\Sigma$.
    Then $X_w = w + \eps$ is a $m$-dimensional Gaussian random variable with mean $w$ and covariance matrix $\Sigma$.
    We have the random chain feature map $\eta_{S(K)} \circ X_w$, the random simplex feature map $\eta_K \circ X_w$, and the random pixel/voxel feature map $\eta_\mathcal{C} \circ X_w$. 
\end{example}

\subsection{Convergence of our random feature maps}

From the previous section, we have 
a random chain feature map, 
a random simplex feature map, and
a random pixel/voxel feature map.
Let $V: (\Omega,\mathcal{F},P) \to \mathcal{B}^S$ denote any of these three random variables.
We apply well known results from probability in Banach spaces; see~\cite{ledoux-talagrand:book} for more details.

The composition $\norm{V}: (\Omega,\mathcal{F},P) \xrightarrow{V} \mathcal{B}^S \xrightarrow{\norm{-}} \R$ is a Borel random variable.
If $E[\norm{V}] < \infty$ then $V$ has a Pettis integral $E[V]$ and $\norm{E[V]} \leq E[\norm{V}]$, 
where $E$ denotes expectation.
Furthermore, we have the following strong law of large numbers.

\begin{theorem} \label{thm:slln}
    Let $(V_k)_{k \in \N}$ be a sequence of independent copies of $V$ and for $n \geq 1$ let  $S_n = \sum_{k=1}^n V_k$.
    If $E[\norm{V}] < \infty$ then the sequence $(\frac{1}{n}S_n)$ of $\mathcal{B}^S$-valued random variables converges to $E[V]$ almost surely. 
\end{theorem}

We also have the following central limit theorem.

\begin{theorem} \label{thm:clt}
    Assume that $\mathcal{B}$ has type 2 (e.g. $\mathcal{B} = L^p(\N \times \R)$, where $p \geq 2$).
    If $E[\norm{V}] < \infty$ and $E[\norm{V}^2] < \infty$ then $(\sqrt{n}(\frac{1}{n} S_n - E[V]))$ converges in distribution to a Gaussian random variable with the same covariance structure as $V$.
\end{theorem}

\subsection{Distribution of our random feature maps}

Here we prove a main property of the expectations of our random feature maps.

\begin{theorem} \label{thm:distribution}
    Assume that $E(F(\alpha(X)))$ exists. Then
    $E(F(\alpha(X))) = \sum_K E(\eta_K(X)) = \sum_{\mathcal{C}} E(\eta_{\mathcal{C}}(X))$.
\end{theorem}

\begin{proof}
    Extending \eqref{cd:heatmap-simplices} and \eqref{cd:heatmap-pixels}, we have the following commutative diagrams.

    \begin{equation*}
    \begin{tikzcd}
        \Omega \ar[r,"X"] & \R^m \ar[r,"F \circ \alpha"] & \mathcal{B}\\
        \Omega \ar[r,"X"] \ar[u,"\Id"] & \R^m \ar[u,"\Id"] \ar[r,"\eta_K"] & \mathcal{B}^K \ar[u,"\Sigma_K"']
    \end{tikzcd} \quad
    \begin{tikzcd}
        \Omega \ar[r,"X"] & \R^m \ar[r,"F \circ \alpha"] & \mathcal{B}\\
        \Omega \ar[r,"X"] \ar[u,"\Id"] & \R^m \ar[u,"\Id"] \ar[r,"\eta_{\mathcal{C}}"] & \mathcal{B}^{\mathcal{C}} \ar[u,"\Sigma_{\mathcal{C}}"']
    \end{tikzcd}
    \end{equation*}

    That is, $F(\alpha(X)) = \sum_K \eta_K X = \sum_{\mathcal{C}} \eta_{\mathcal{C}} X$.
    Since these random variables are equal, they have the same expectation.
    By the linearity of expectation, the result follows.
\end{proof}

\section{Stabilizing the chain, simplex, and pixel/voxel feature maps}
\label{sec:stabilization}

In this section, we take the expectation of the random feature maps of the previous section to obtain functions that are Lipschitz-continuous with respect to changes of the weight.

\subsection{Means of random feature maps}

Let $K$ be a finite simplicial complex with position function $v: K_0 \to \R^d$.
Let $F: \mathcal{D} \to \mathcal{B}$ be a feature map with a labeled decomposition $\hat{F}: \hat{\mathcal{D}}_K \to \mathcal{B}^{S(K)}$.
Choose a persistence algorithm $\mathcal{A}$ and a linear extension $\trianglelefteq$ of the partial order $\subseteq$ on $K$.
Following the previous section, we have 
a chain feature map $\eta_{S(K)}: \R^m \to \mathcal{B}^{S(K)}$,
a simplex feature map $\eta_K: \R^m \to \mathcal{B}^K$, and 
a pixel/voxel feature map $\eta_\mathcal{C} : \R^m \to \mathcal{B}^\mathcal{C}$.
Let $X: (\Omega,\mathcal{F},P) \to \R^m$ be a random (parameter for a) weight.
Then we have 
the random chain feature map $\eta_{S(K)} \circ X$,
the random simplex feature map $\eta_K \circ X$, and 
the random pixel/voxel feature map $\eta_\mathcal{C} \circ X$.

Let $g: \mathcal{B} \to \R$ be a linear functional.
Then we have a random variable $g^{S(K)} \circ \eta_{S(K)} \circ X: (\Omega,\mathcal{F},P) \to \R^{S(K)}$,
where $g^{S(K)}$ is given by applying $g$ to each coordinate.
We have the following commutative diagram.
\begin{equation*}
    \begin{tikzcd}
        \Omega \ar[r,"X"] & \R^m \ar[r,"F \circ \alpha"] & \mathcal{B} \ar[r,"g"] & \R \\
        \Omega \ar[r,"X"] \ar[u,"\Id"] & \R^m \ar[u,"\Id"] \ar[r,"\eta_{S(K)}"] & \mathcal{B}^{S(K)} \ar[u,"\Sigma_{S(K)}"'] \ar[r,"g^{S(K)}"] & \R^{S(K)} \ar[u,"\Sigma_{S(K)}"']
    \end{tikzcd}
\end{equation*}
Similarly, we have the random variables
$g^K \circ \eta_K \circ X: (\Omega,\mathcal{F},P) \to \R^K$ and
$g^\mathcal{C} \circ \eta_\mathcal{C} \circ X: (\Omega,\mathcal{F},P) \to \R^\mathcal{C}$.
When they are defined, these maps have means
$E[g^{S(K)} \circ \eta_{S(K)} \circ X]$,
$E[g^K \circ \eta_K \circ X]$, and
$E[g^\mathcal{C} \circ \eta_\mathcal{C} \circ X]$.

Now, let us consider the random feature maps and their means coordinatewise.
Recall that $\eta_{S(K)} = \hat{F} \circ \hat{\alpha}: \R^m \to \mathcal{B}^{S(K)}$.
For $\gamma \in S(K)$, the $\gamma$ coordinate of $\eta_{S(K)}$ is given by $(\eta_{S(K)})_\gamma = \hat{F}_\gamma \circ \hat{\alpha}: \R^m \to \mathcal{B}$.
Therefore the $\gamma$ coordinate of 
$g^{S(K)} \eta_{S(K)} X$ is given by $g \hat{F}_\gamma \hat{\alpha} X$.
\begin{equation*}
        (\Omega,\mathcal{F},P) \xrightarrow{X} \R^m \xrightarrow{(\eta_{S(K)})_\gamma = \hat{F}_\gamma \circ \hat{\alpha}} \mathcal{B} \xrightarrow{g} \R
\end{equation*}
Hence the $\gamma$ coordinate of 
$E[g^{S(K)} \eta_{S(K)} X]$ is given by $E[g \hat{F}_\gamma \hat{\alpha} X]$.
Similarly, for $\sigma \in K$, the $\sigma$ coordinate of 
$E[g^K \eta_K X]$ is given by $E[g \bar{F}_\sigma \hat{\alpha} X]$
and for $c \in \mathcal{C}$, the $c$ coordinate of
$E[g^{\mathcal{C}} \eta_{\mathcal{C}} X]$ is given by $E[g \tilde{F}_c \hat{\alpha} X]$.

\subsection{Convolution with a kernel}
\label{sec:convolution}

For $\gamma \in S(K)$, $\sigma \in K$, or $c \in \mathcal{C}$,
let $h: \R^m \to \R$ be given by 
$g \hat{F}_\gamma \hat{\alpha}$,
$g \bar{F}_\sigma \hat{\alpha}$, or
$g \tilde{F}_c \hat{\alpha}$, respectively.

Let $\kappa$ be a kernel on $\R^m$.
Let $\eps \sim \kappa$, 
i.e. $\eps$ is a random variable $\eps: (\Omega,\mathcal{F},P) \to \R^m$, with probability distribution $\kappa$.
Let $w$ be a weight on $K$.
Let $Y_w: (\Omega,\mathcal{F},P) \to \R$ be the random variable given by $Y_w = h(w-\eps)$, i.e.
\begin{equation*}
    Y_w: (\Omega,\mathcal{F},P) \xrightarrow{w-\eps} \R^m \xrightarrow{h} \R.
\end{equation*}
Then 
\begin{equation} \label{eq:convolution}
    E[Y_w] = \int_{\R^d} h(w-x)\kappa(x) dx = (h * \kappa)(w),
\end{equation}
where $*$ denotes convolution, assuming this integral exists.


\subsection{Mean feature maps} \label{sec:mean-feature-maps}

Assume that the integral in \eqref{eq:convolution} exists for all $w \in \R^m$. 
Then we define the \emph{mean chain feature map}
$\hat{H}: \R^m \to \R^{S(K)}$ by 
$\hat{H}_\gamma = (g \hat{F}_\gamma \hat{\alpha}) * \kappa$, where $\gamma \in S(K)$.
Similarly, we define the \emph{mean simplex feature map}
$\bar{H}: \R^m \to \R^K$ by 
$\bar{H}_\sigma = (g \bar{F}_\sigma \hat{\alpha}) * \kappa$, where $\sigma \in K$, and
we define the \emph{mean pixel/voxel feature map}
$\tilde{H}: \R^m \to \R^{\mathcal{C}}$ by 
$\tilde{H}_c = (g \tilde{F}_c \hat{\alpha}) * \kappa$, where $c \in \mathcal{C}$.

\subsection{Stability of the mean feature maps} \label{sec:stability}

We will show that under suitable hypotheses, the mean feature maps of the previous section are Lipschitz.

Let $h: \R^m \to \R$ denote 
$g \hat{F}_\gamma \hat{\alpha}$, where $\gamma \in S(K)$,
$g \bar{F}_\sigma \hat{\alpha}$, where $\sigma \in K$, or
$g \tilde{F}_c \hat{\alpha}$, where $c \in \mathcal{C}$.
In this section we do not assume that $g$ is linear.
Let $\kappa$ be a kernel on $\R^m$.
Since $\kappa$ is integrable, if $h$ is essentially bounded then $h * \kappa$ is defined everywhere and is uniformly continuous~\cite[255K]{MR2462280}.

\begin{theorem}[{\cite[473D(d)]{MR2462372}}] 
    If $\norm{h}_1 = a$ and $\kappa$ is $b$-Lipschitz, then $h * \kappa$ is $ab$-Lipschitz.
\end{theorem}

For $w \in \R^m$ and $a > 0$, let $B_{a}(w)$ denote the ball centered at $w$ with radius $a$.
Let $V$ denote the volume of $B_1(0)$.

\begin{theorem}[{\cite[Theorem 4]{BBW:2020}}]
    Let $w \in \R^m$ and $a > 0$.
    If $\norm{h}_\infty \leq M$ on $B_{2a}(w)$, 
    $\kappa$ is $b$-Lipschitz and $\supp(K) \subseteq B_a(w)$,
    then $h * \kappa$ is $2Mba^{\abs{K}}V$-Lipschitz in $B_a(w)$.
\end{theorem}

\begin{theorem}[{\cite[Theorem 5]{BBW:2020}}]
    If $\norm{h}_\infty \leq M$ and $\int \norm{\kappa(s+t)-\kappa(s)} ds \leq b\norm{t}$ for all $t \in \R^m$ then
    $h*\kappa$ is $Mb$-Lipschitz.
\end{theorem}

Next we consider the cases of the Gaussian, triangular, and Epanechnikov kernels. 
Let $a > 0$.

\begin{theorem} \label{thm:stability-gaussian}
    Let $\kappa_a$ denote the (isotropic) Gaussian kernel on $\R^m$ \eqref{eq:gaussian-kernel}.
    Assume that $\norm{h}_{\infty} \leq M$.
    Then $h * \kappa_a$ is $\frac{2M}{a\sqrt{2\pi}}$-Lipschitz.
\end{theorem}

\begin{proof}
    This result follows directly from \cite[Corollary 3]{BBW:2020}.
\end{proof}

\begin{theorem} \label{thm:stability-triangular}
    Let $\kappa_a$ denote the triangular kernel on $\R^m$ \eqref{eq:triangular-kernel}.
    Assume that $\norm{h}_{\infty} \leq M$ on $B_{2a}(w)$.
    Then $h * \kappa_a$ is $\frac{2M(\abs{K}+1)}{a}$-Lipschitz.
\end{theorem}

\begin{proof}
    This result follows directly from \cite[Corollary 1]{BBW:2020}.
\end{proof}

\begin{theorem} \label{thm:stability-epanechnikov}
    Let $\kappa_a$ denote the Epanechnikov kernel on $\R^m$ \eqref{eq:epanechnikov-kernel}.
    Assume that $\norm{h}_{\infty} \leq M$ on $B_{2a}(w)$.
    Then $h * \kappa_a$ is $\frac{2M(\abs{K}+2)}{a}$-Lipschitz.
\end{theorem}

\begin{proof}
    This result follows directly from \cite[Corollary 2]{BBW:2020}.
\end{proof}

\subsection{Stable feature maps in practice} \label{sec:in-practice}

Let $\kappa$ denote either the Gaussian, triangular, or Epanechnikov kernel.

Consider 
the mean chain feature map $\hat{H}: \R^m \to \R^{S(K)}$, where 
for $\gamma \in S(K)$,
$\hat{H}_\gamma: \R^m \to \R$ is given by $\hat{H}_\gamma = (g \hat{F}_\gamma \hat{\alpha}) * \kappa$,
the mean simplex feature map $\bar{H}: \R^m \to \R^K$, where 
for $\sigma \in K$,
$\bar{H}_\sigma: \R^m \to \R$ is given by $\bar{H}_\sigma = (g \bar{F}_\sigma \hat{\alpha}) * \kappa$, or
the mean pixel/voxel feature map $\tilde{H}: \R^m \to \R^{\mathcal{C}}$, where 
for $c \in \mathcal{C}$,
$\tilde{H}_c: \R^m \to \R$ is given by $\tilde{H}_c = (g \tilde{F}_c \hat{\alpha}) * \kappa$.
In the previous section, we showed that under mild hypotheses, these maps are stable.

Let $h$ denote $g \hat{F}_\gamma \hat{\alpha}$, $g \bar{F}_\sigma \hat{\alpha}$, or $g \tilde{F}_c \hat{\alpha}$.
Using \eqref{eq:convolution}
to obtain an analytic expression of $h * \kappa$ is not feasible.
However, we may approximate $g = h * \kappa$ by the following simple procedure.

We wish to approximate $E[Y_w]$.
Draw a sample $\eps_1,\eps_2,\ldots$, where $\eps_i \sim \kappa$ are independent. 
Then, by the strong law of large numbers,
the empirical mean $\frac{1}{m} \sum_{i=1}^m h(w-\eps_i)$,
converges almost surely to $E[Y_w]$.
That is, we have the following.

\begin{theorem} \label{thm:empirical-mean}
    Let $w$ be a weight on $K$, and let $\eps_1,\eps_2,\ldots$ be drawn independently from $\kappa$. Then
    \begin{equation*}
        \frac{1}{m} \sum_{i=1}^m h(w-\eps_i) \to g(w)
    \end{equation*}
    almost surely.
\end{theorem}

\subsection{Averaging over a family of examples} 
\label{sec:average}

Above, we showed how to produce a stable summary of the part of a particular test case responsible for a given classification. 
However, in some applications, we would instead like an average summary over all samples in a class. 
We will now show how to do this.

As in \cref{sec:in-practice},
for $\gamma \in S(K)$, $\sigma \in K$, or $c \in \mathcal{C}$, let $h: \R^m \to \R$ be one of the maps 
$g \hat{F}_\gamma \hat{\alpha}$,
$g \bar{F}_\sigma \hat{\alpha}$, or
$g \tilde{F}_c \hat{\alpha}$, respectively.
Previously, we fixed a weight $w \in \R^d$, took $\eps \sim \kappa$, and showed that $E[h(w-\eps)] = (h * \kappa)(w)$.
Now, instead assume that we have a random weight $W: (\Omega,\mathcal{F},P) \to \R^d$. 
Then we have the composite random variable $h(W) = h \circ W: (\Omega,\mathcal{F},P) \to \R$.
Assume that $E[h(W)]$ exists.

\begin{theorem} \label{thm:empirical-mean-family}
    Let $w_1,w_2,\ldots \in \R^d$, where each $w_k$ is drawn independently from $W$.
    Then, with probability $1$,
    \begin{equation}
        \lim_{n \to \infty} \left(\frac{1}{n} \sum_{k=1}^n h(w_k)\right) = E[h(W)].
    \end{equation}
\end{theorem}

\begin{proof}
    Let $W_1,W_2,\ldots$ be a sequence of independent copies of $W$.
    For $n \geq 1$, let $S_n = \sum_{k=1}^n h(W_k)$.
    Then by the strong law of large numbers, $(\frac{1}{n} h(S_n))$ converges to $E[h(W)]$ almost surely.
\end{proof}

That is, 
the empirical mean is stable in the sense that
with probability $1$ 
it converges to $E[h(W)]$ independently of the choice of draw from $W^\infty$.




\section{Computational examples}
\label{sec:examples}

In this section we provide a few computational examples to demonstrate our methods.
We use standard methods of TDA to solve various classification problems and then use mean feature maps to explain the classification.

\subsection{Classification for uniform samples on topologically distinct spaces}
\label{sec:annulus-and-double-annulus}

In our first example, we consider a classification problem for points sampled from the uniform distribution on two compact regions in the plane, an annulus and a double annulus, together with Gaussian noise.

\begin{figure}[!htb]
    \centering
    \includegraphics[height=35mm]{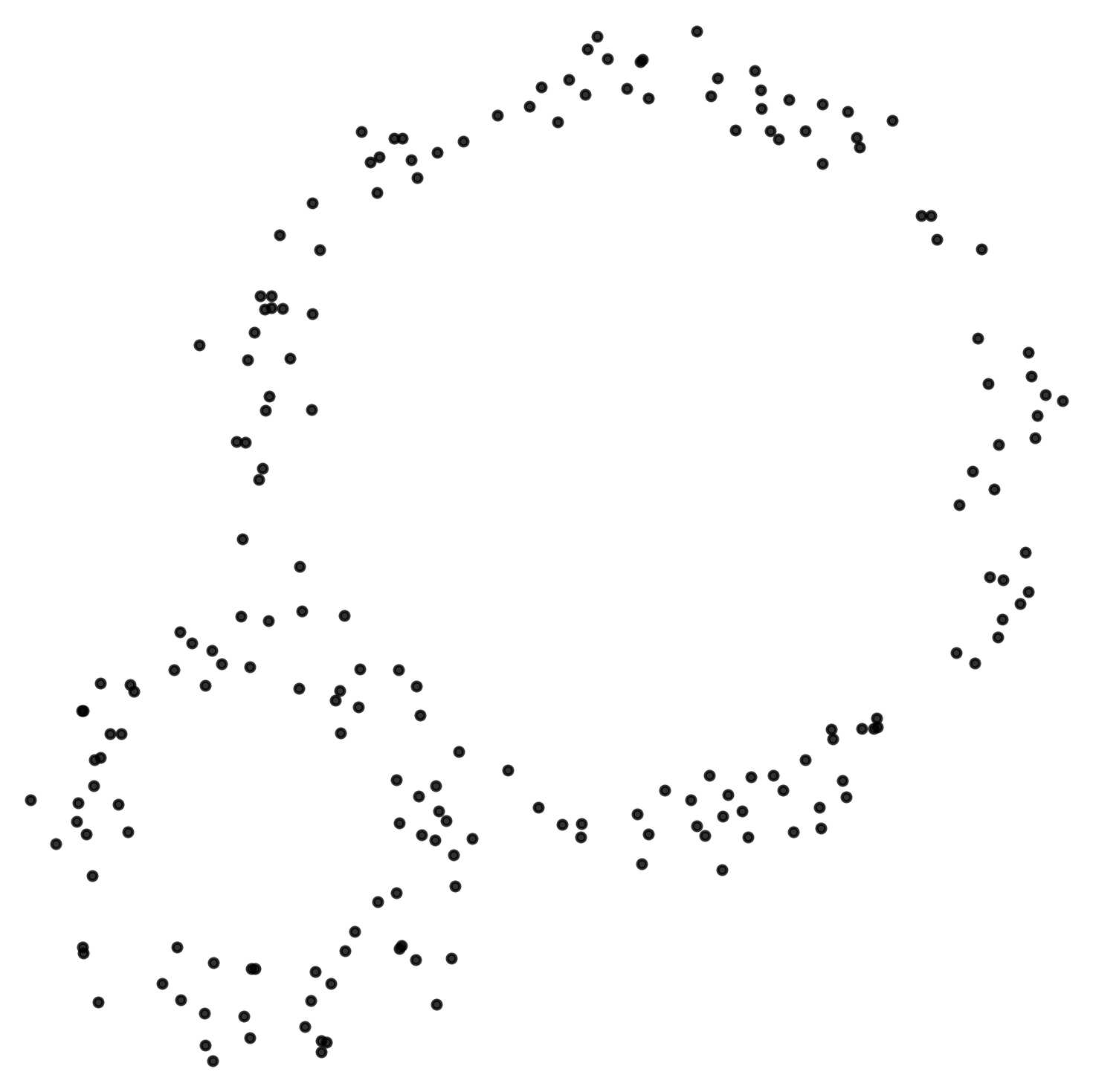}
    \quad \quad
    \includegraphics[height=35mm]{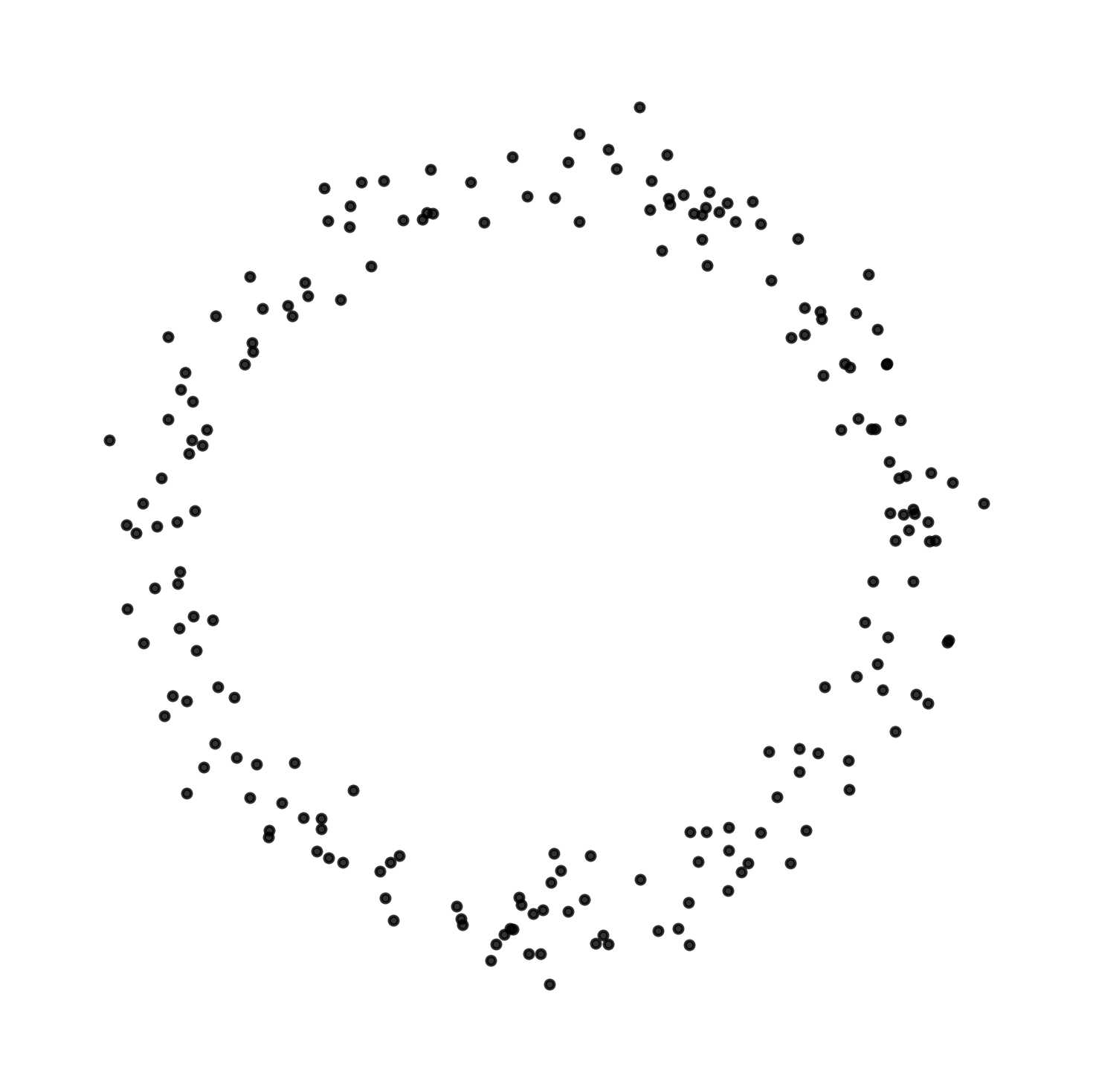}
    \caption{
    200 points sampled independently from the uniform distribution on a double annulus (left) and annulus (right) together with Gaussian noise.}
    \label{fig:annulus}
\end{figure}

The double annulus is obtained from the annulus by taking a union with a smaller annulus. 
In each case, we sample $200$ points independently from the uniform area measure together with fixed Gaussian noise. See \cref{fig:annulus}.

Our training data consists of $100$ point clouds from each of the classes.
For each of these point clouds, $X$, we construct the \emph{Delaunay complex}. 
In detail, $X$ has a corresponding \emph{Voronoi} tessellation of the plane into closed polyhedrons, consisting of points whose nearest neighbor in $X$ includes a specified point of $X$.
The Delaunay complex, $K$, is the simplicial complex that is dual to the the Voronoi tessellation, and it provides a canonical triangulation of the convex hull of $X$.
Note that we are assuming that the points in $X$ are in general position, which occurs with probability $1$.
Then, we consider \emph{Voronoi balls} given by intersecting balls centered at the points of $X$ and their corresponding Voronoi polyhedrons.
We weight each simplex in the Delaunay complex by the smallest radius for which the Voronoi balls of the vertices of this simplex intersect.
For each Delaunay complex, $(K,w)$, we compute the persistence diagram for homology in degree $1$ and the corresponding persistence landscape.


We apply support vector machines (SVM) with a linear kernel to the persistence landscapes to produce a classifier consisting of a normal vector and a bias for the separating hyperplane.

For testing data, we compute the persistence landscape and take the inner product with the normal vector and subtract the bias and the resulting sign of the real number determines which class we judge the sample to belong to.
We then use the results of this paper to visualize the simplices that produce this real number.

Our construction of $X$ induces a total order on $K$ as follows.
The points in $X$ were chosen one at a time, so they have a total order.
This order induces an order on the simplices of $K$, which are subsets of $X$, given by lexicographic order, where the simplices of $K$ are written using the order on $X$.
We also fix a choice of persistence algorithm $\mathcal{A}$.
See \url{https://github.com/Himanshu484/persistence_heatmap}.

Next, we compute the persistence diagrams and persistence landscapes, labeled by the birth simplices, death simplices, representative cycles and bounding chains.
We then distribute the number obtained from the SVM classifier to the simplices of $K$ and furthermore to the pixels of a rectangle containing all of the sampled points.

\begin{figure}[!htbp]
    \centering
    \includegraphics[height=35mm]{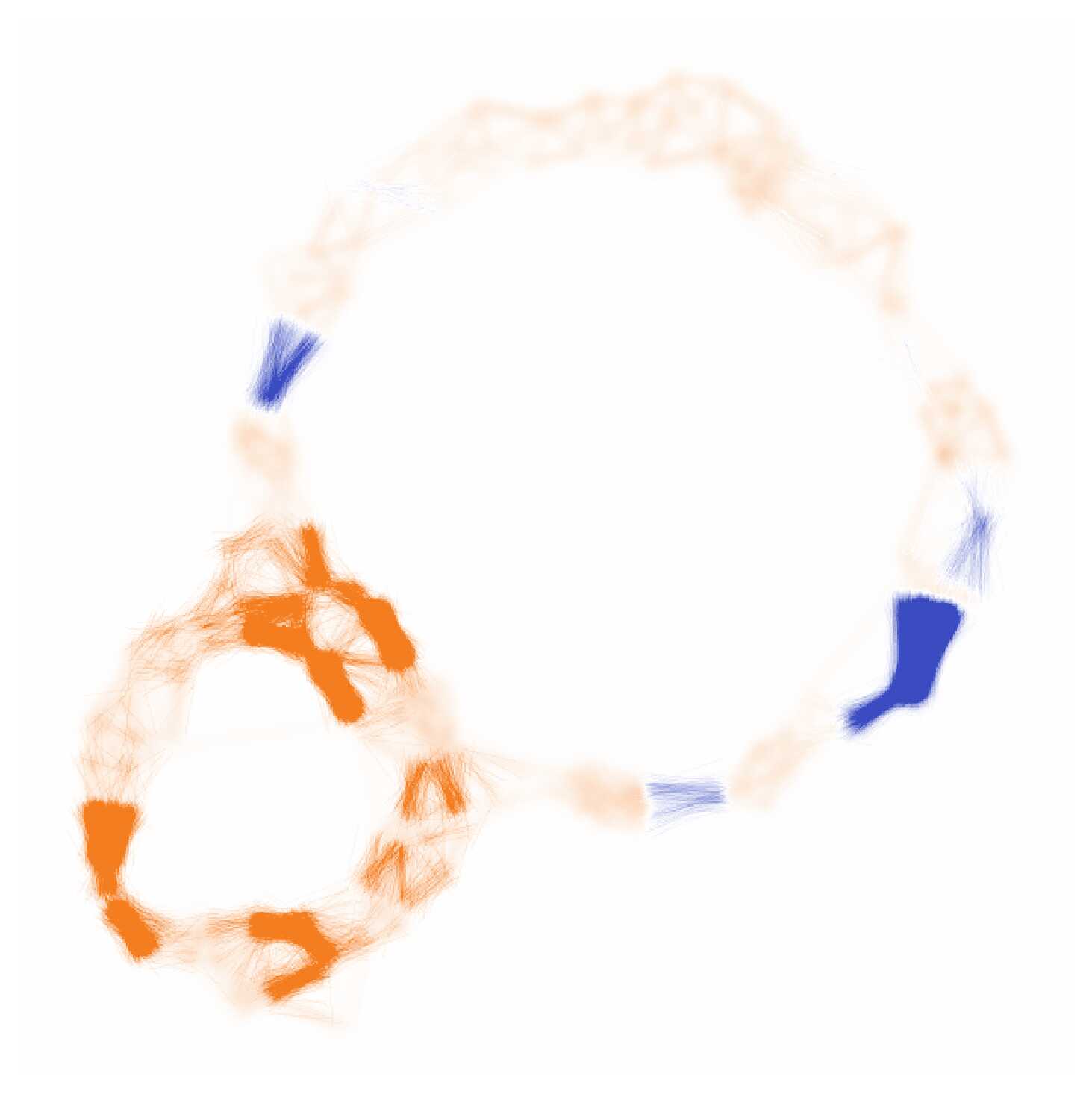}
    \includegraphics[height=35mm]{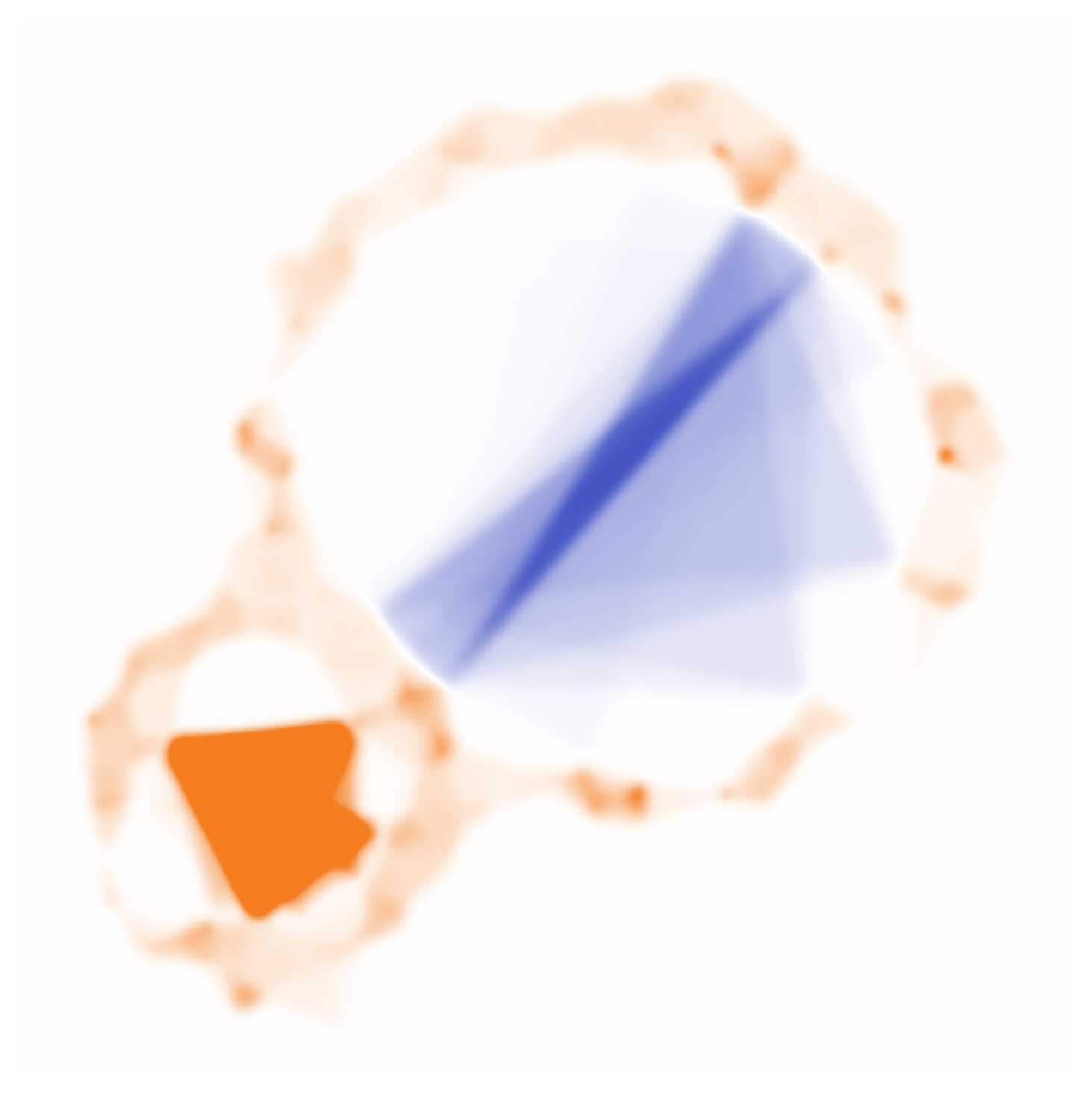}
    \includegraphics[height=35mm]{figures/RP_dual_double.jpg}
    \includegraphics[height=35mm]{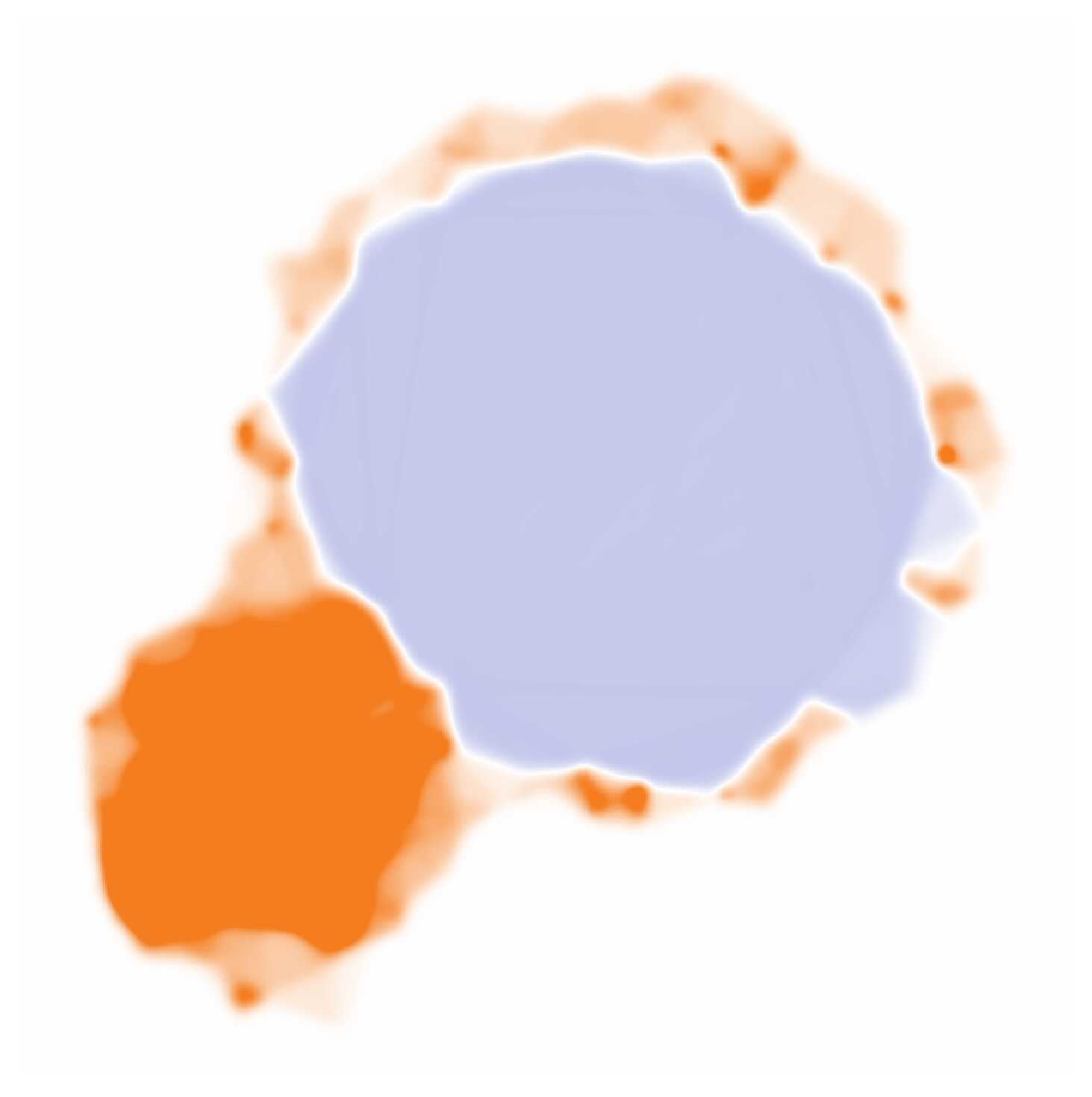}\\
    \includegraphics[height=35mm]{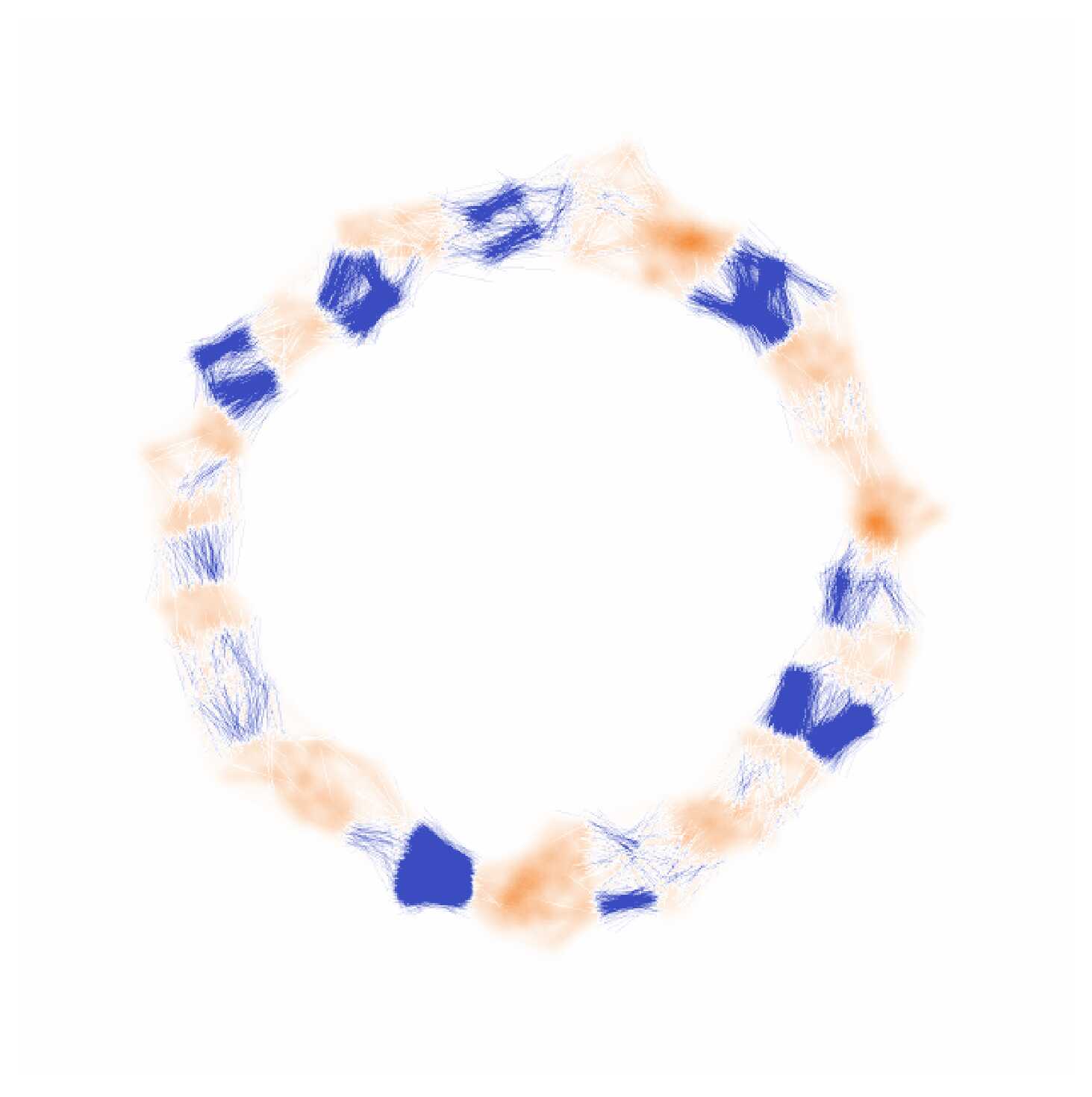}
    \includegraphics[height=35mm]{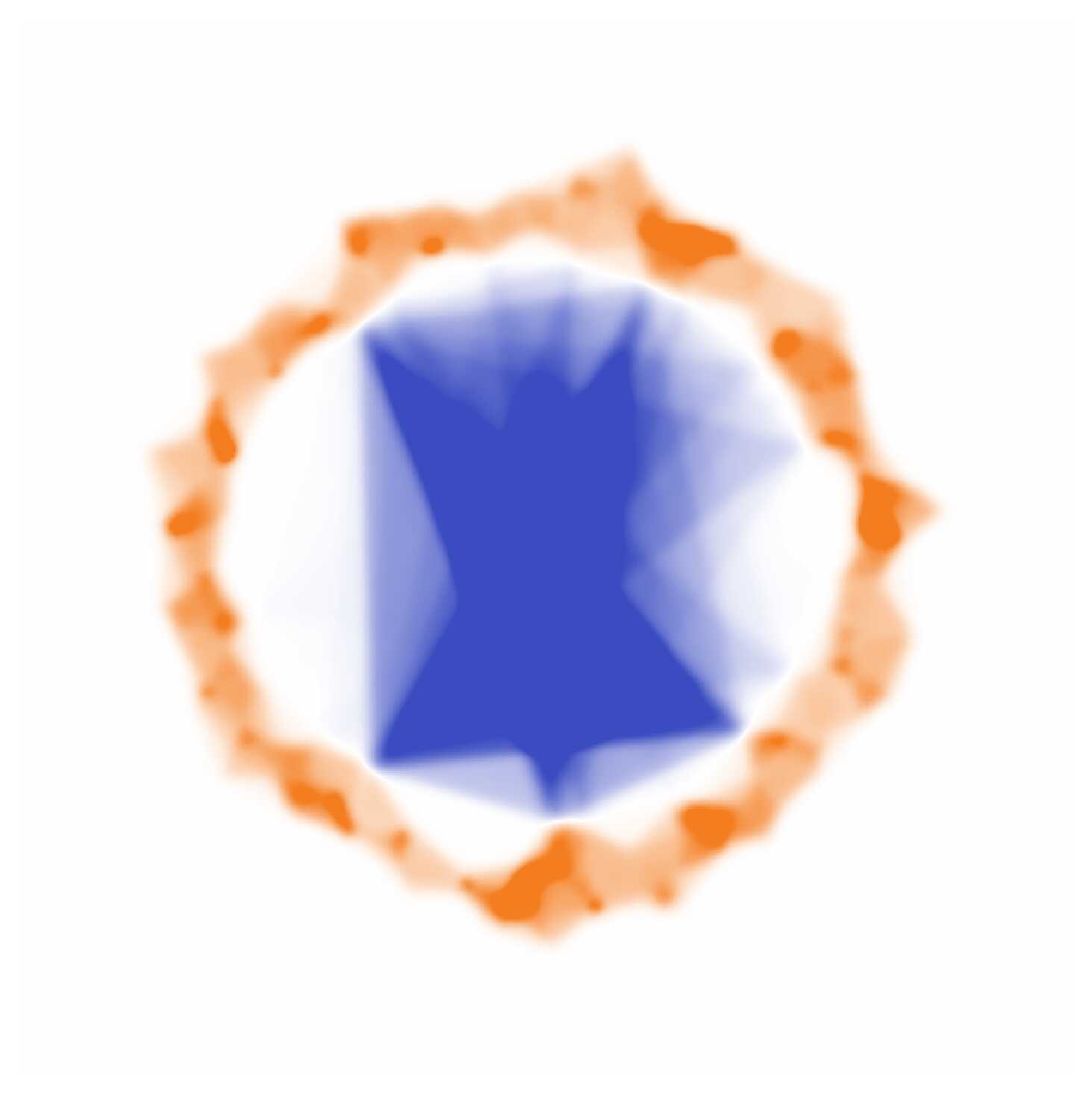}
    \includegraphics[height=35mm]{figures/RP_dual_single.jpg}
    \includegraphics[height=35mm]{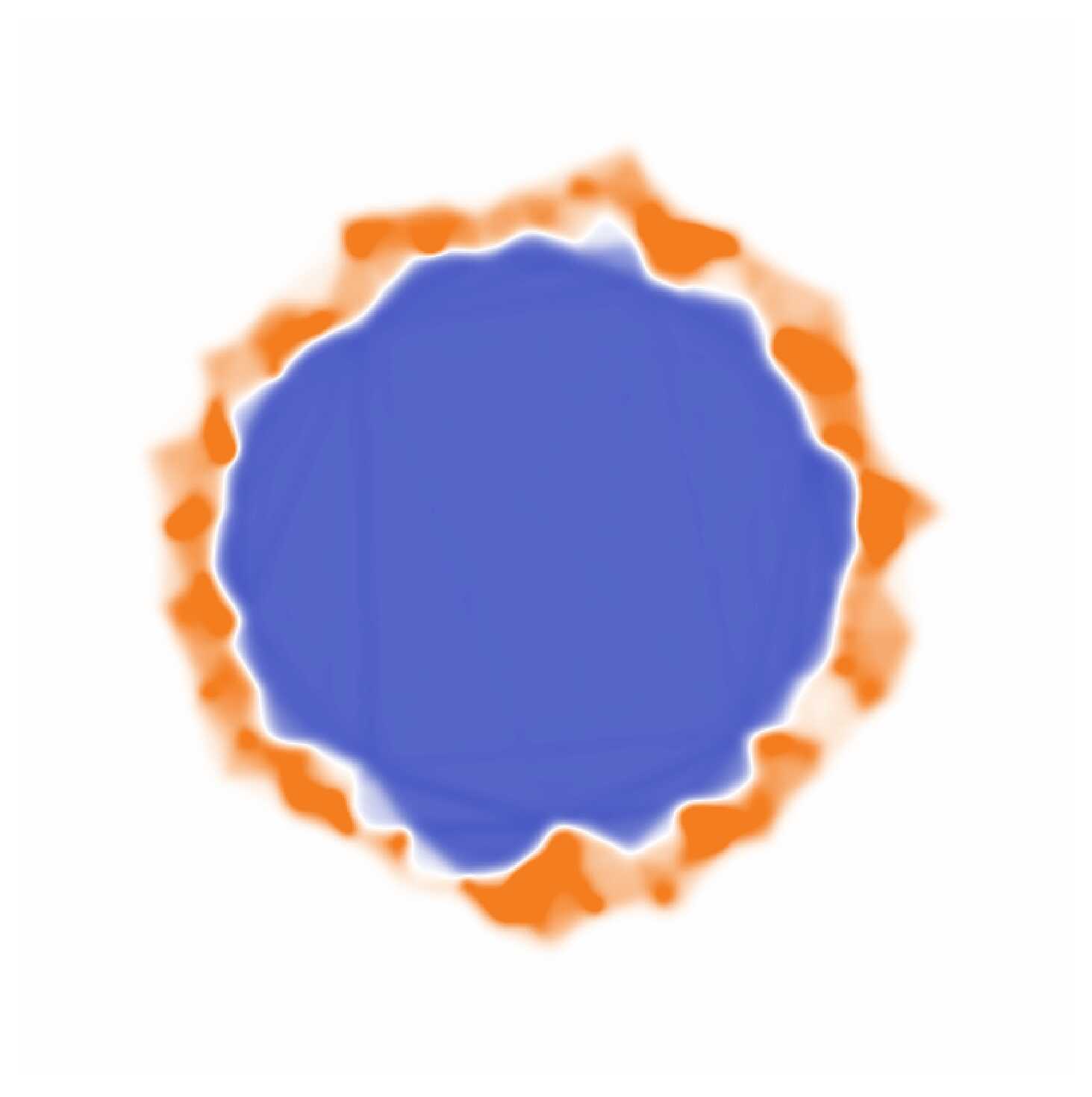}
    \caption{Visualization of the learned classifier that differentiates between points sampled from the double annulus (top row) and points sampled from the annulus (bottom row), using homology in degree $1$ and labeled persistence landscapes.
    We use a orange-to-blue color gradient in which orange corresponds to the double annulus and blue corresponds to the annulus.
    From left to right, the columns use birth edges, death triangles, representative cycles, and bounding chains, respectively.
    In each of the images, the sum of the values of the pixels equals the empirical mean value of the classifier.
    }
    \label{fig:annulus-visualization-single}
\end{figure}

\begin{figure}[!htbp]
    \centering
    \includegraphics[height=35mm]{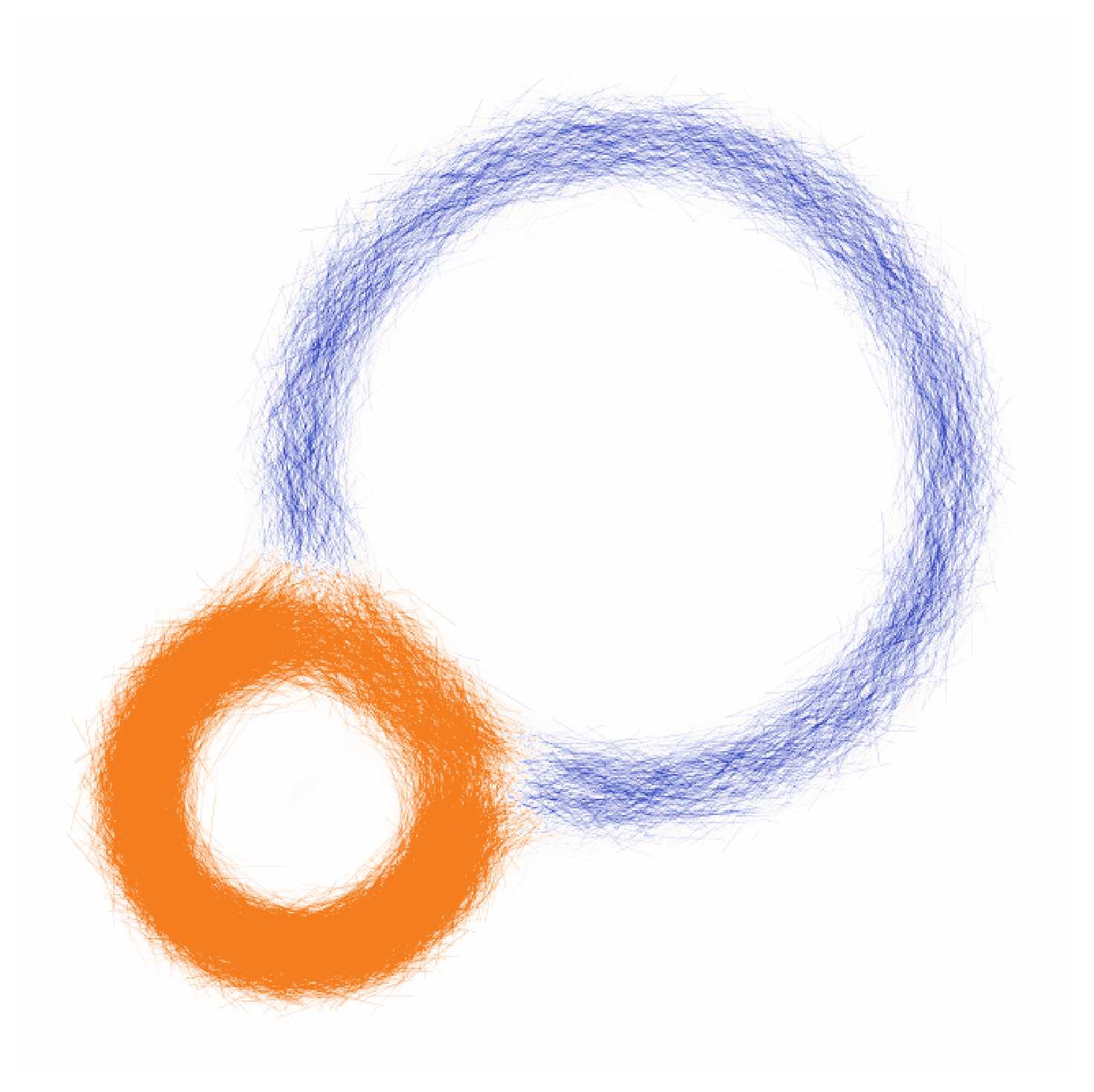}
    \includegraphics[height=35mm]{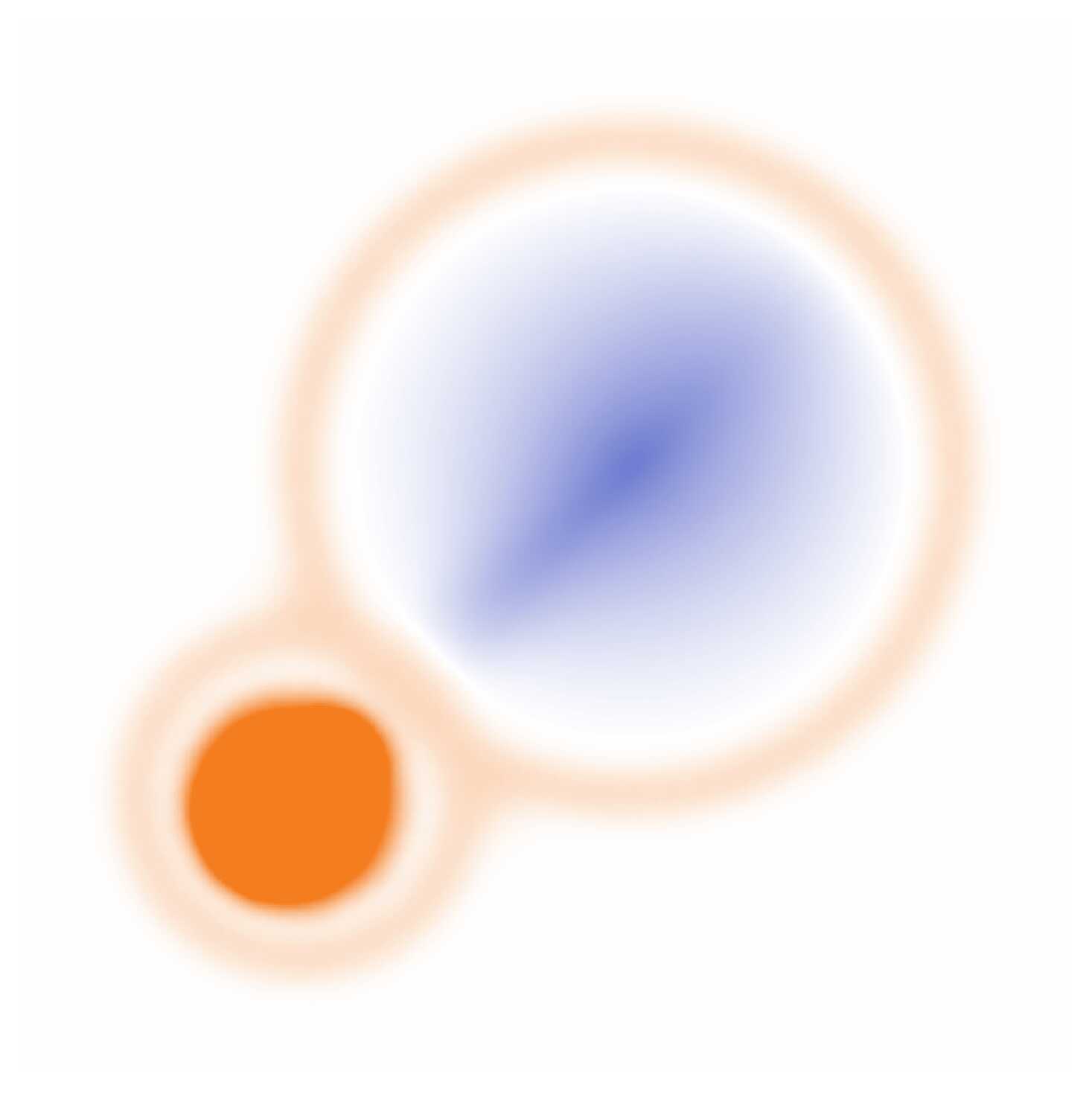}
    \includegraphics[height=35mm]{figures/RP_dual_double_all.jpg}
    \includegraphics[height=35mm]{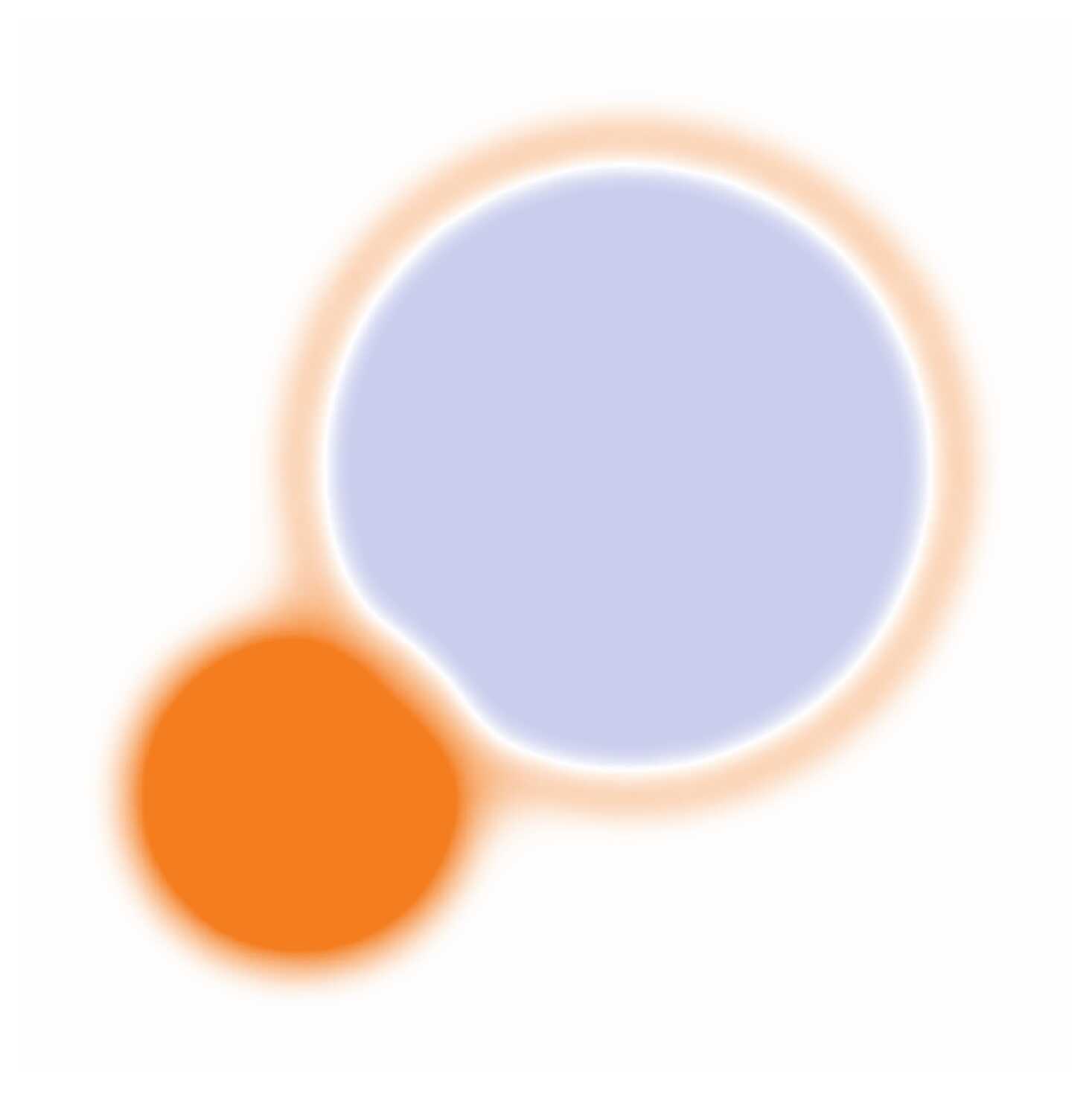}\\
    \includegraphics[height=35mm]{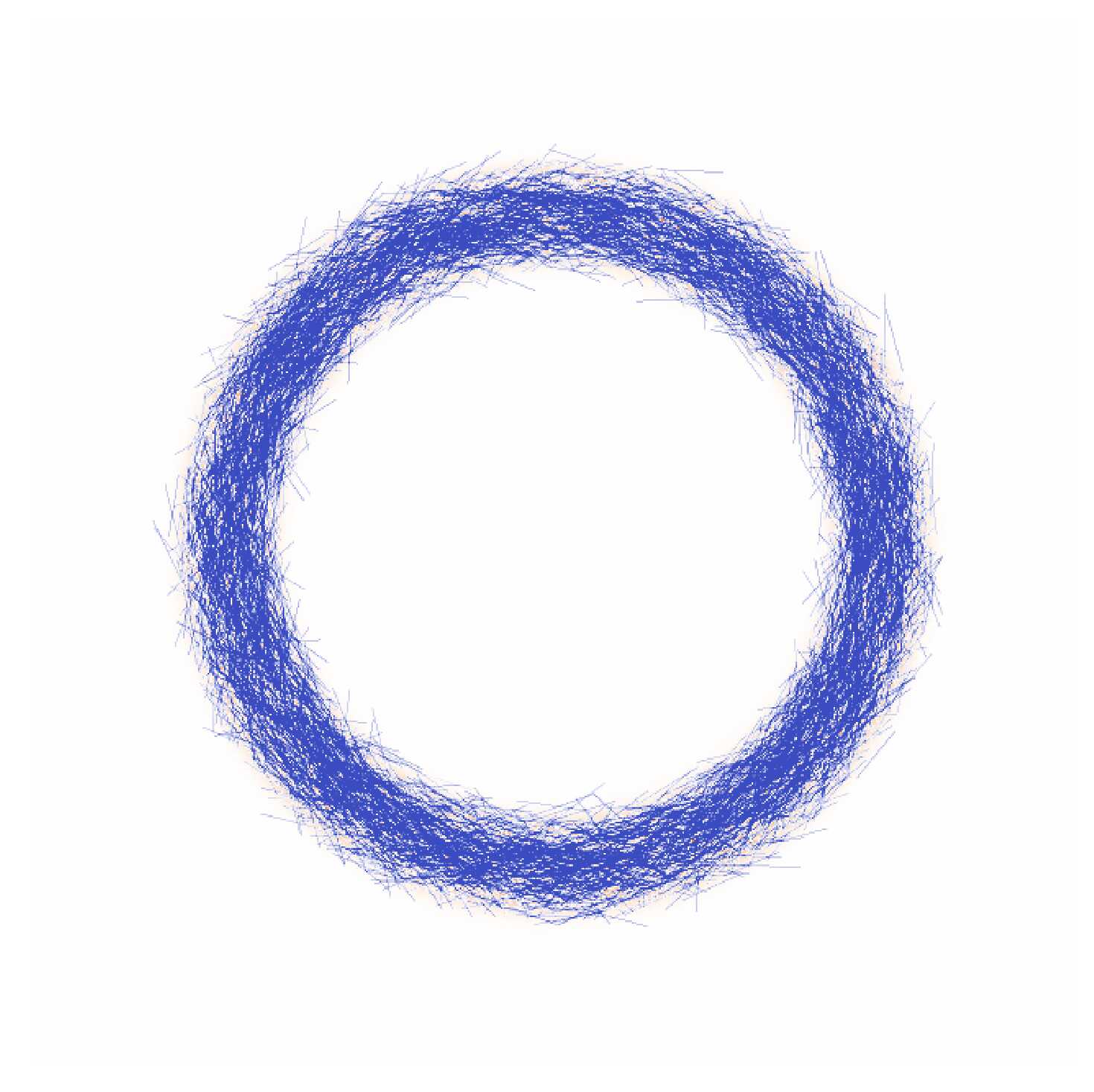}
    \includegraphics[height=35mm]{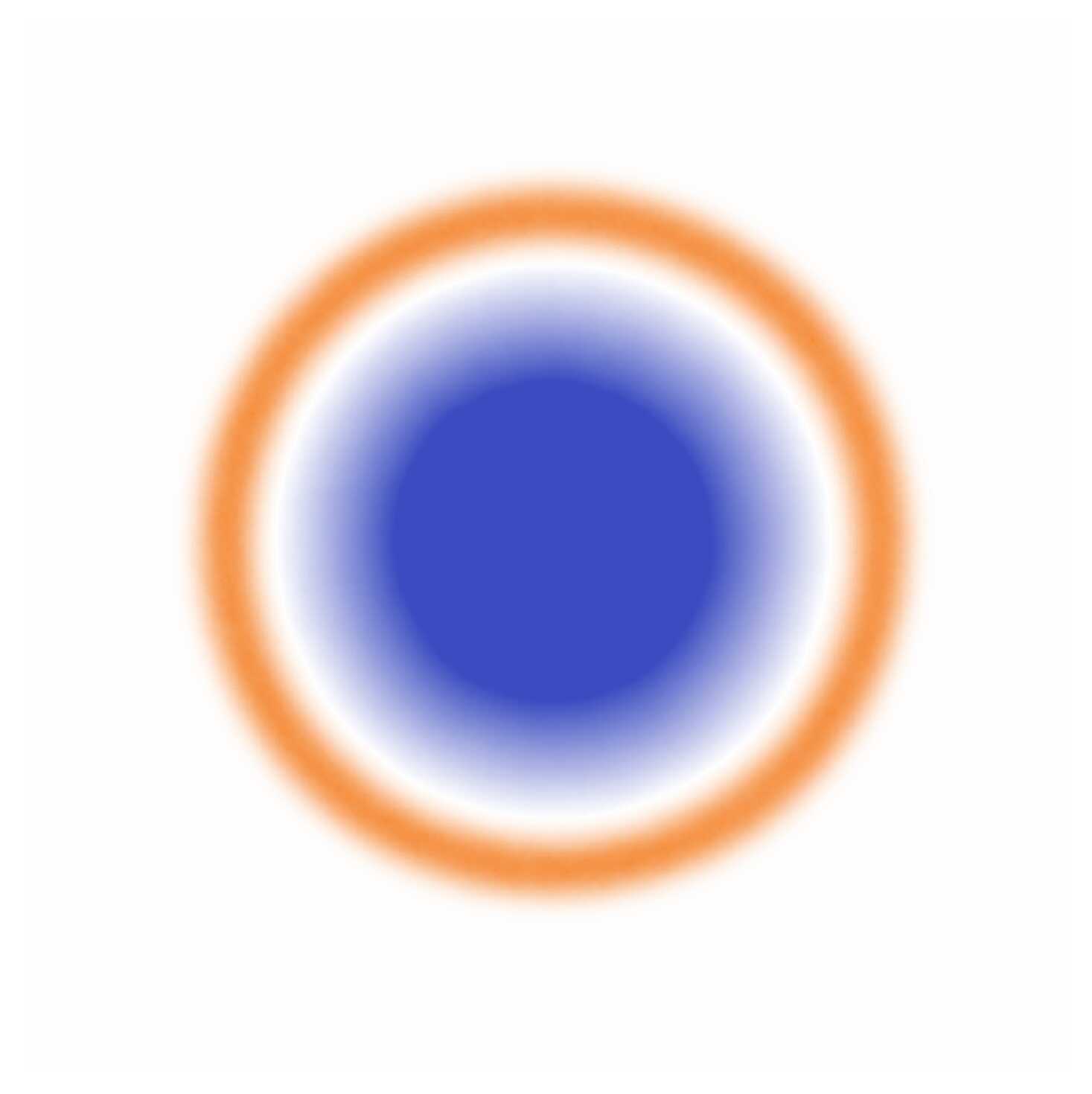}
    \includegraphics[height=35mm]{figures/RP_dual_single_all.jpg}
    \includegraphics[height=35mm]{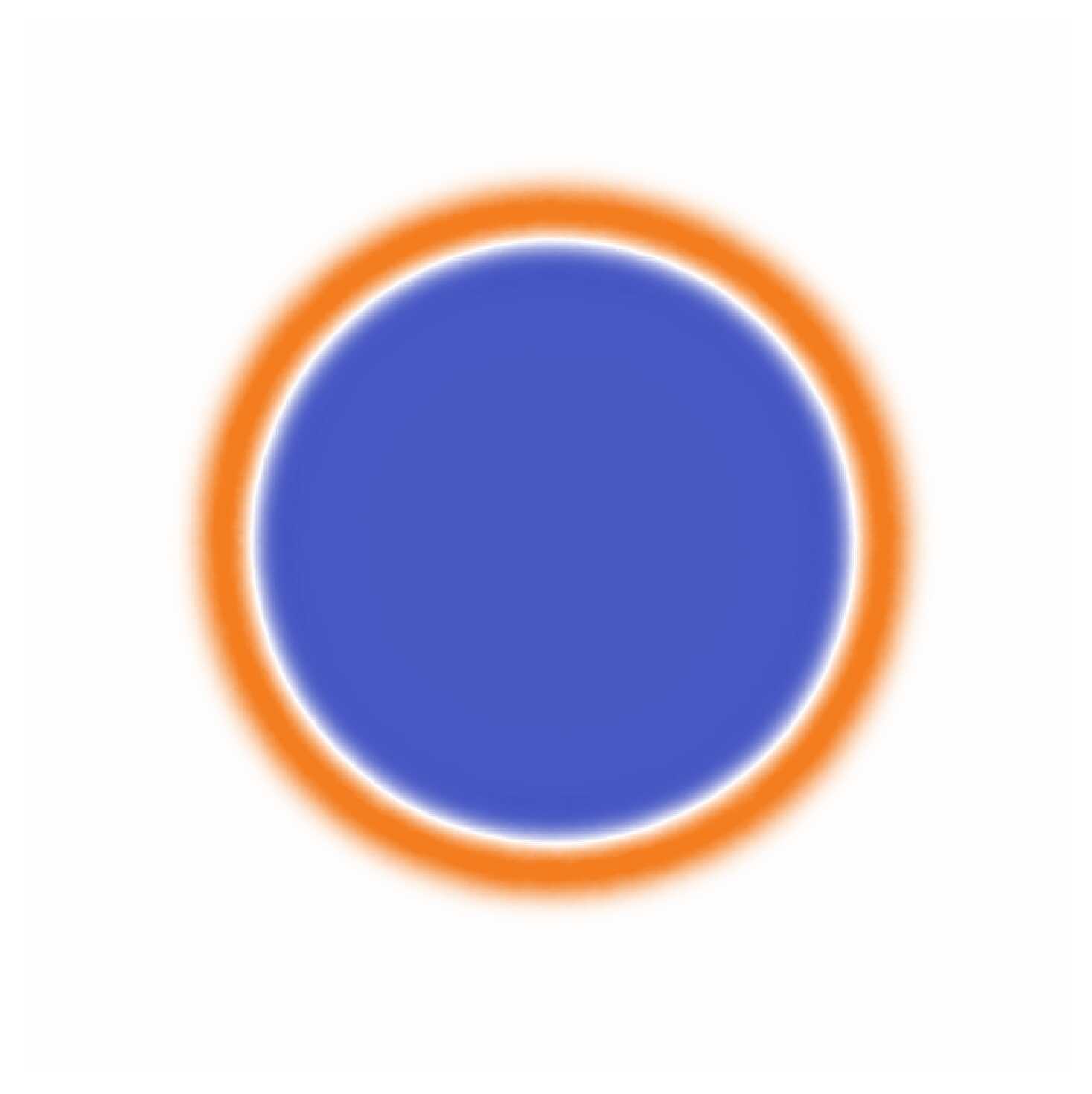}
    \caption{
    Same as \cref{fig:annulus-visualization-single} but instead of averaging over perturbations of a single point cloud, we average over point clouds.
    }
    \label{fig:annulus-visualization-mean}
\end{figure}

Finally, we produce empirical means of two sets of random feature maps.
For the first
(see \cref{fig:annulus-visualization-single}),
we fix a single point cloud sampled from each class and perturb it $10000$ times:
each perturbation adds independent Gaussian noise, with mean $0$ and standard
deviation $\sigma = 0.05$, to each Cartesian coordinate of every point in the
cloud.
In each of the images, the sum of the values of the pixels equals the empirical mean value of the classifier.
For the second
(see \cref{fig:annulus-visualization-mean}),
we instead draw $10000$ independent point clouds from each class and average their
feature maps.
We remark that the empirical means for $N=100$ look similar to the ones given here.

Note that the visualizations for the larger annulus, which is common to both classes, is  similar for both classes. However, the visualizations for the smaller annulus in the double annulus, which is not present in the annulus samples, is strongly colored in the color of the double annulus class. 




\subsection{Classification for samples from different distributions on a topologically trivial space}
\label{sec:unit-disc}

For our second example, we consider a classification problem for points sampled from two distributions concentrated on or close to the unit disk. 

First, we recall the (complex) Ginibre ensemble~\cite{Meckes:2020}.
Consider the random $N \times N$ matrix, $G_N$, whose entries  are sampled independently from the standard complex normal distribution, i.e. distributed according to $\frac{1}{\sqrt{2}}Z_1 + i\frac{1}{\sqrt{2}}Z_2$, where $Z_1,Z_2$ are independent standard normal variables.
Take its $N$ (complex) eigenvalues, scaled by $\frac{1}{\sqrt{N}}$.
Let $\mu_N$ denote the empirical measure on these $N$ points.
Then $\mu_N$ converges weakly almost surely to the uniform distribution on the unit disk as $N \to \infty$.

\begin{figure}[!htb]
    \centering
    \includegraphics[height=35mm]{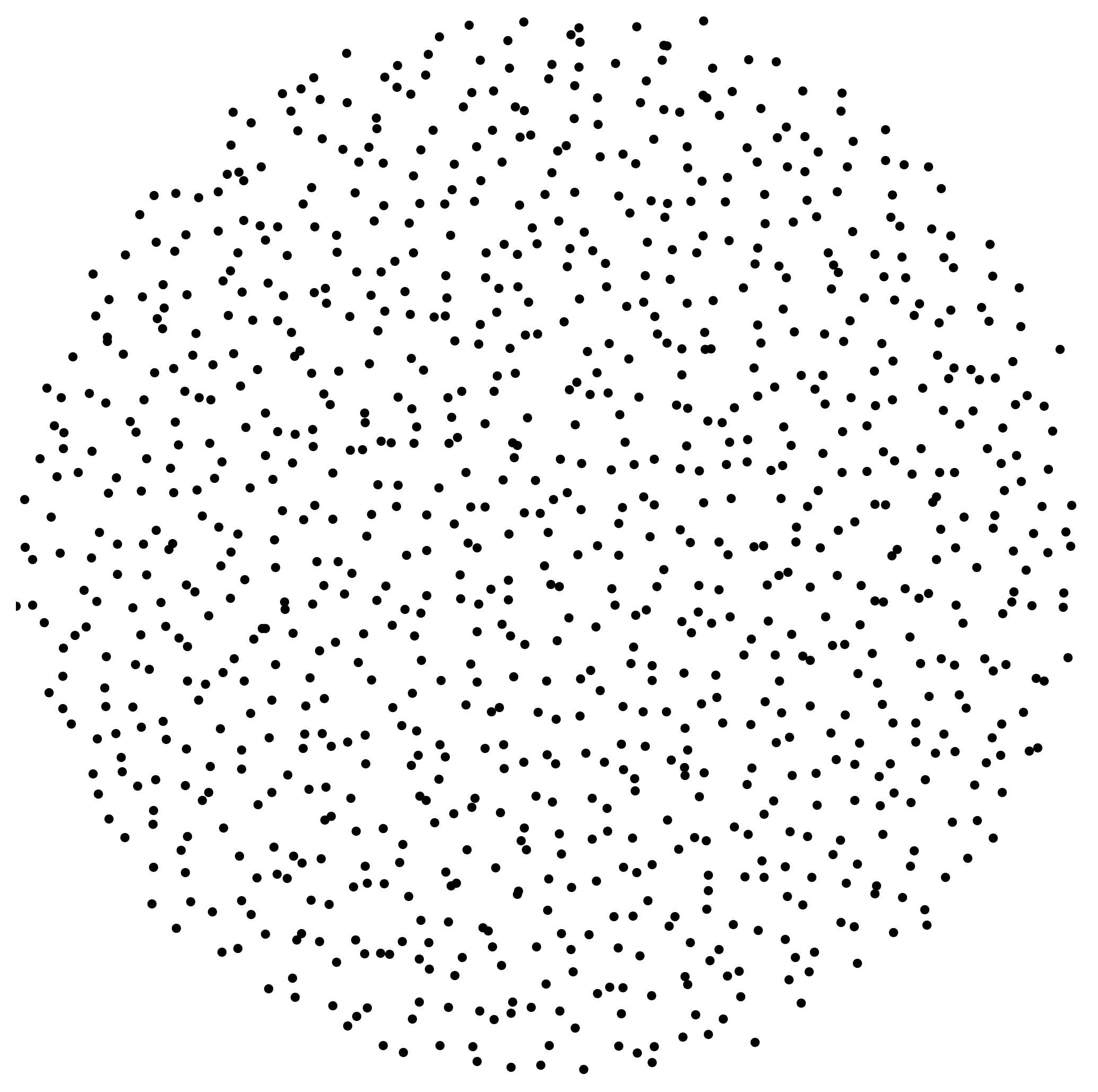}
    \quad \quad \quad
    \includegraphics[height=35mm]{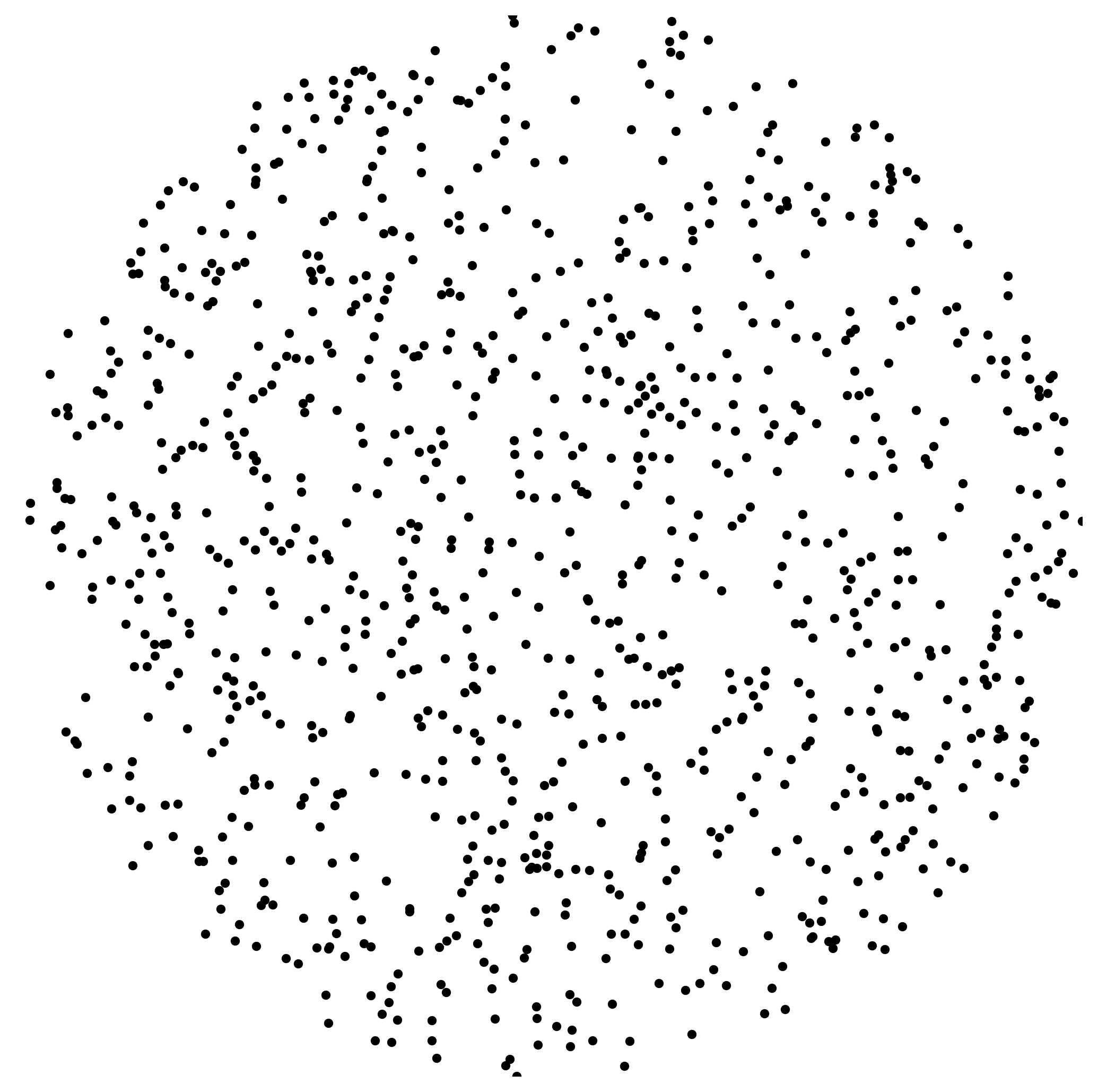}
    \caption{
    Left: scaled complex eigenvalues of the complex Ginibre ensemble.
    Right: Points sampled independently from the uniform distribution on the unit disk. 
    }
    \label{fig:ginibre}
\end{figure}

Let $N=1000$.
For the first distribution, 
we take the $N$ (complex) eigenvalues of $G_N$ scaled by $\frac{1}{\sqrt{N}}$. 
For the second distribution, 
we take $N$ points sampled independently from the uniform distribution on the unit disk.
See \cref{fig:ginibre}.

\begin{figure}[!htbp]
    \centering
    \includegraphics[height=35mm]{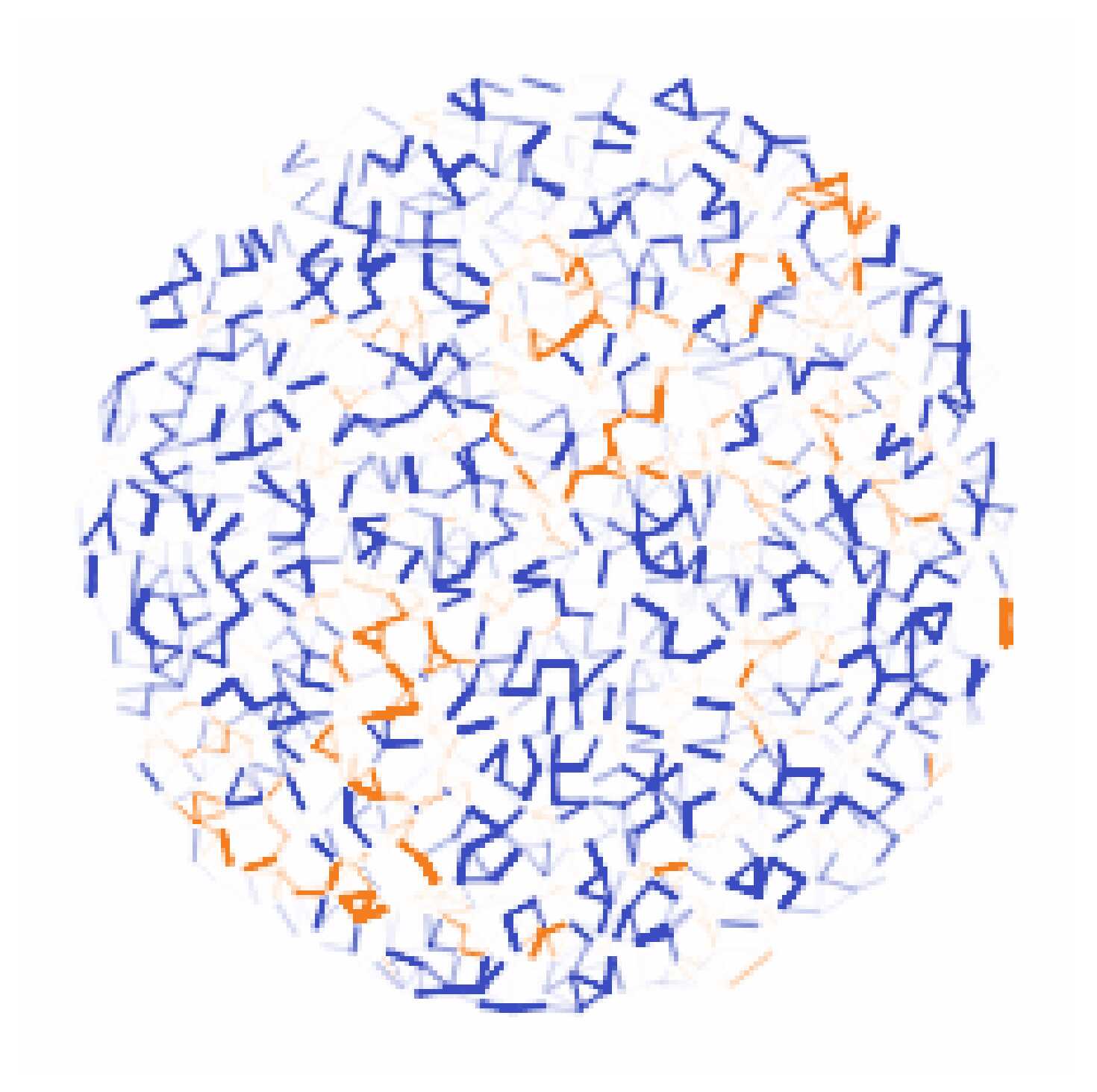}
    \includegraphics[height=35mm]{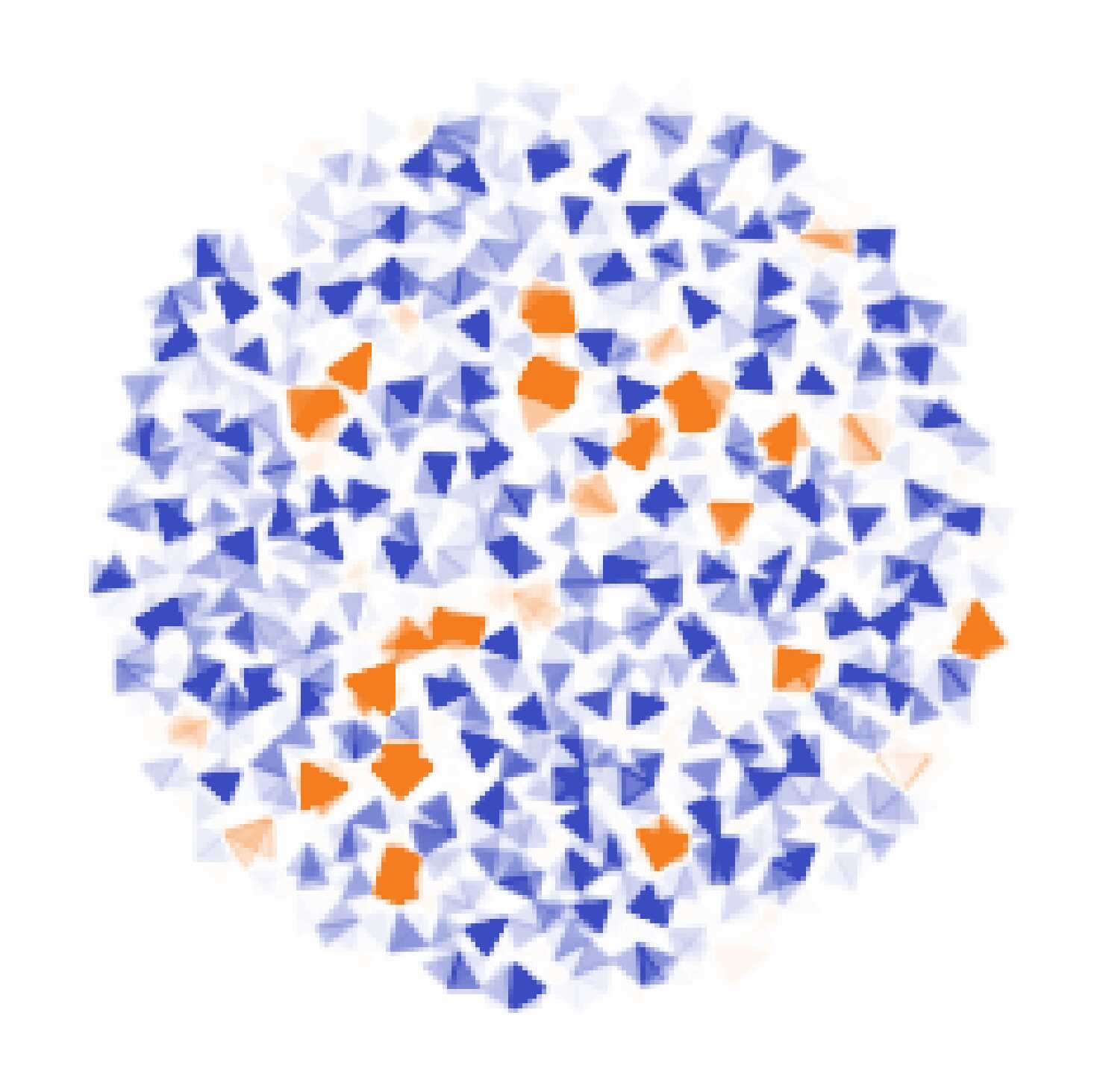}
    \includegraphics[height=35mm]{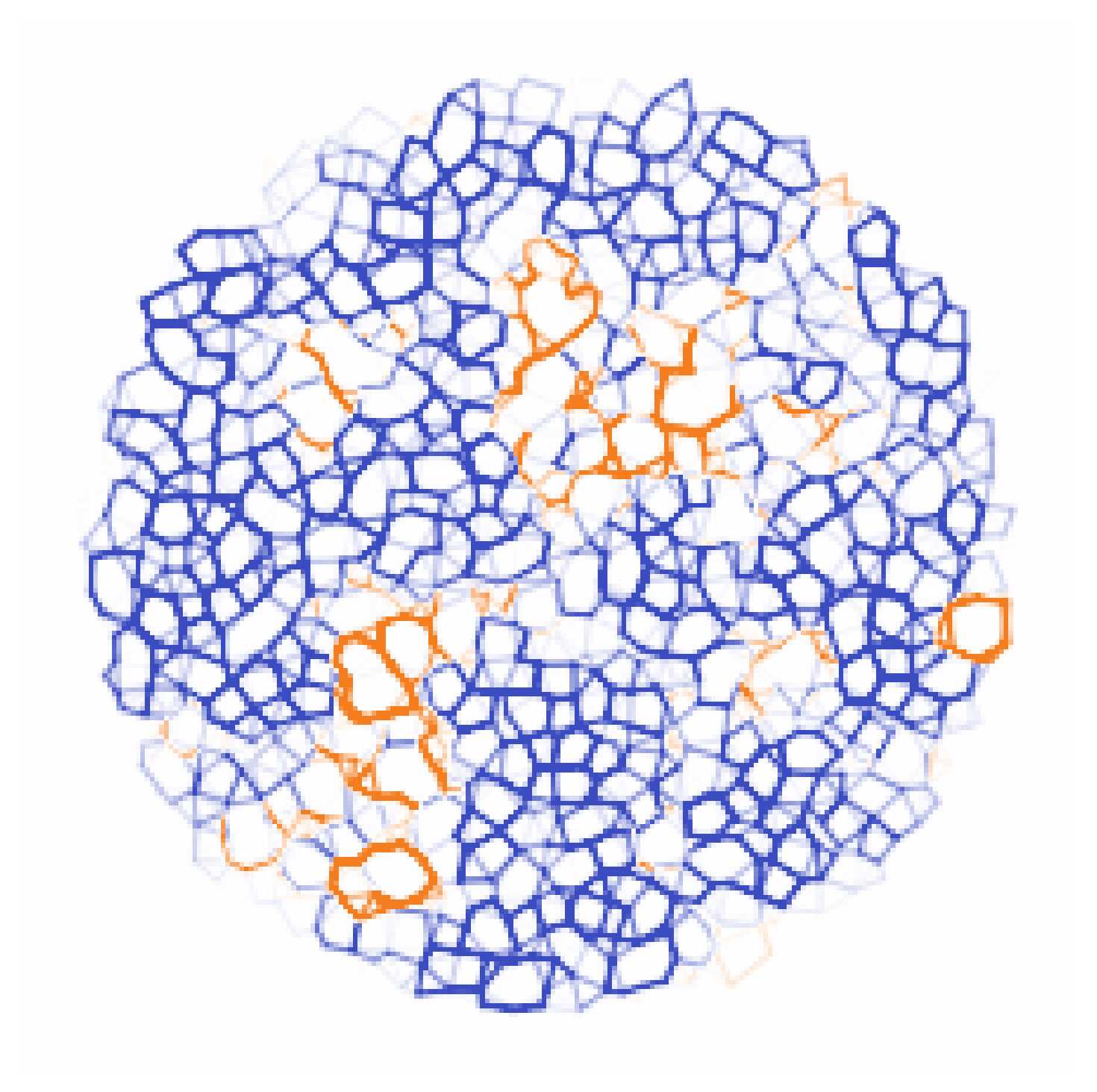}
    \includegraphics[height=35mm]{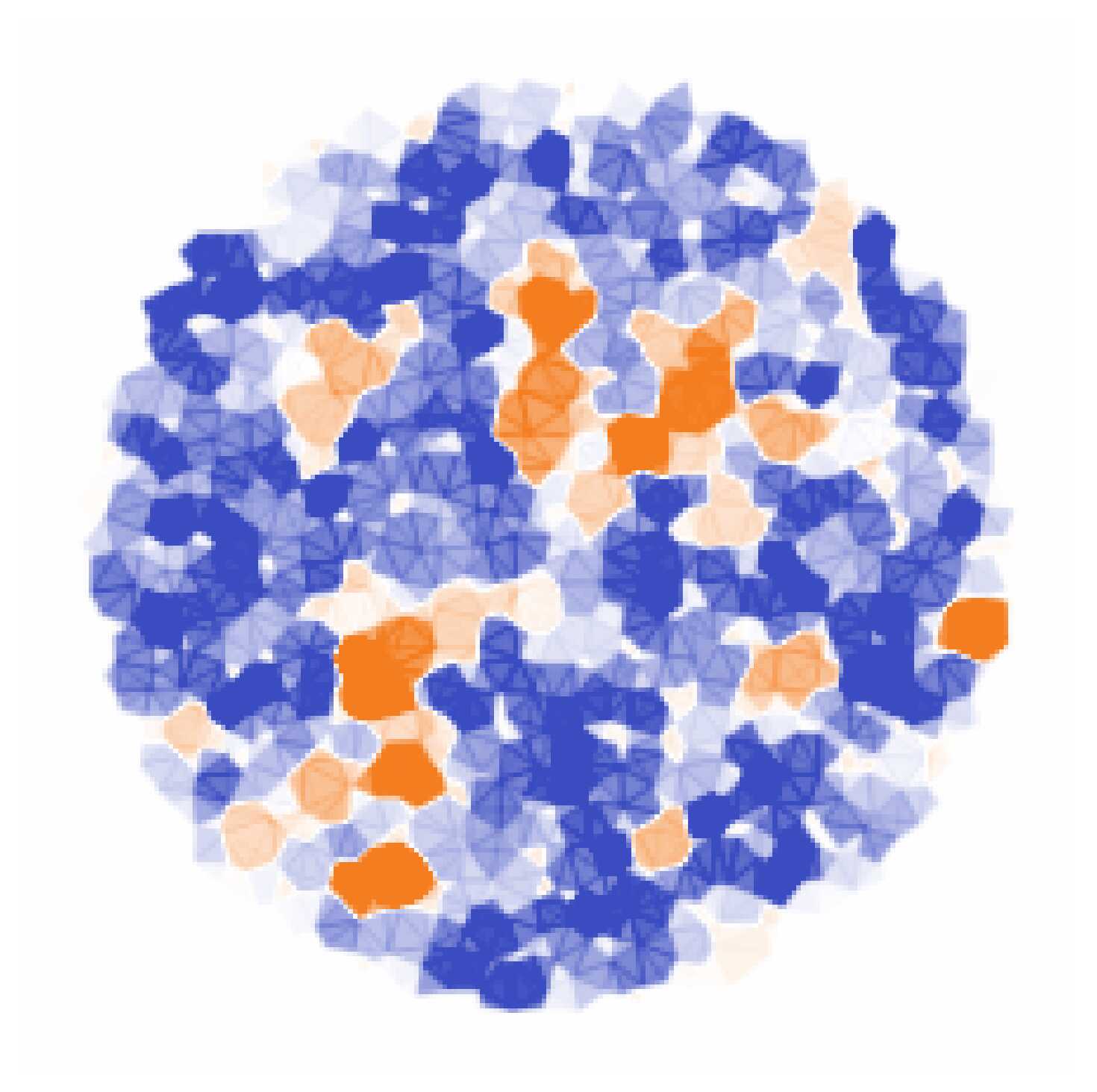}\\
    \includegraphics[height=35mm]{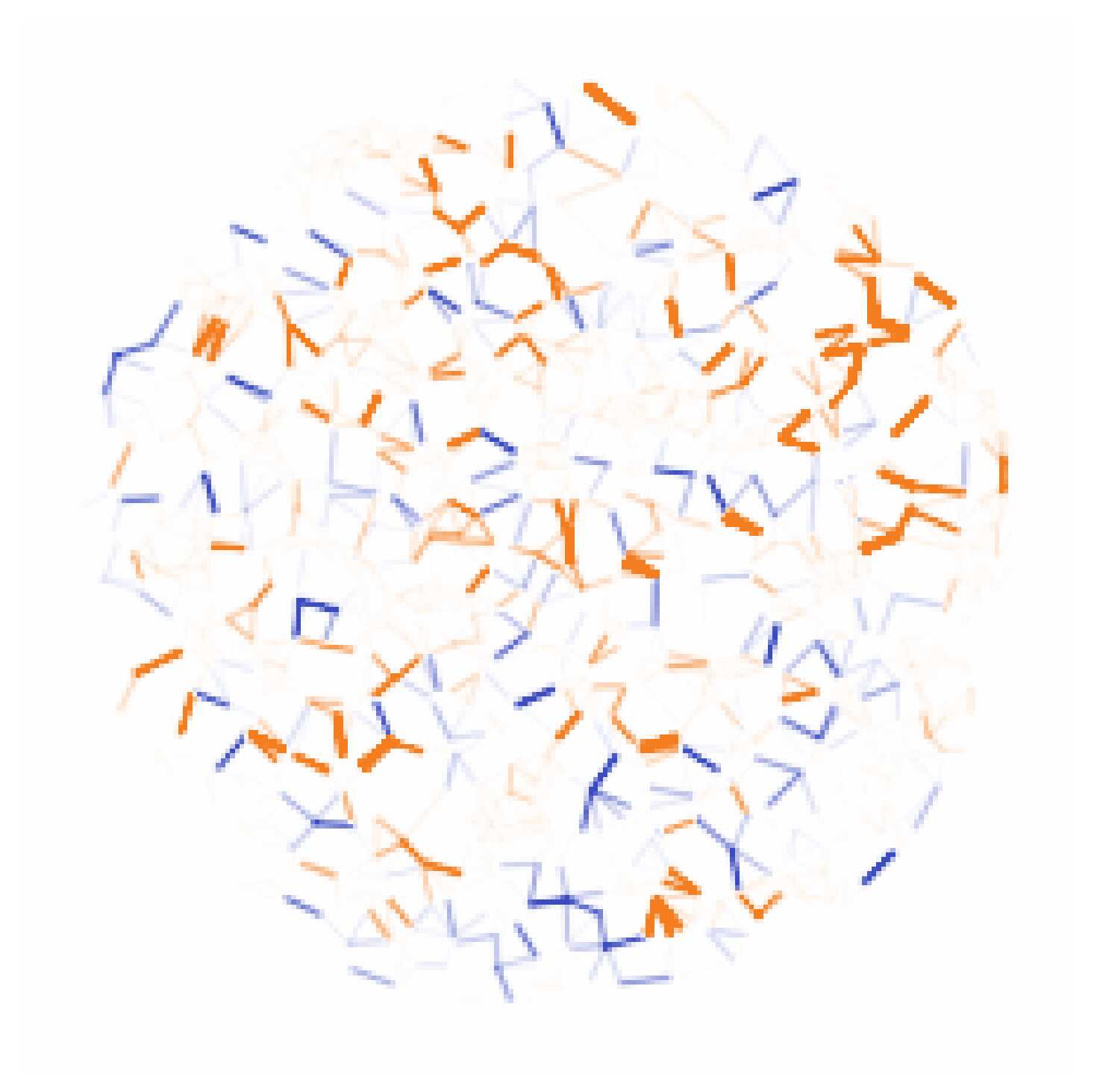}
    \includegraphics[height=35mm]{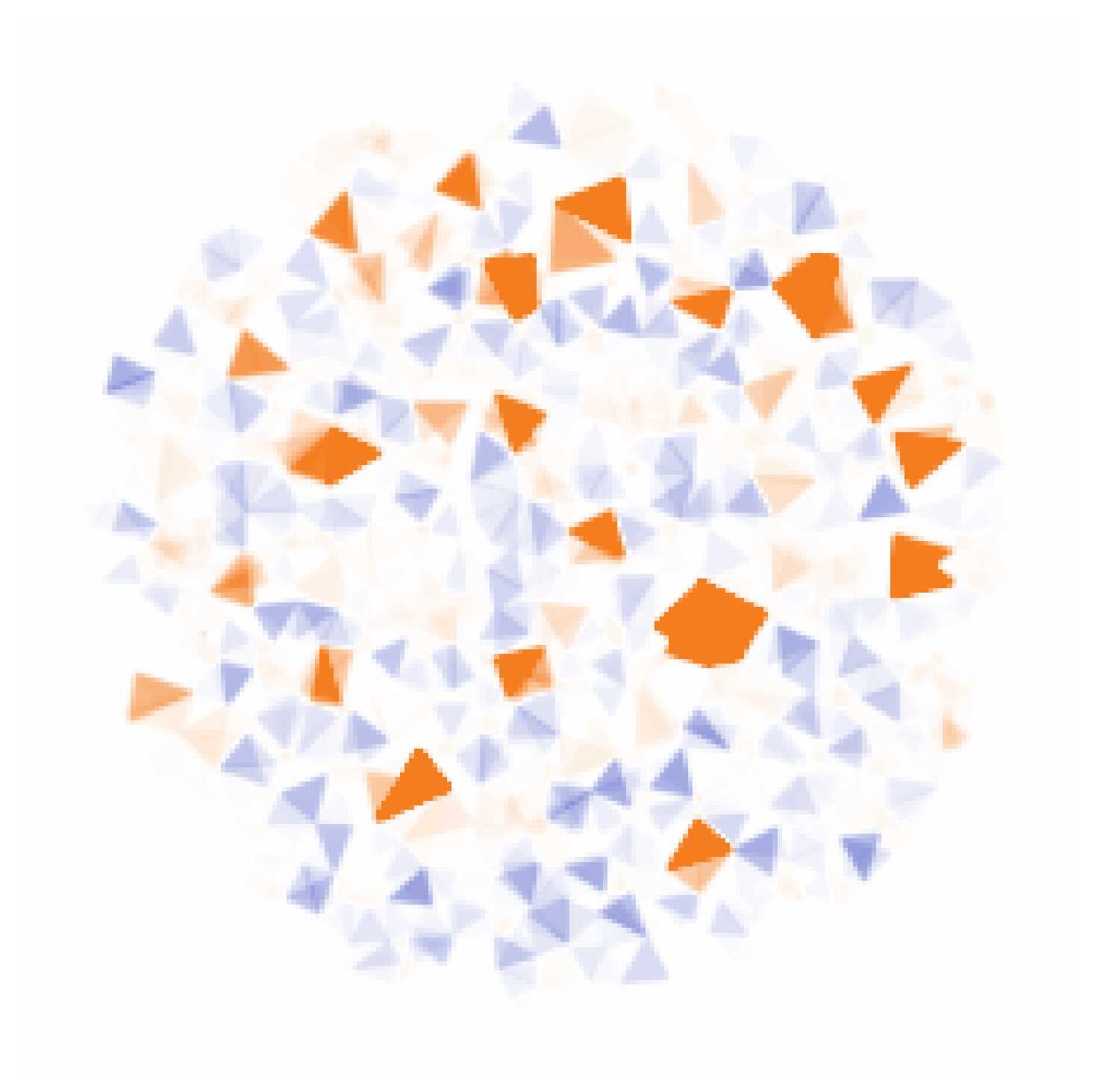}
    \includegraphics[height=35mm]{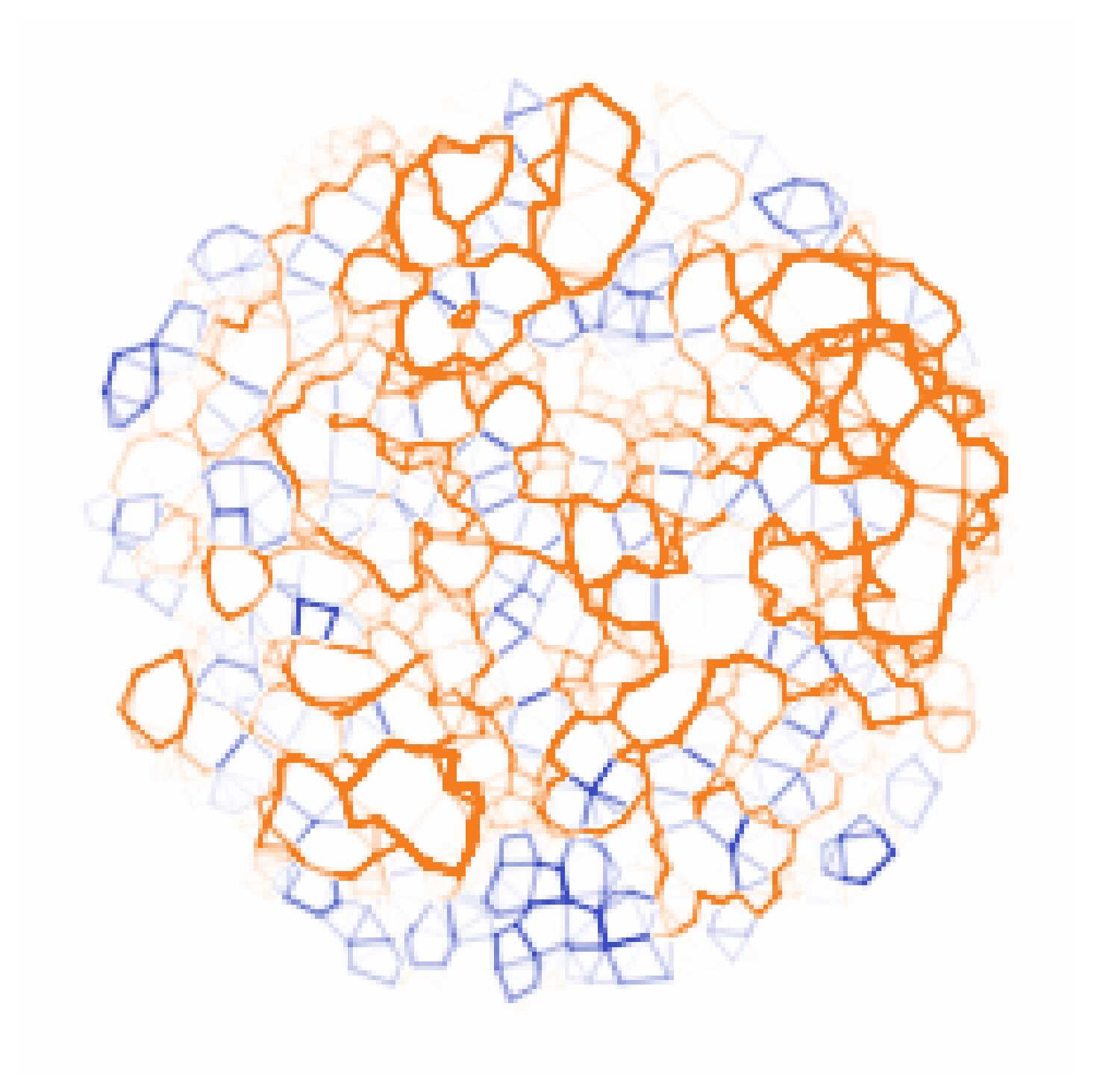}
    \includegraphics[height=35mm]{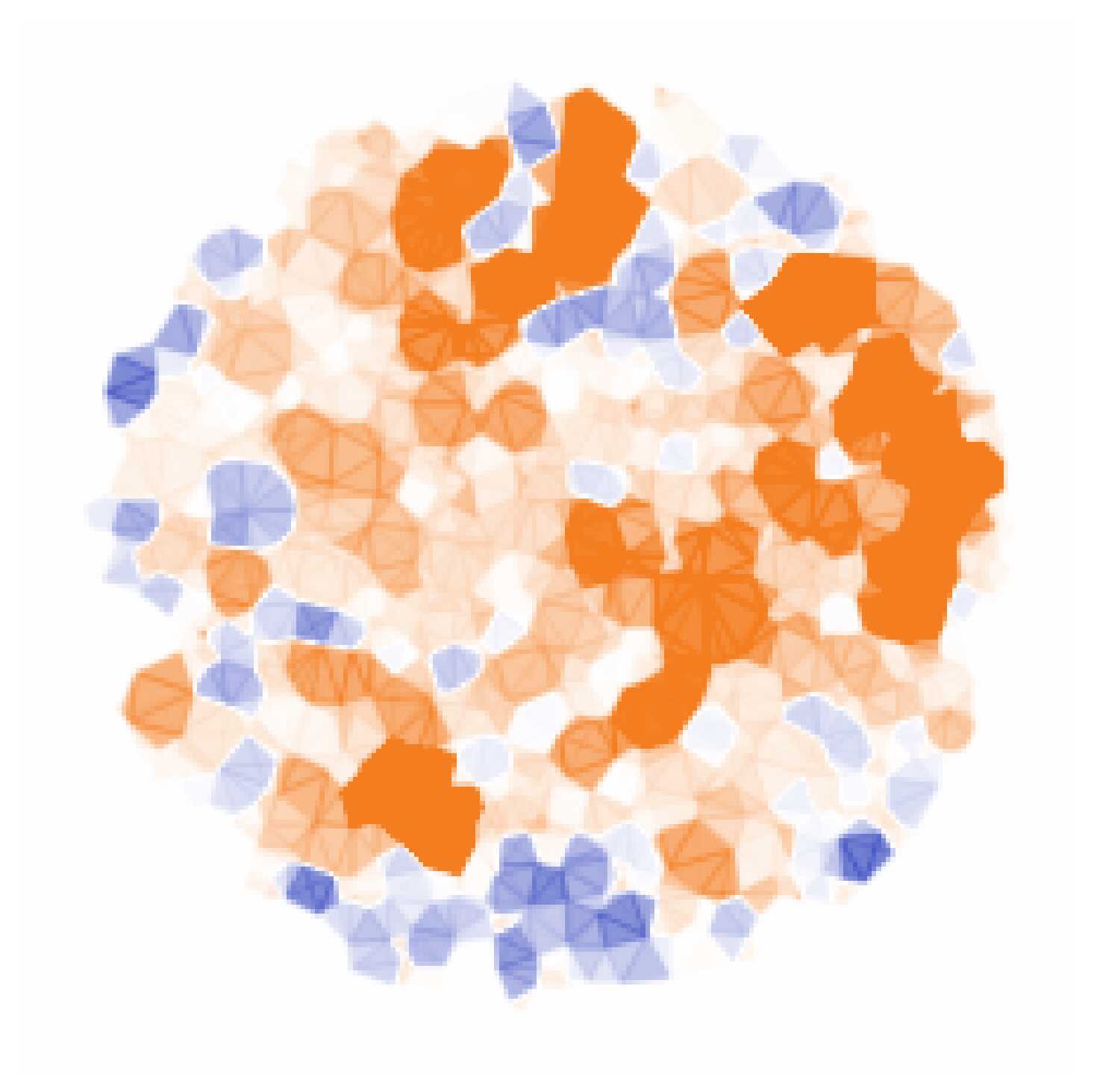}
    \caption{
    Visualization of the learned classifier that differentiates between 
    eigenvalues of the Ginibre ensemble (top) and
    points sampled from the uniform distribution on a disk (bottom), 
    using homology in degree $1$ and labeled persistence landscapes. 
    We use a orange-to-blue color gradient in which orange corresponds to the uniform distribution and blue corresponds to the Ginibre ensemble.
    From left to right we have birth vertices, death edges, representative cycles and bounding chains, respectively.
    In each image, the sum of the values of the pixels equals the empirical mean value of the classifier.
    }
    \label{fig:ginibre-visualization-deg1-single}
\end{figure}

As in the previous example, we take 100 point clouds from each class.
For each of these we consider its corresponding Delaunay complex.
We compute the persistence diagrams for homology 
in degree $1$, labeled by birth simplices, death simplices, representative cycles, and bounding chains. 
For the degree-$1$ labeled persistence diagram, we consider its labeled persistence landscape. 
In each of these cases, we apply SVM and use the normal vector and bias to map each feature vector to a real number.
As in the previous example take one sample from each class and compute the empirical mean feature map under 10000 perturbations.
Each perturbation adds independent Gaussian noise, with mean $0$ and standard deviation $\sigma = 0.005$,
to each Cartesian coordinate of every point in the cloud.
See 
\cref{fig:ginibre-visualization-deg1-single}.

Note that the points in the Ginibre ensemble are distributed more uniformly (since eigenvalues `repel' each other) and the uniform distribution has larger `voids'. It is these voids which are most strongly colored with the color corresponding to the uniform distribution.




\subsection{Classification for 3-d data}

Here we consider the classification problem for two classes of point clouds in $\R^3$.

The first class is a \emph{cup}: the union of an (open) cylinder and a bottom disk
(and hence open at the top). 
With axis radius $R=1$ and height $h=2$, it is
\[
    \{(R\cos\theta,\,R\sin\theta,\,z) : \theta\in[0,2\pi),\ z\in[0,h]\}
    \;\cup\;
%
    \{(\rho\cos\theta,\,\rho\sin\theta,\,0) : \rho\in[0,R],\ \theta\in[0,2\pi)\}
\]

The second class is a \emph{cup with a handle}: the same cup together with a flat
half-annulus handle lying in the plane $y=0$, centered at $(R,0,h/2)$ and opening
in the positive $x$ direction,
\[
    \{(R + r\cos\alpha,\;0,\;h/2 + r\sin\alpha) : r\in[r_{\mathrm{in}}, r_{\mathrm{out}}],\ \alpha\in[-\tfrac{\pi}{2},\tfrac{\pi}{2}]\},
\]
with $r_{\mathrm{in}}=0.5$ and $r_{\mathrm{out}}=1.0$.
The handle attaches to the cup along the line $x=R$, $y=0$ and bulges outward in $+x$;
the loop running through the handle and back along the cup wall contributes a degree-$1$ homology class, which is the only topological feature distinguishing the
two classes.

\begin{figure}[!htb]
    \centering
    \includegraphics[height=35mm]{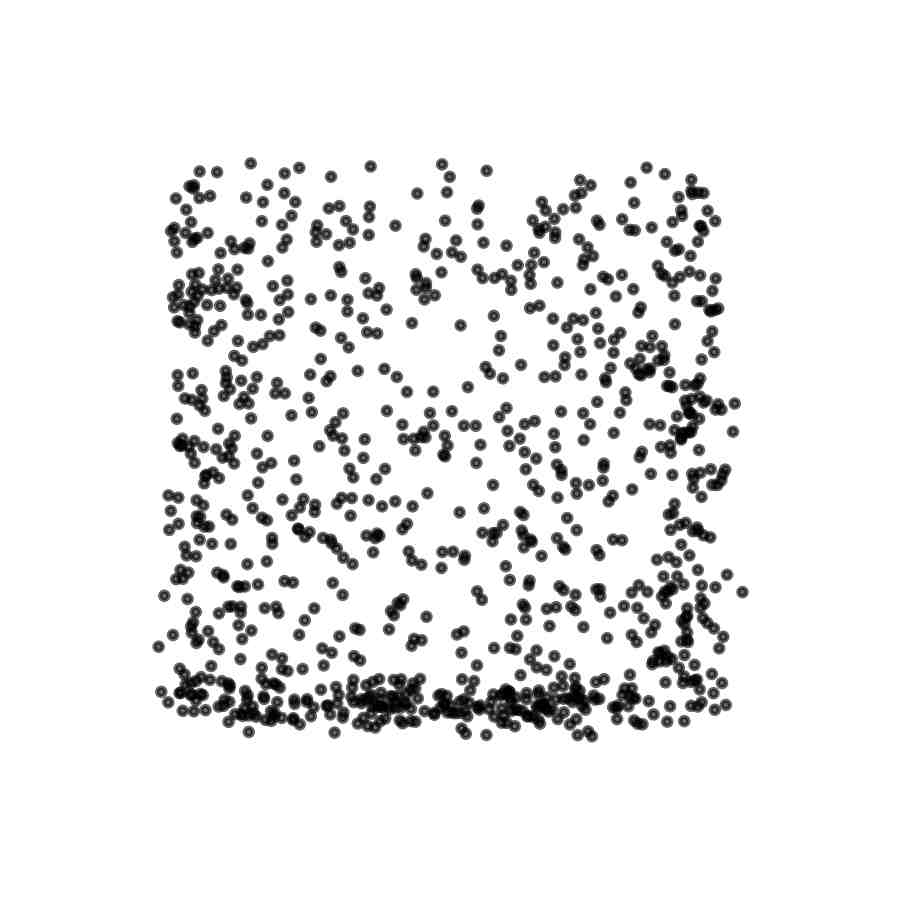}
    \quad \quad
    \includegraphics[height=35mm]{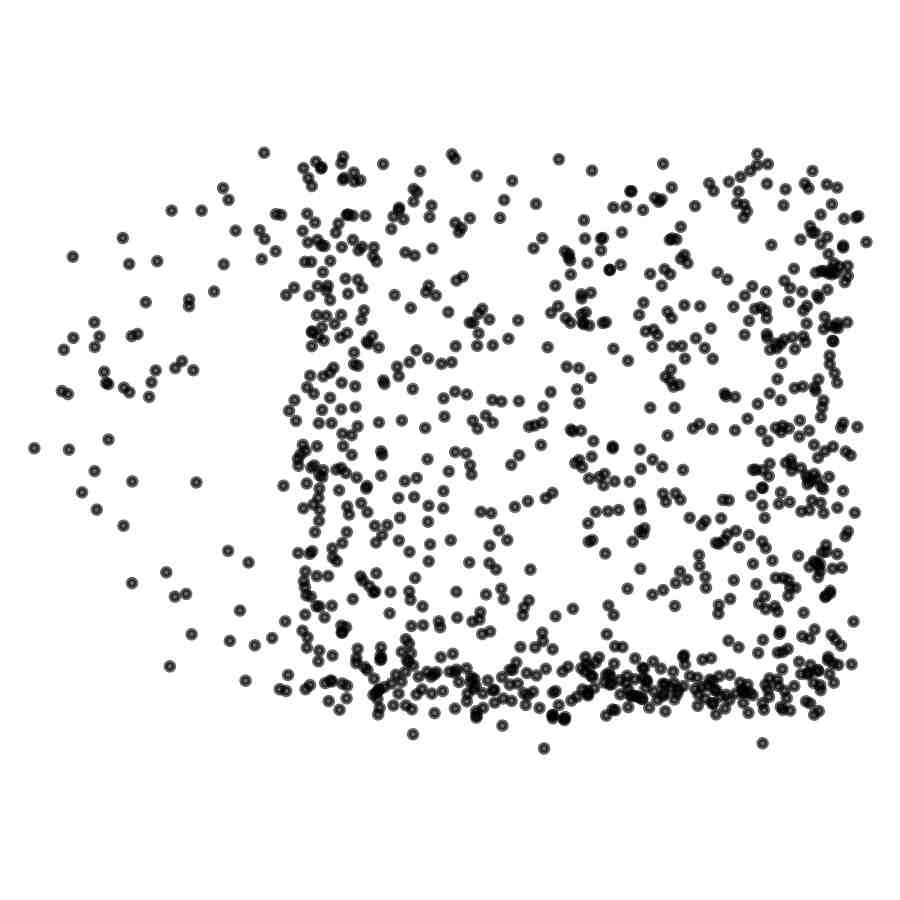}
    \caption{
    $1000$ points sampled uniformly from the surface of a cup (left) and a cup with a
    handle (right), each perturbed by isotropic Gaussian noise with $\sigma=0.05$.
    }
    \label{fig:cup-handle}
\end{figure}

We sample $1000$ points from each surface using the uniform probability measure with
respect to surface area, and then perturb every point by independent isotropic Gaussian
noise $\mathcal{N}(0,\sigma^2 I_3)$ with $\sigma=0.05$.
See \cref{fig:cup-handle}.

\begin{remark}
    To sample uniformly with respect to surface area, we (i) allocate the $1000$ points
    to the constituent pieces in proportion to their areas and (ii) sample each piece
    uniformly.
    The cylinder has area $2\pi R h$ and is sampled by
    $\theta \sim \mathrm{Uniform}[0,2\pi)$, $z \sim \mathrm{Uniform}[0,h]$.
    The bottom disk has area $\pi R^2$ and is sampled by
    $\theta \sim \mathrm{Uniform}[0,2\pi)$, $\rho = R\sqrt{U}$ with $U \sim \mathrm{Uniform}[0,1]$.
    The half-annulus handle has area $\frac{\pi}{2}\bigl(r_{\mathrm{out}}^2 - r_{\mathrm{in}}^2\bigr)$
    and is sampled by
    $\alpha \sim \mathrm{Uniform}[-\tfrac{\pi}{2},\tfrac{\pi}{2}]$,
    $r = \sqrt{r_{\mathrm{in}}^2 + U\,(r_{\mathrm{out}}^2 - r_{\mathrm{in}}^2)}$.
\end{remark}

\begin{figure}[!htbp]
    \centering
    \includegraphics[height=17mm]{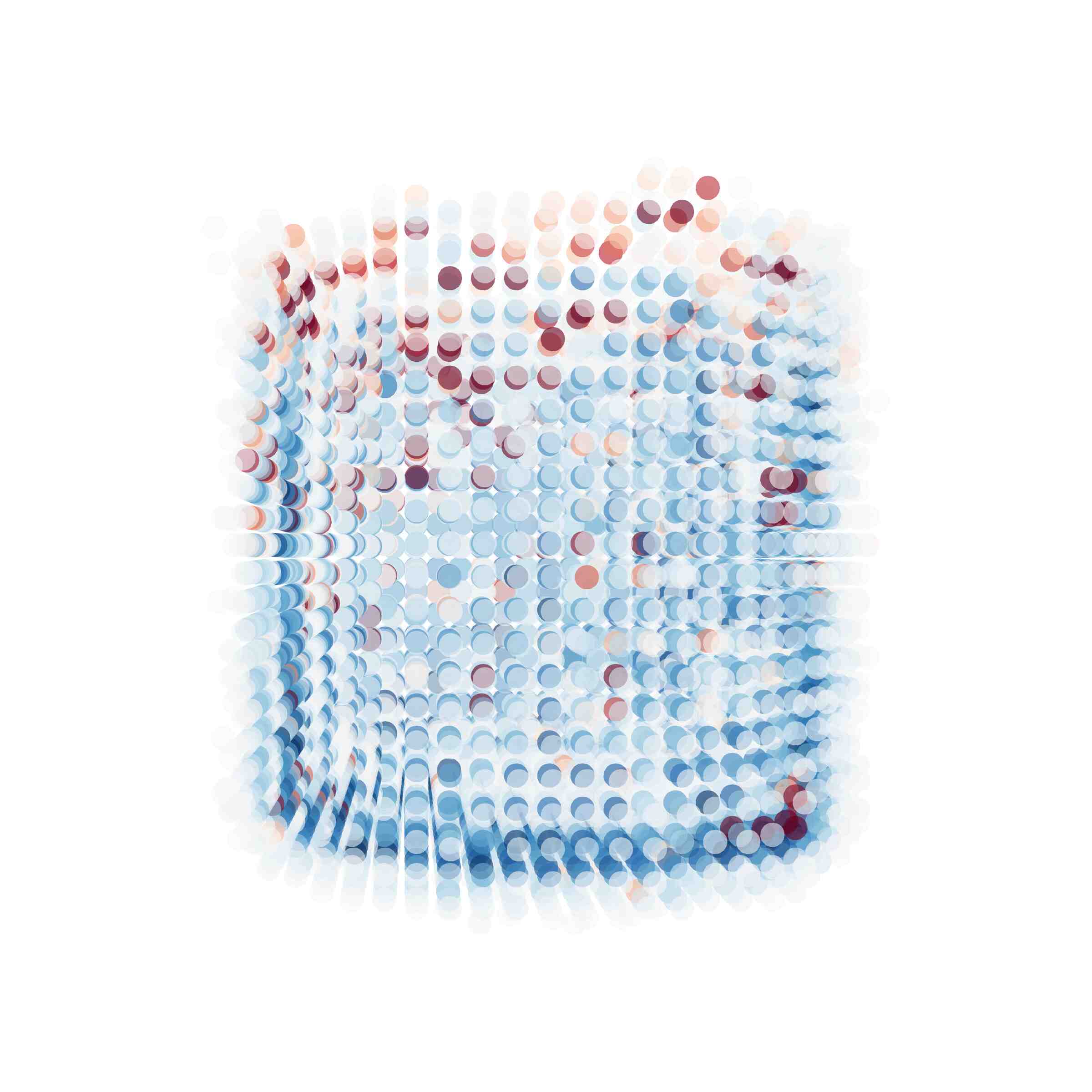}
    \includegraphics[height=17mm]{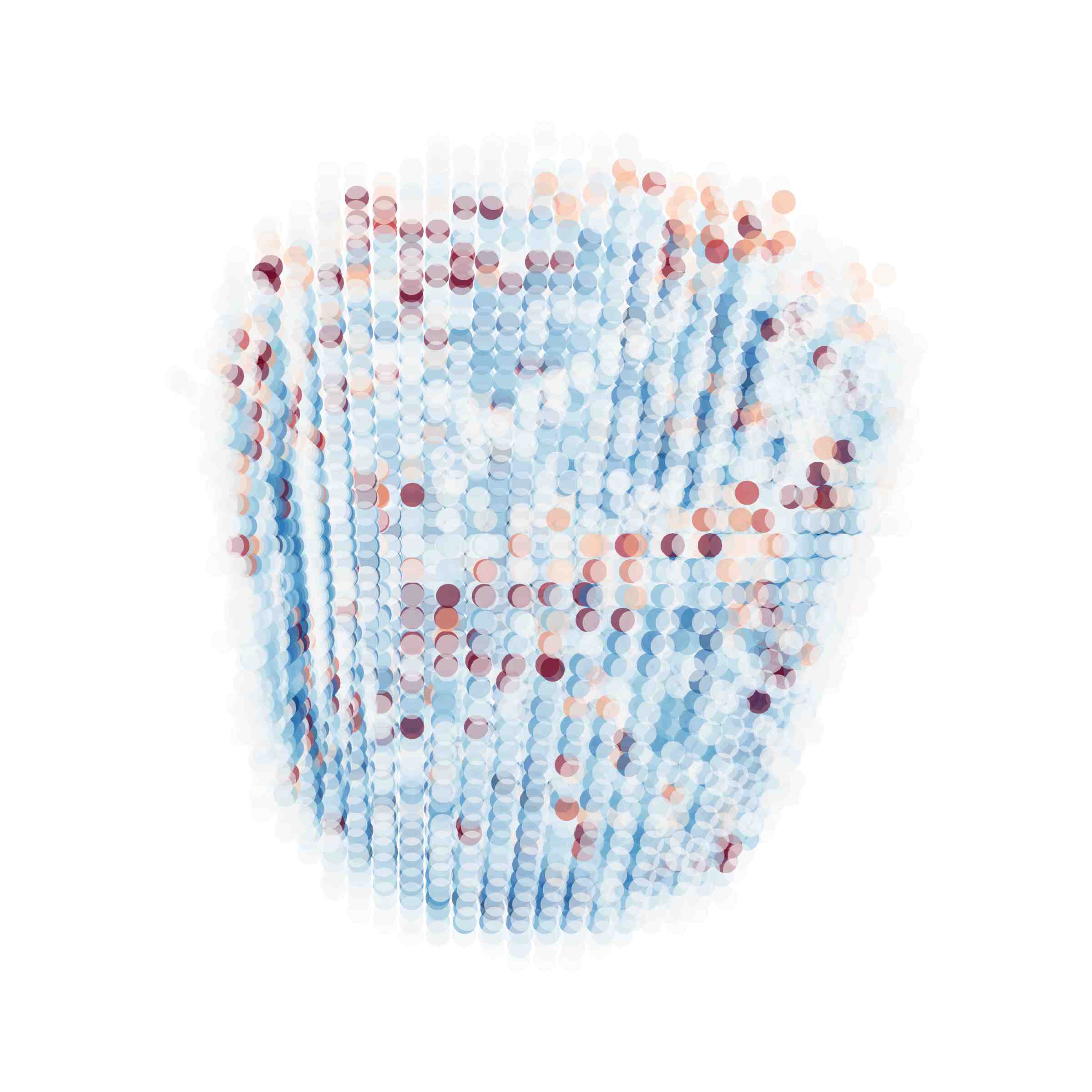}
    \includegraphics[height=17mm]{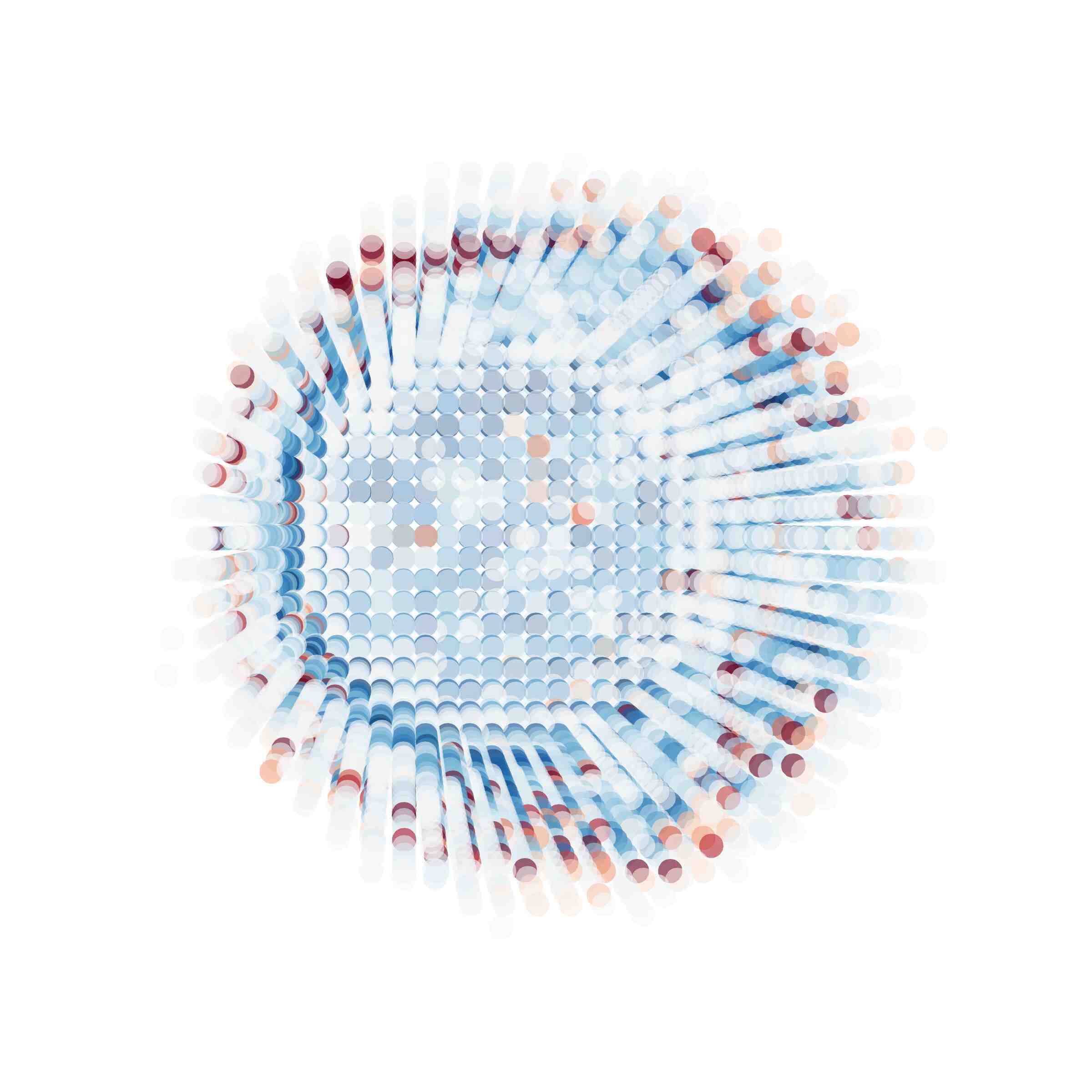}
    \includegraphics[height=17mm]{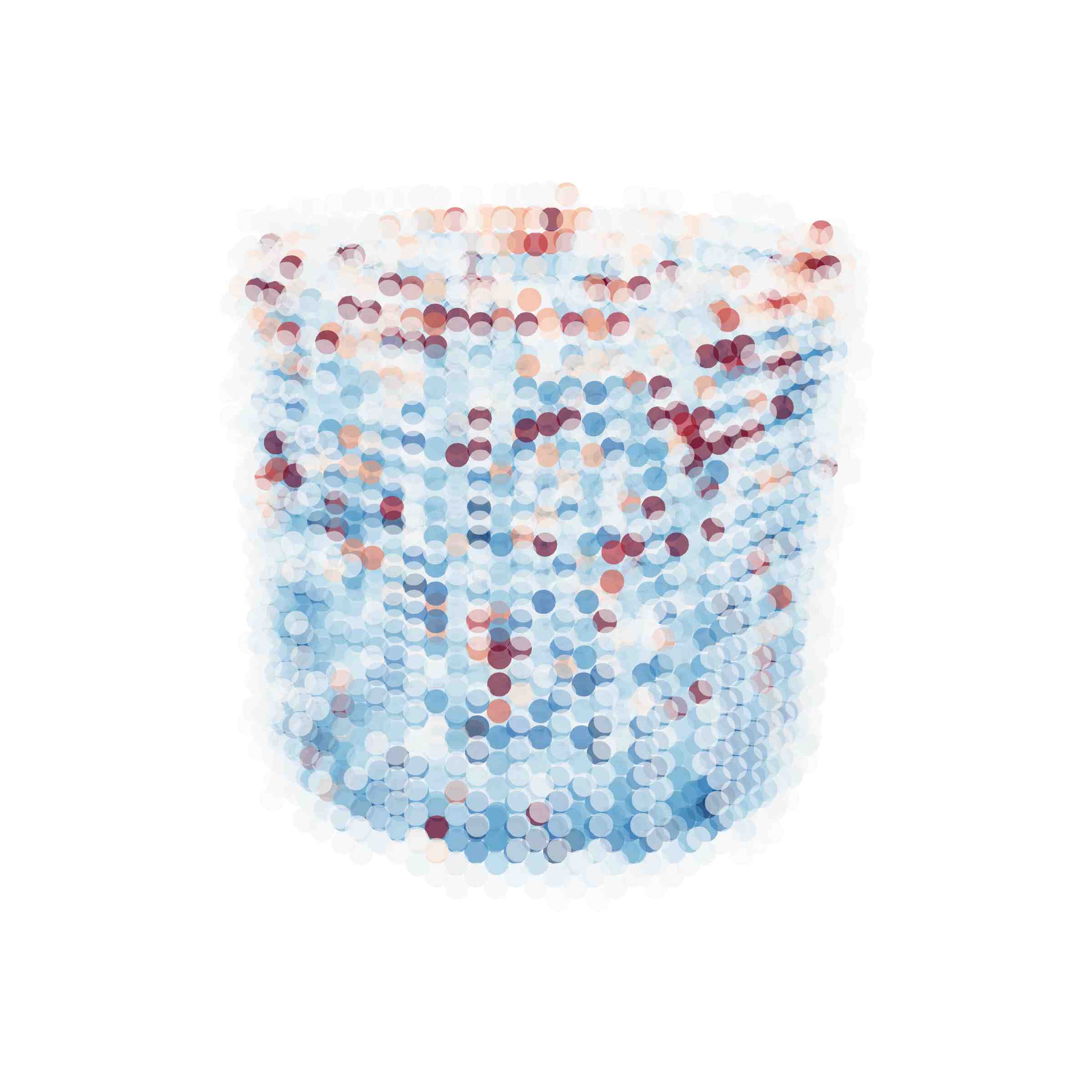}
    \includegraphics[height=17mm]{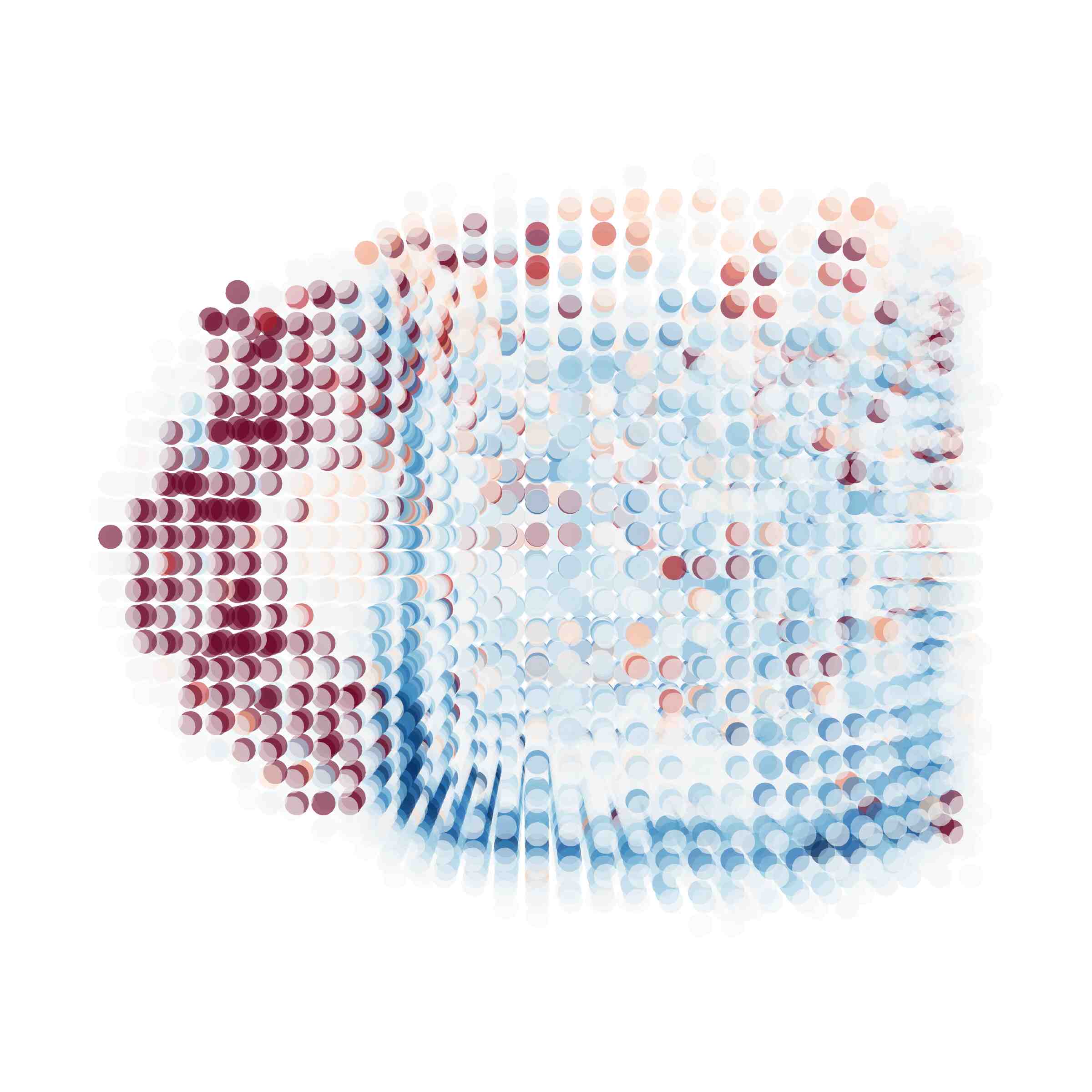}
    \includegraphics[height=17mm]{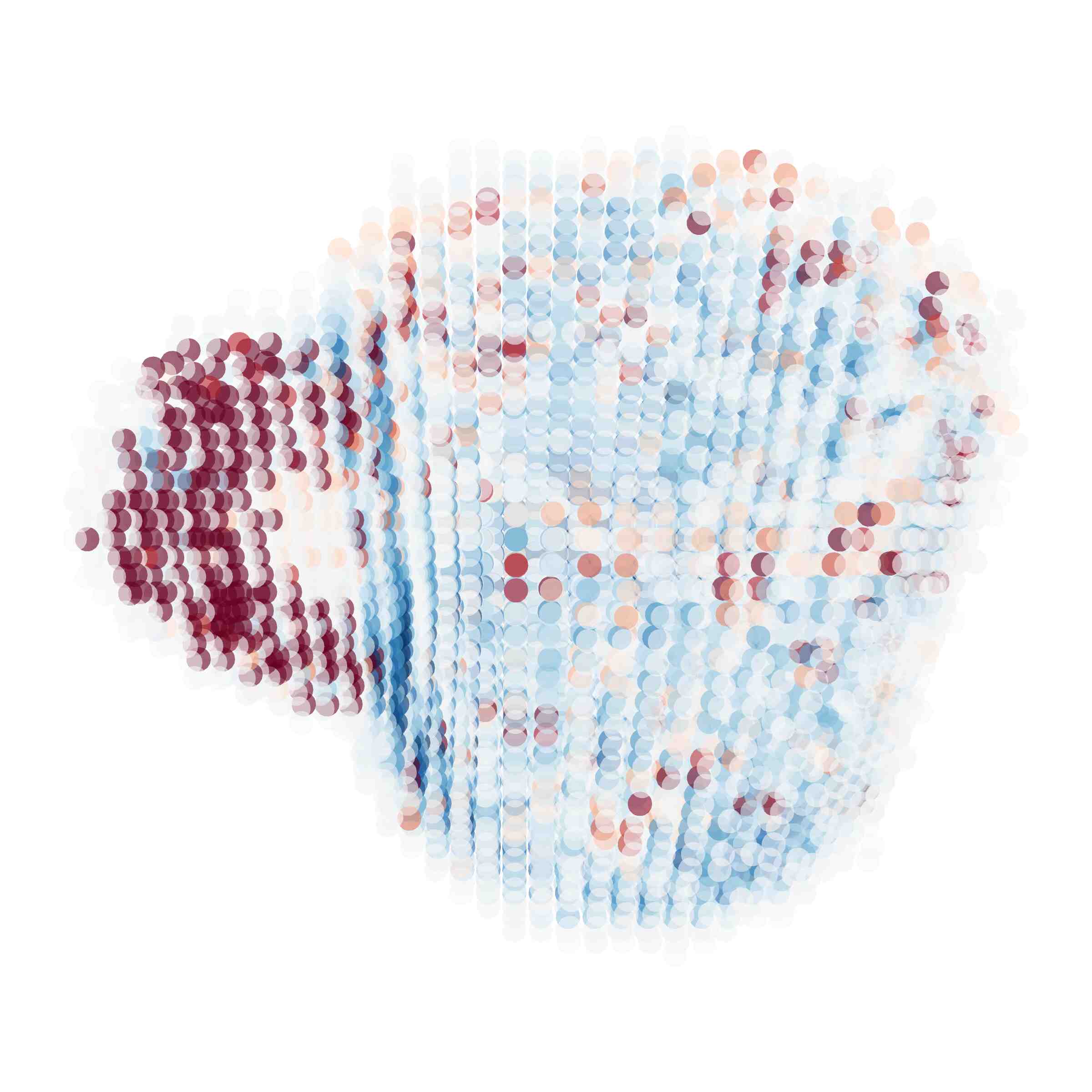}
    \includegraphics[height=17mm]{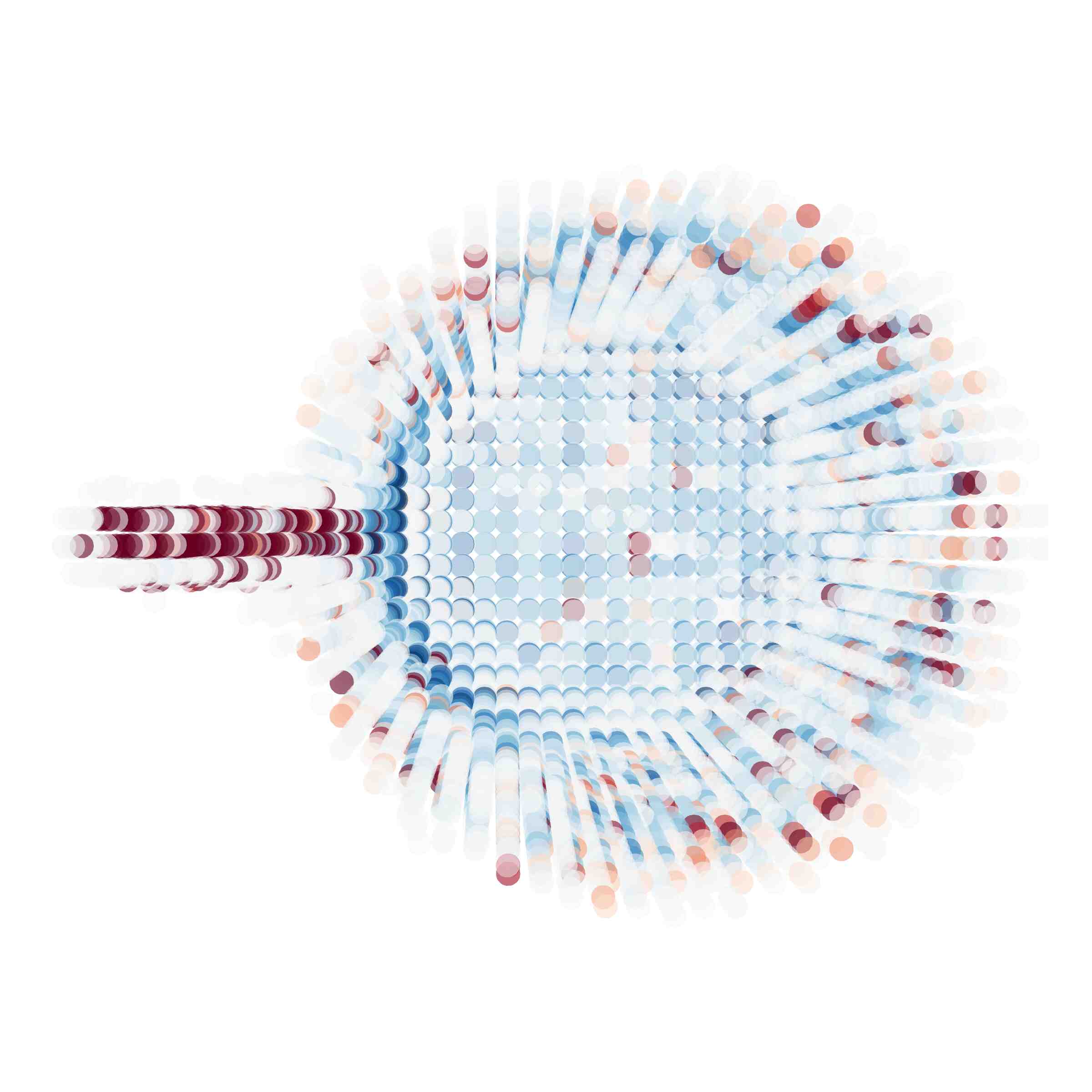}
    \includegraphics[height=17mm]{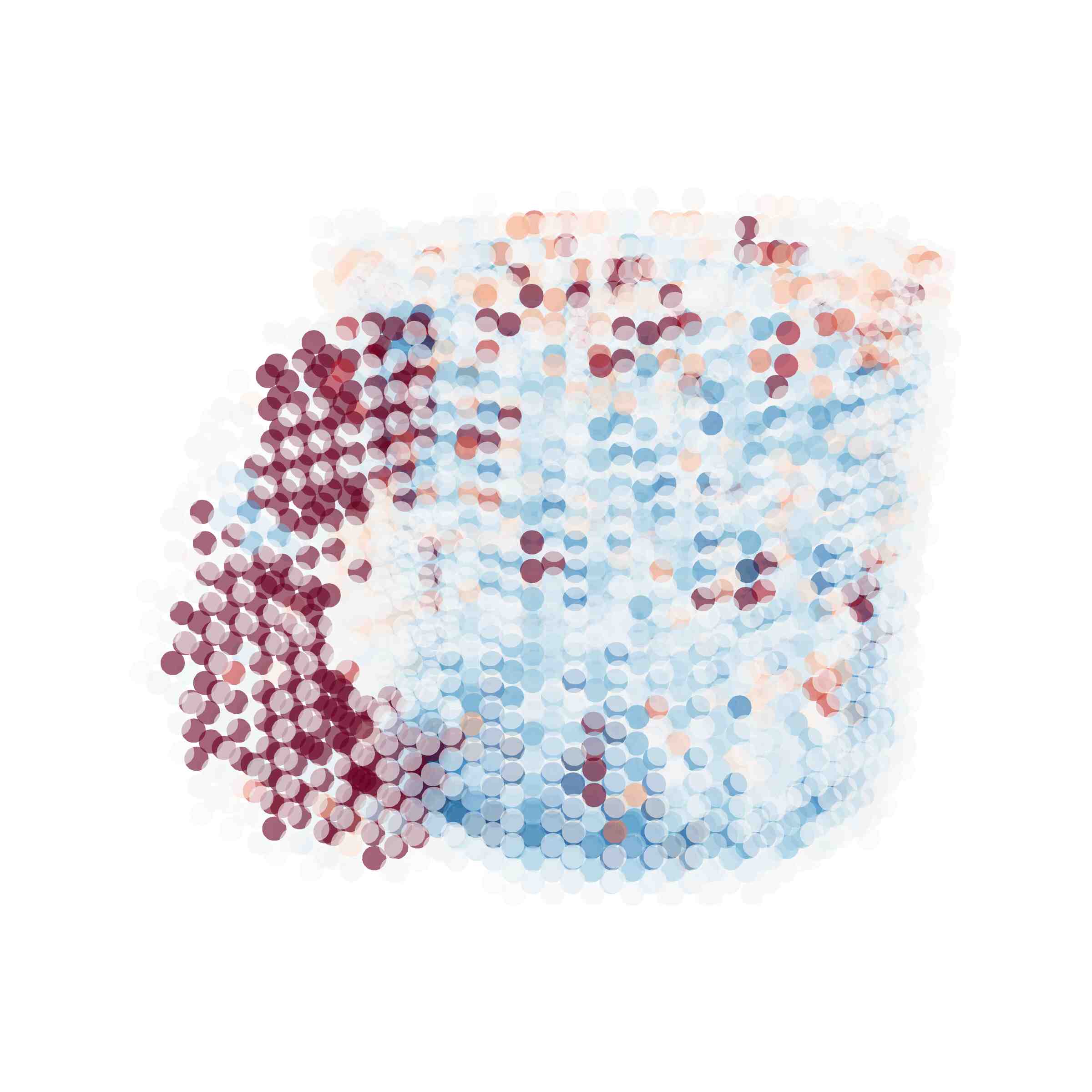}

    \includegraphics[height=17mm]{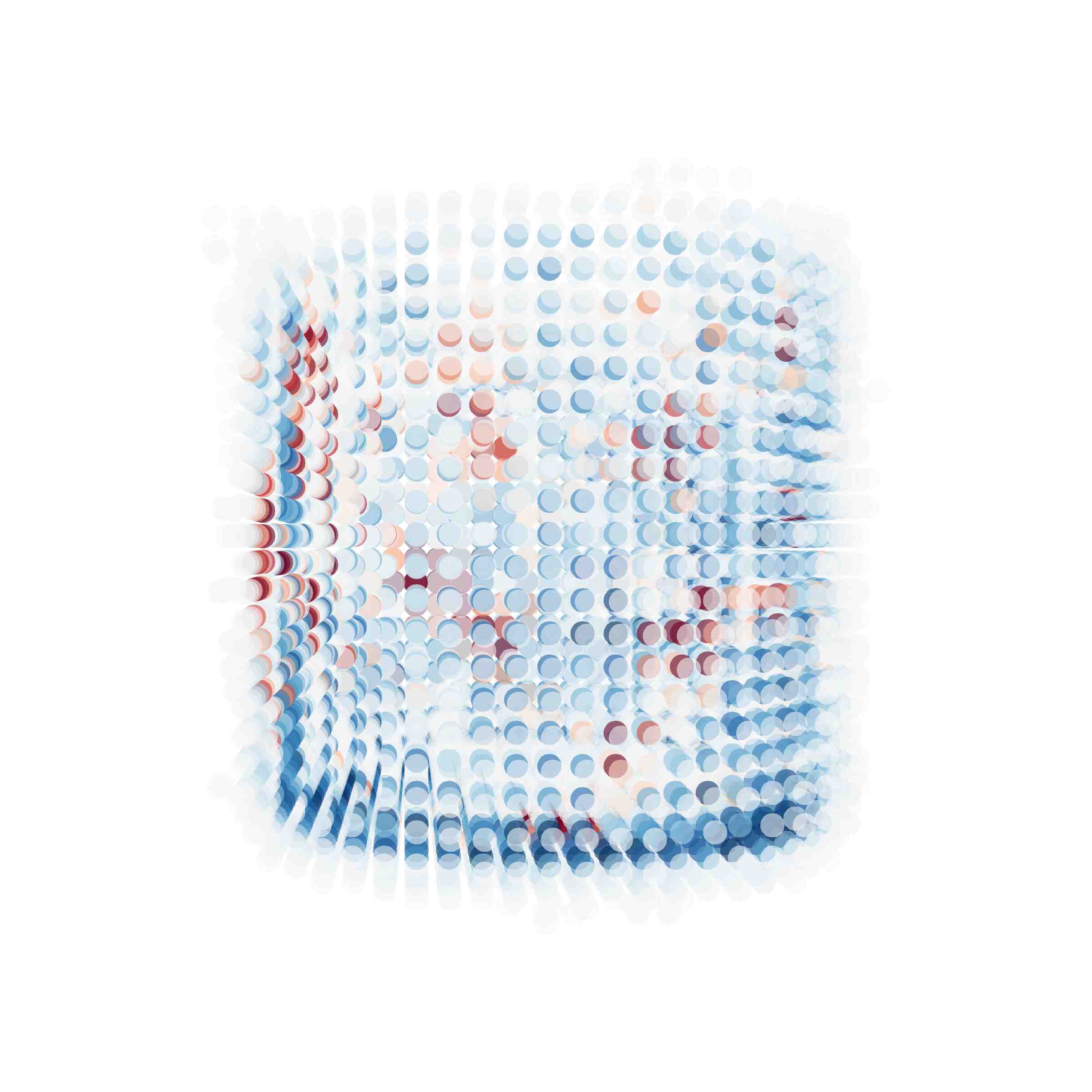}
    \includegraphics[height=17mm]{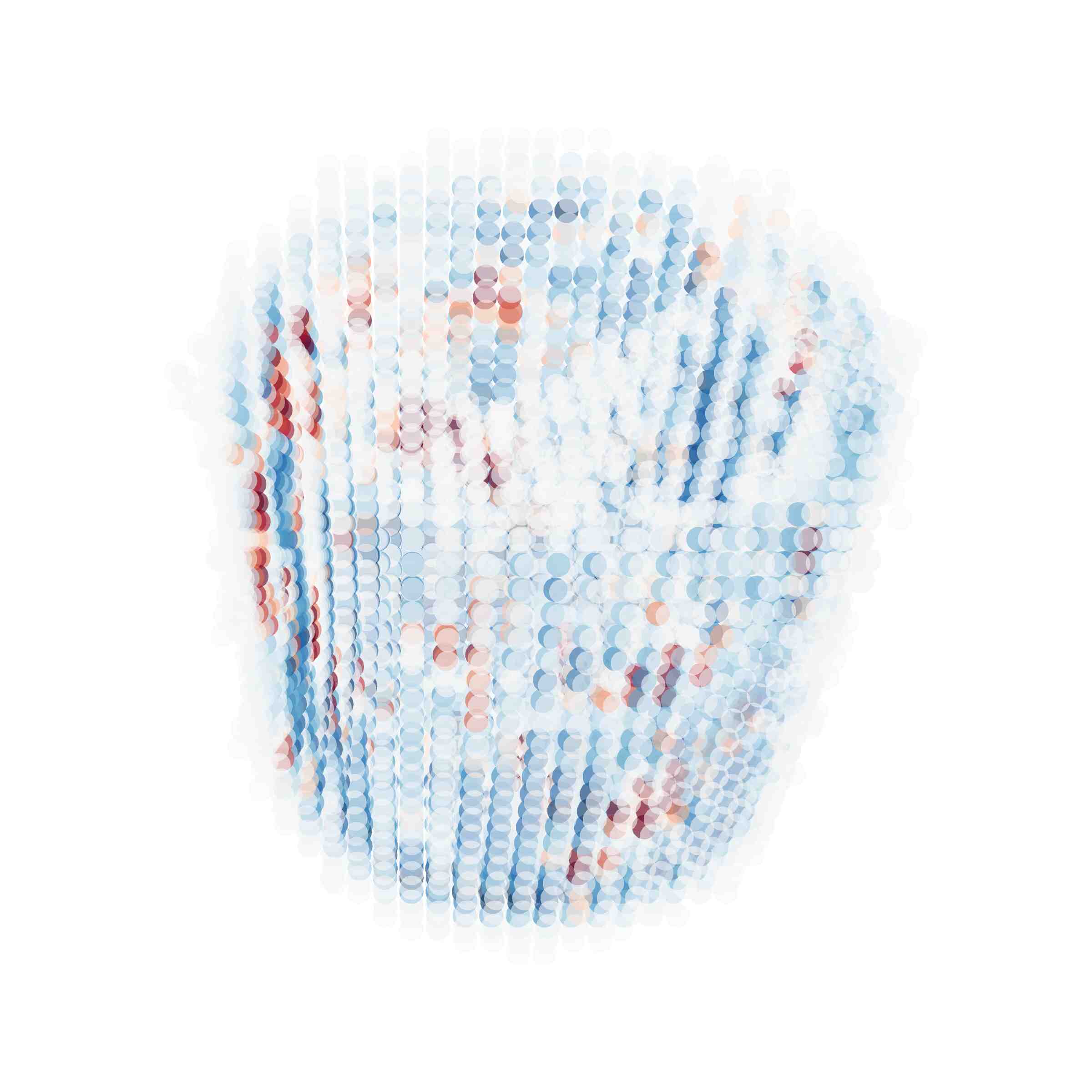}
    \includegraphics[height=17mm]{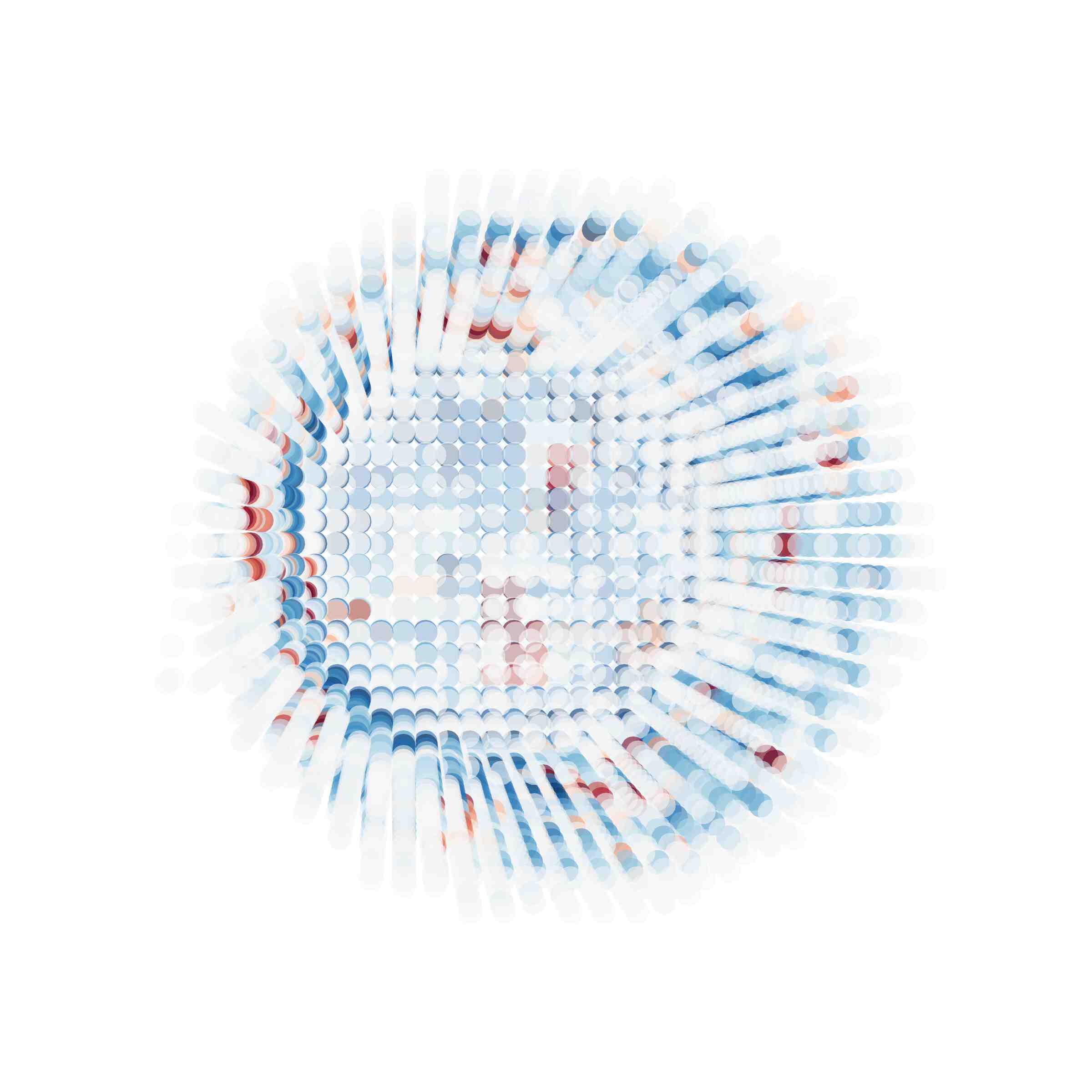}
    \includegraphics[height=17mm]{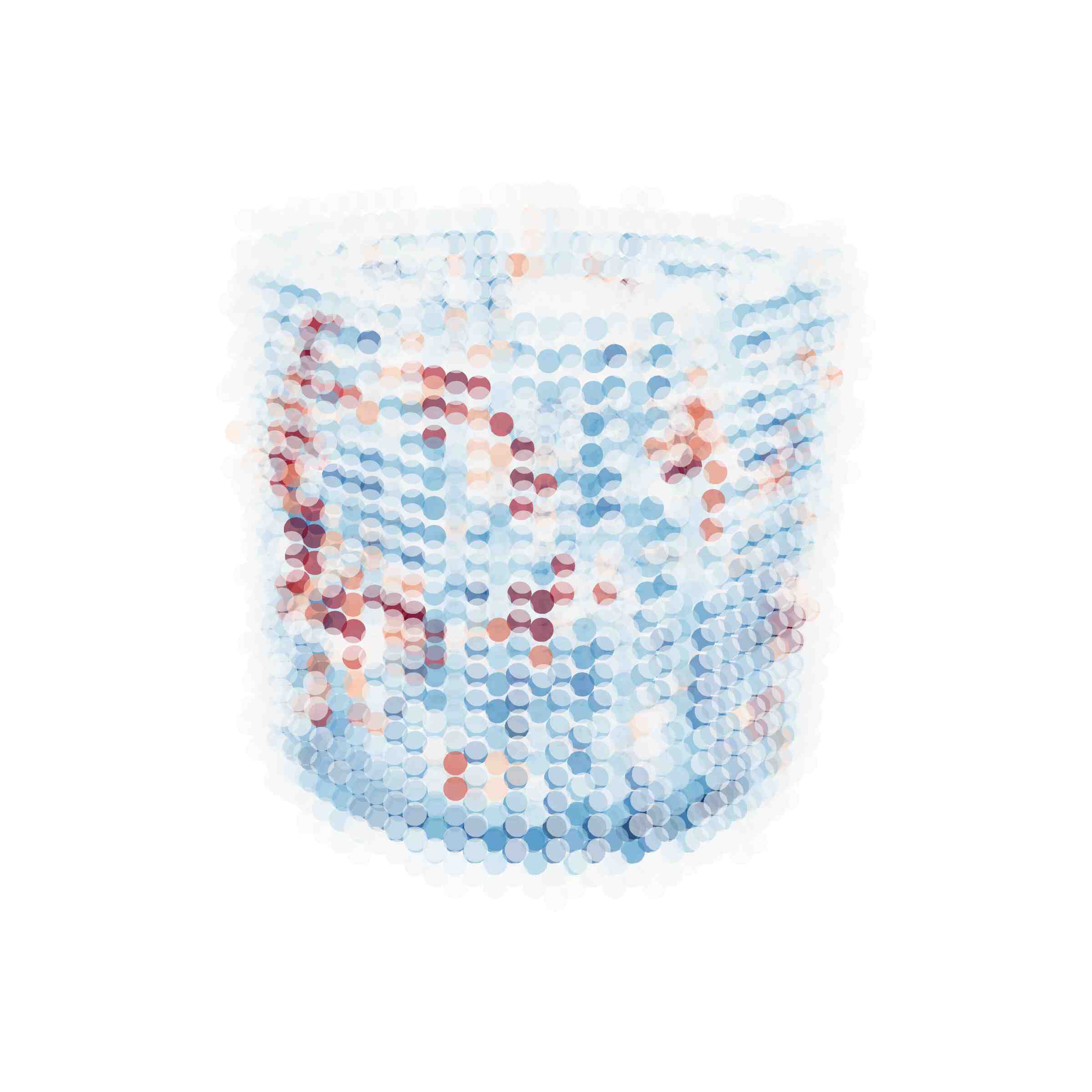}
    \includegraphics[height=17mm]{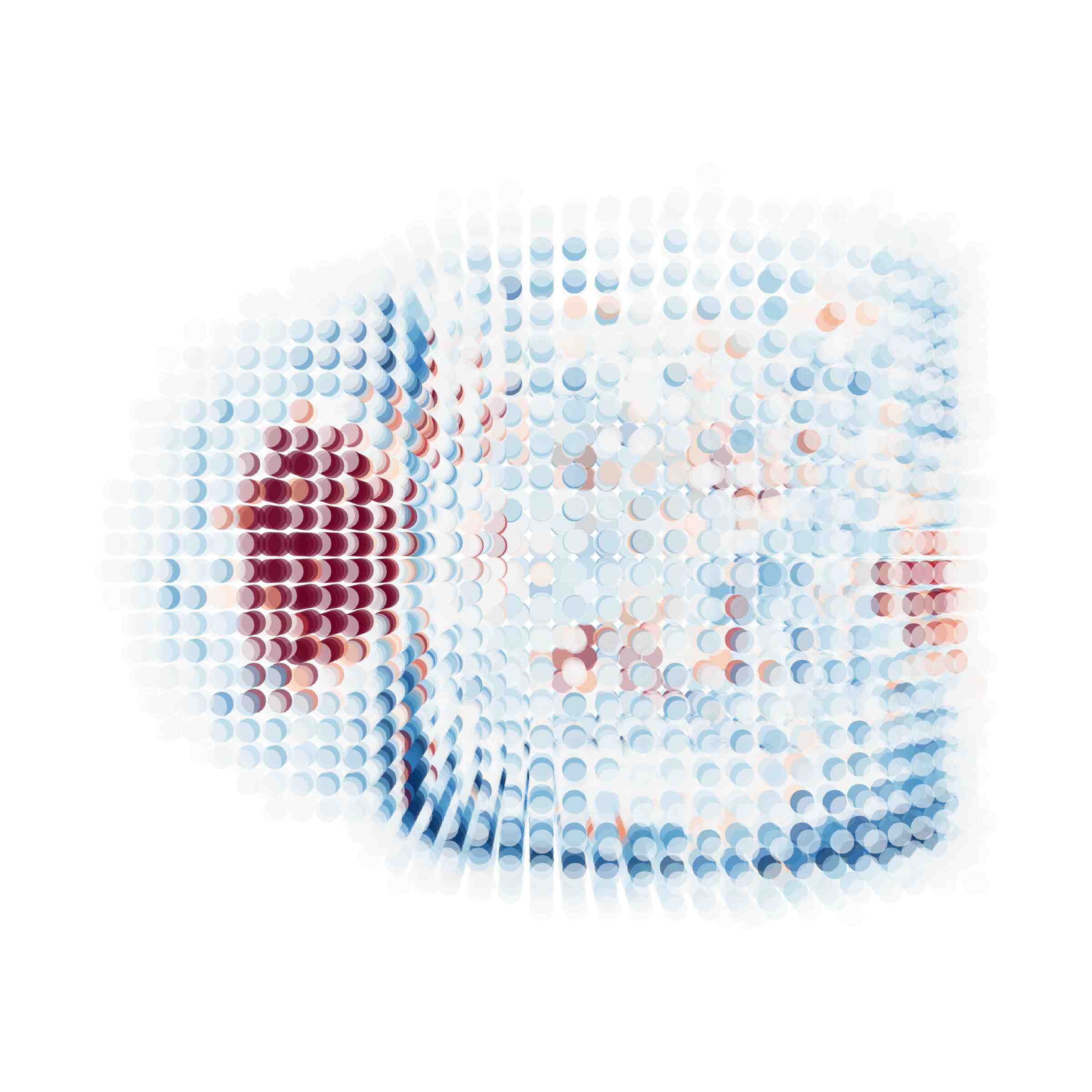}
    \includegraphics[height=17mm]{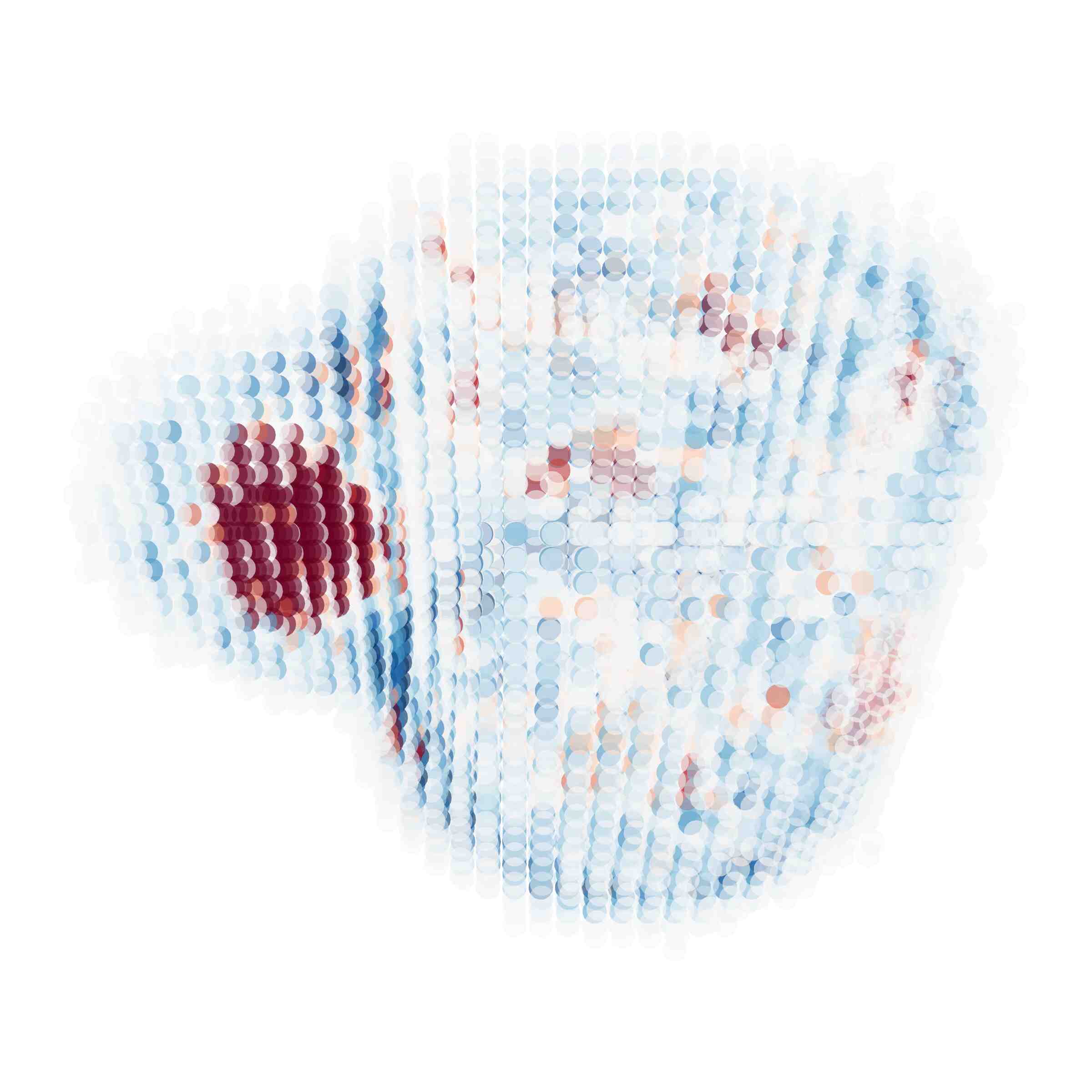}
    \includegraphics[height=17mm]{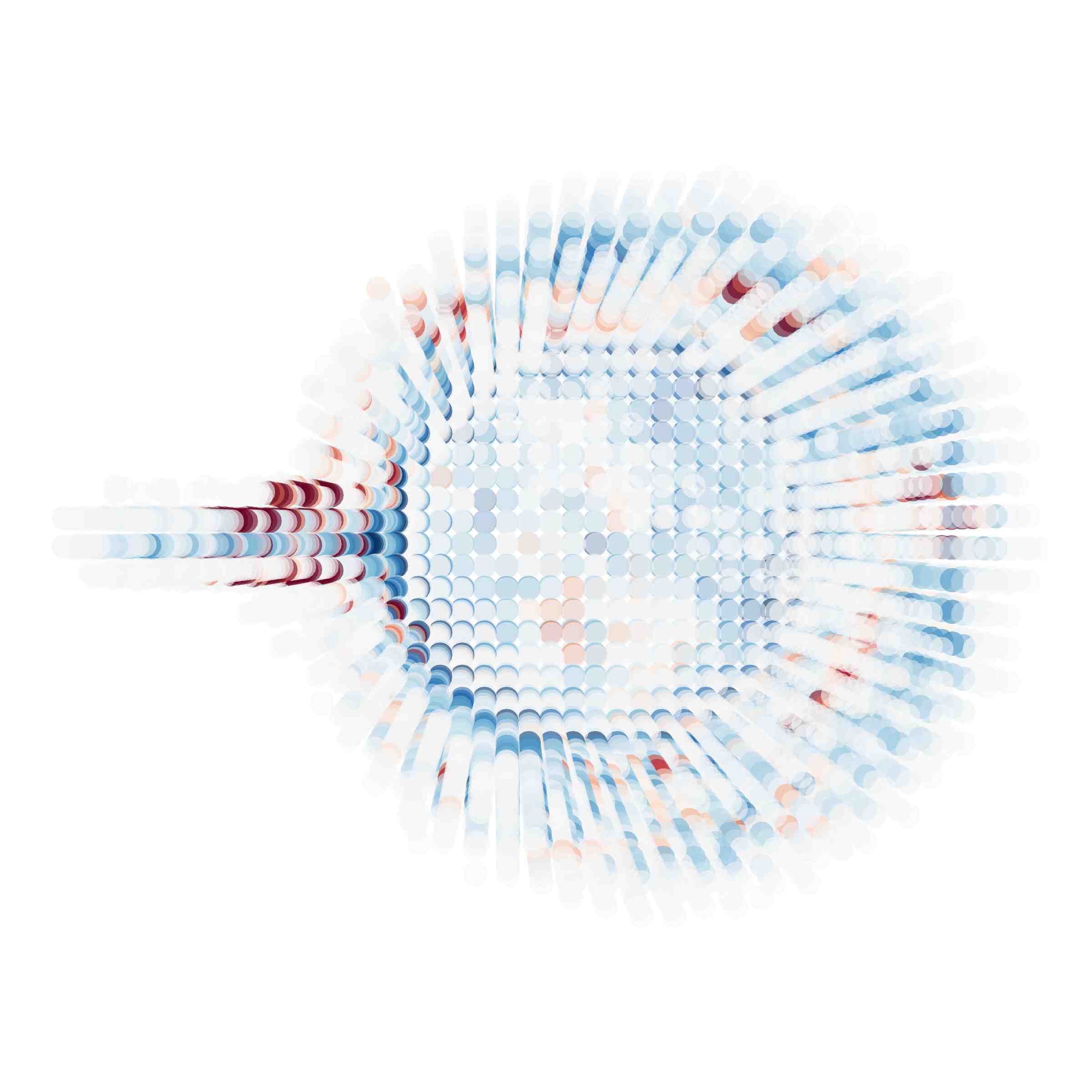}
    \includegraphics[height=17mm]{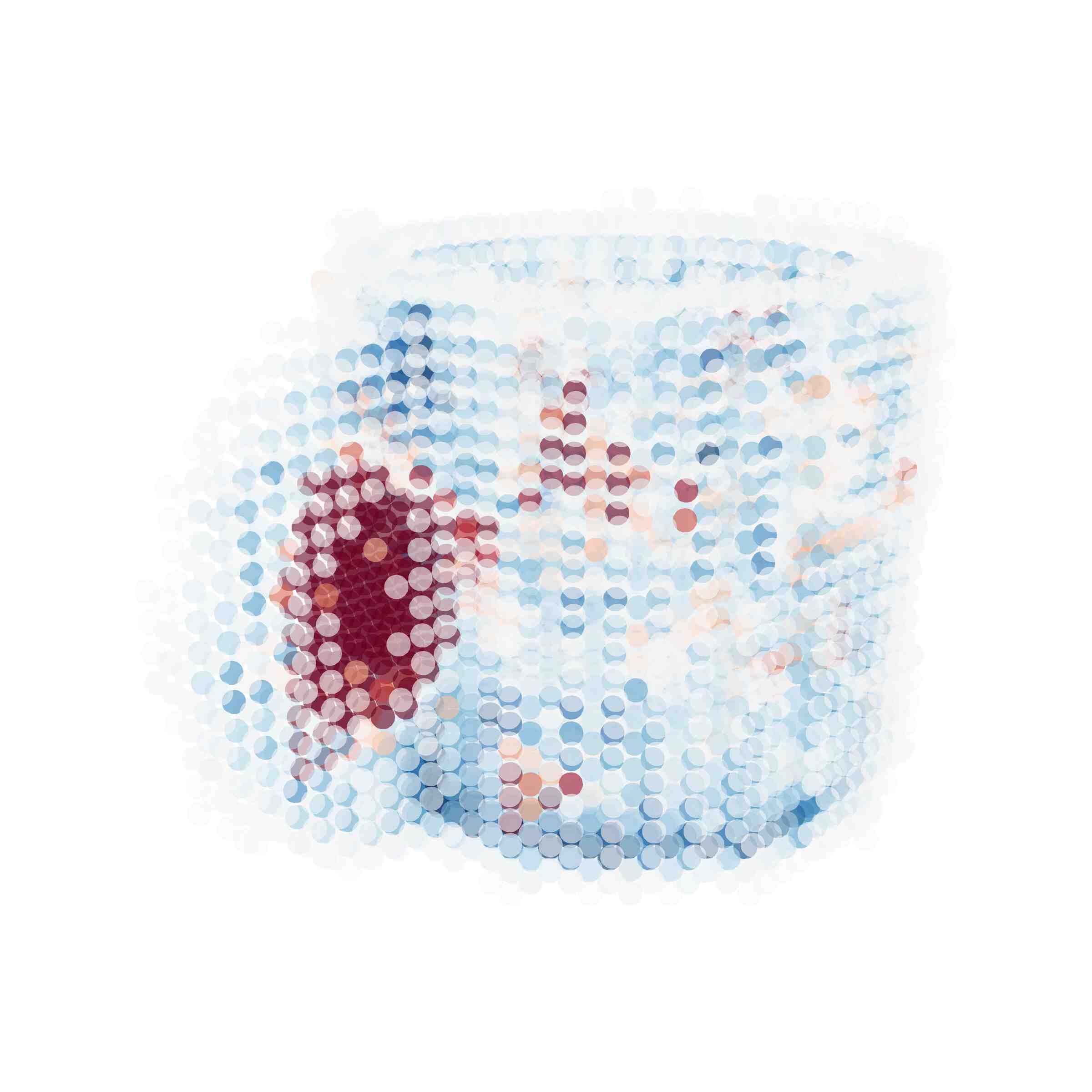}

    \includegraphics[height=17mm]{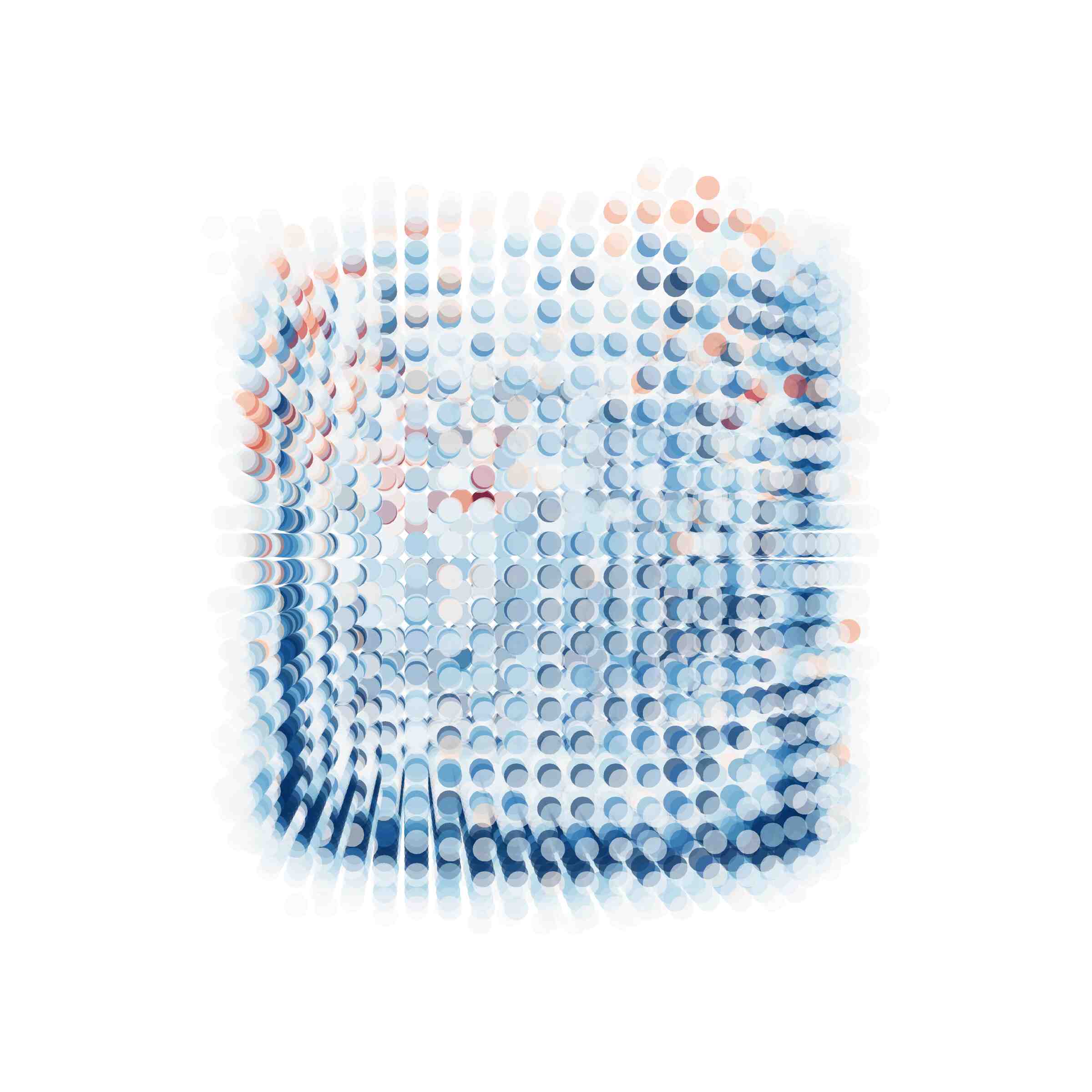}
    \includegraphics[height=17mm]{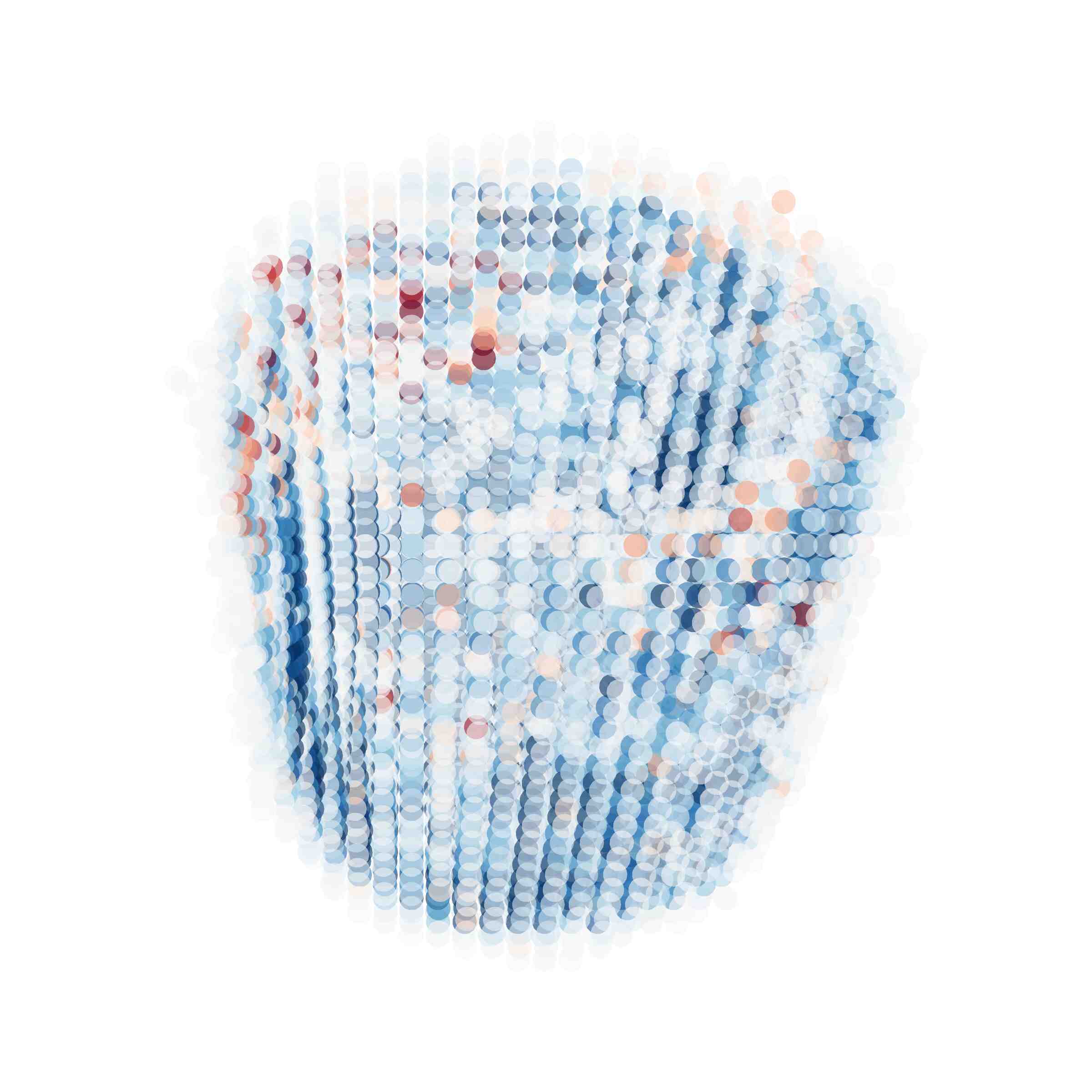}
    \includegraphics[height=17mm]{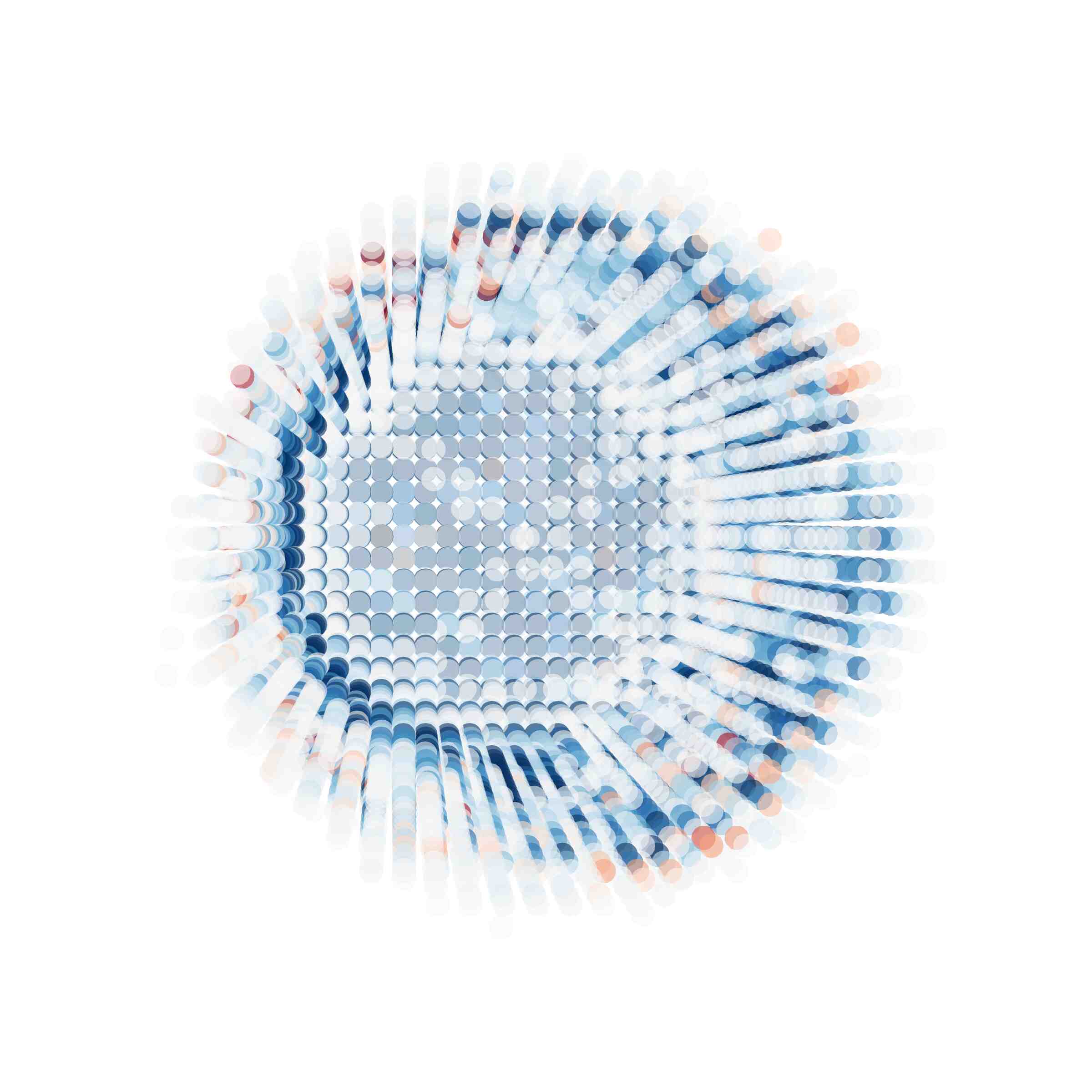}
    \includegraphics[height=17mm]{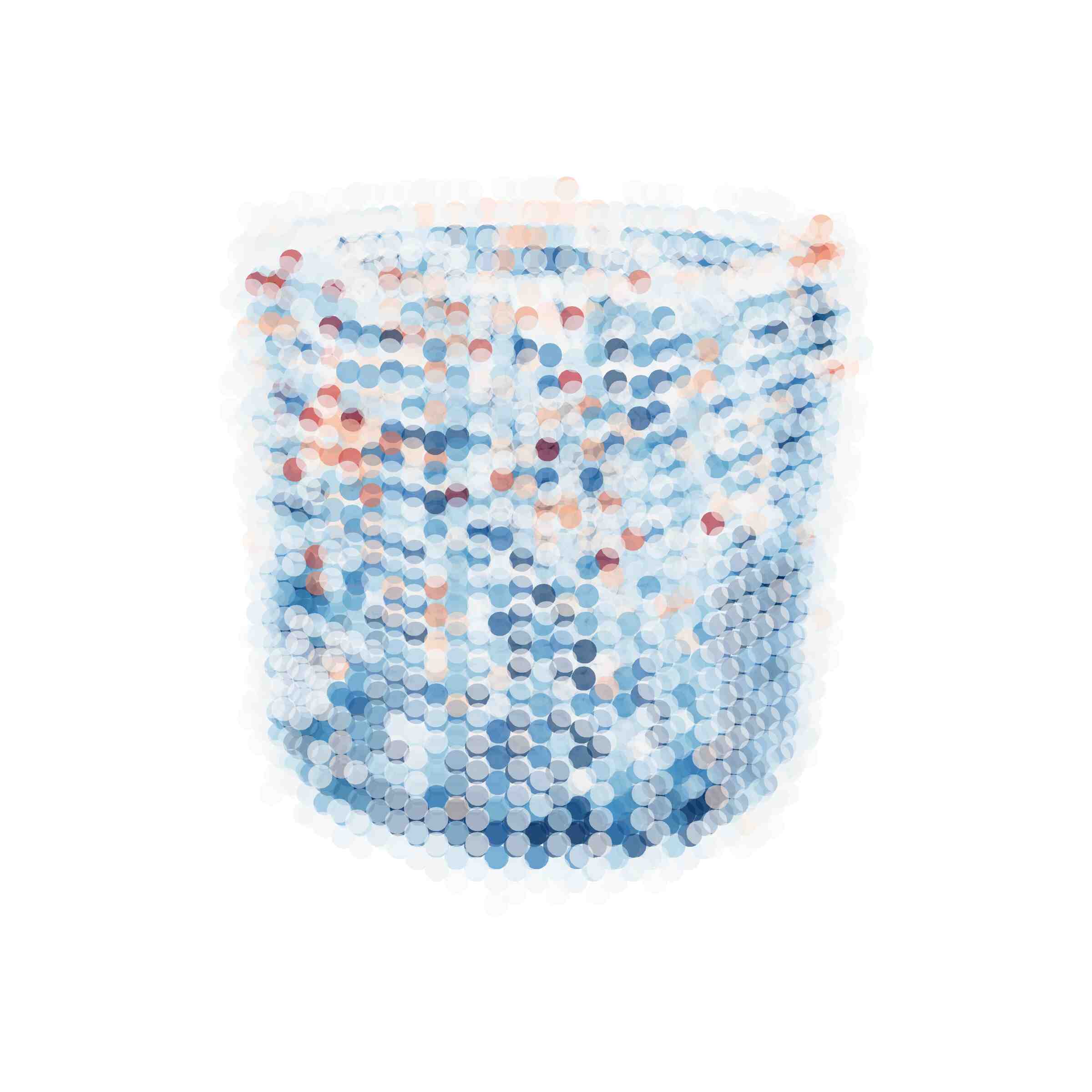}
    \includegraphics[height=17mm]{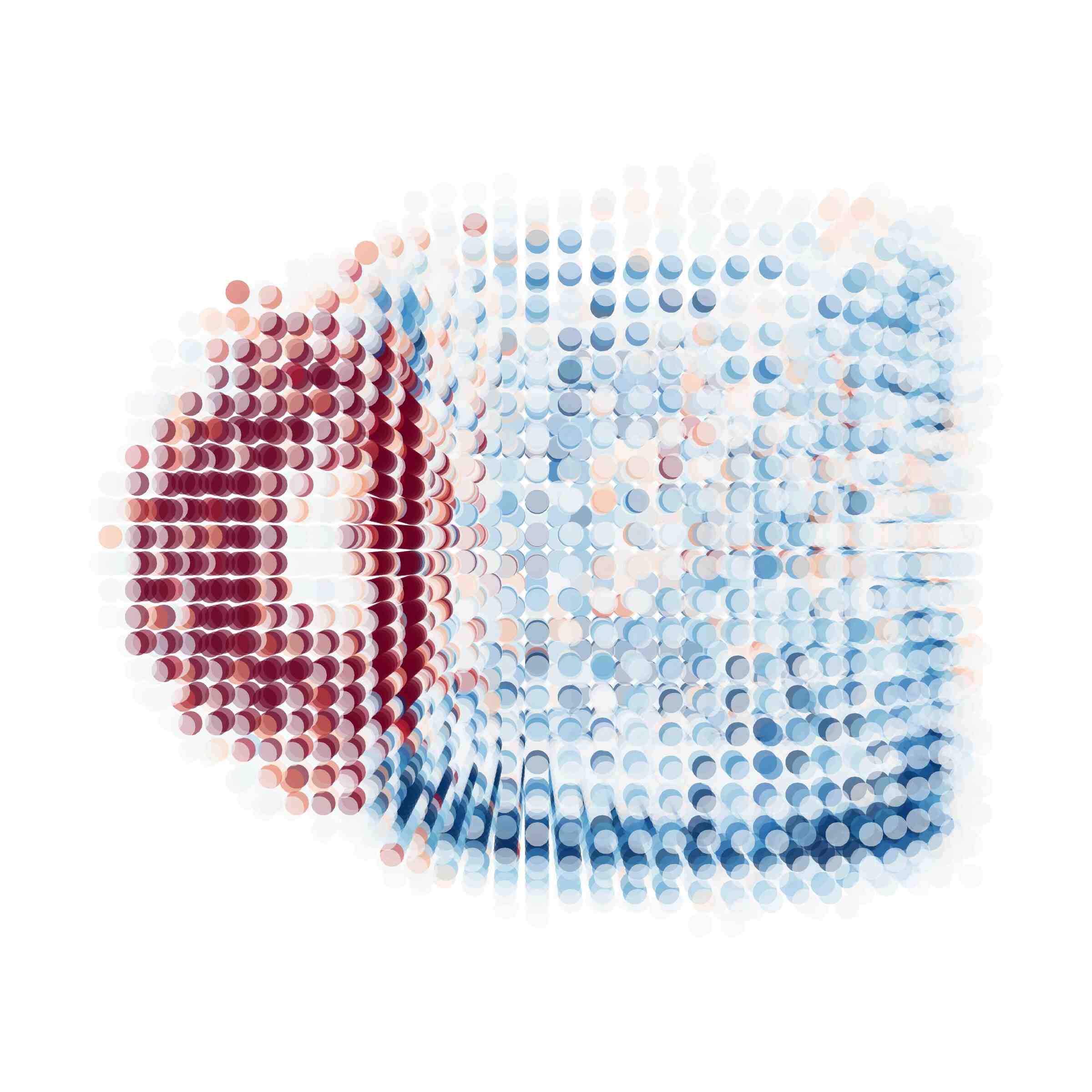}
    \includegraphics[height=17mm]{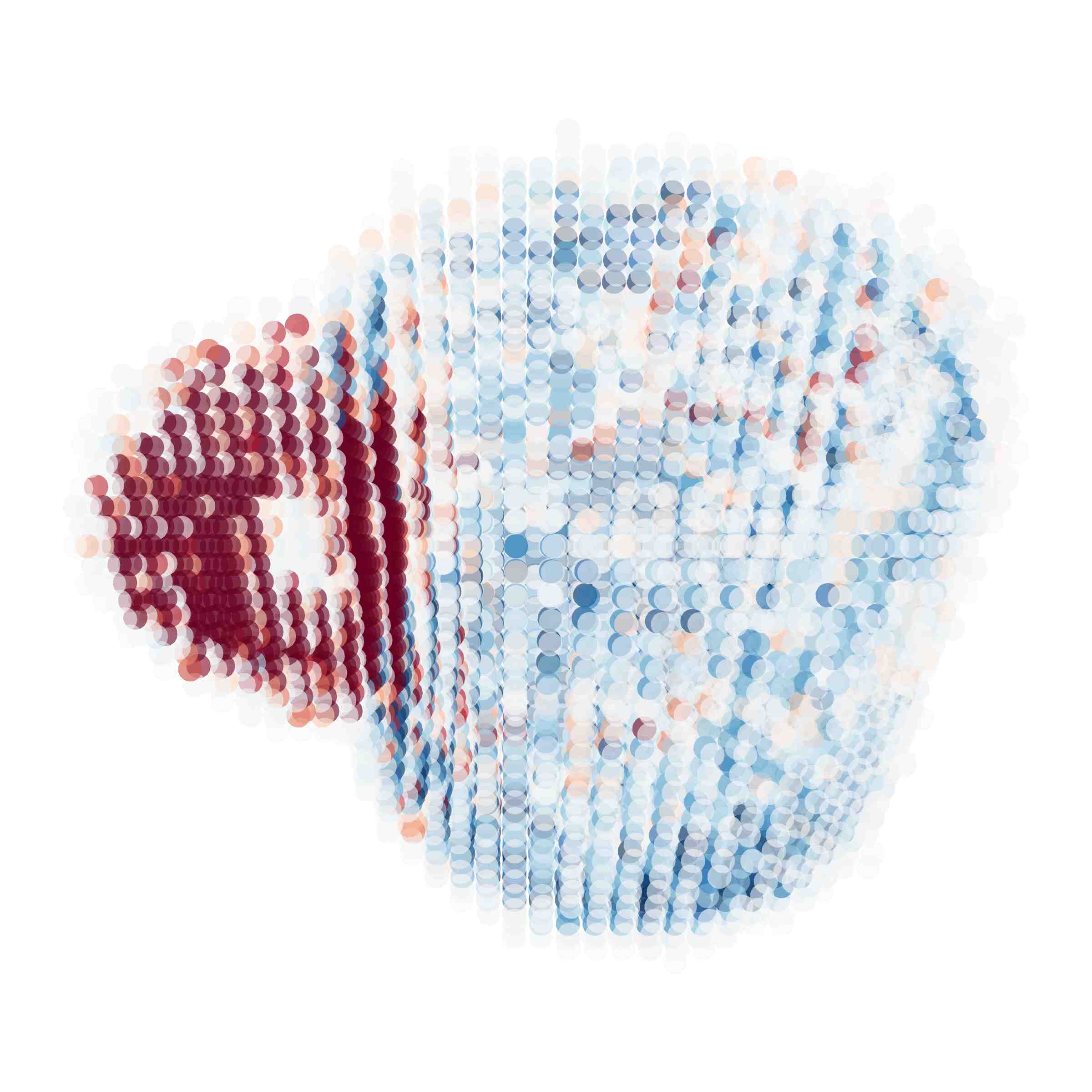}
    \includegraphics[height=17mm]{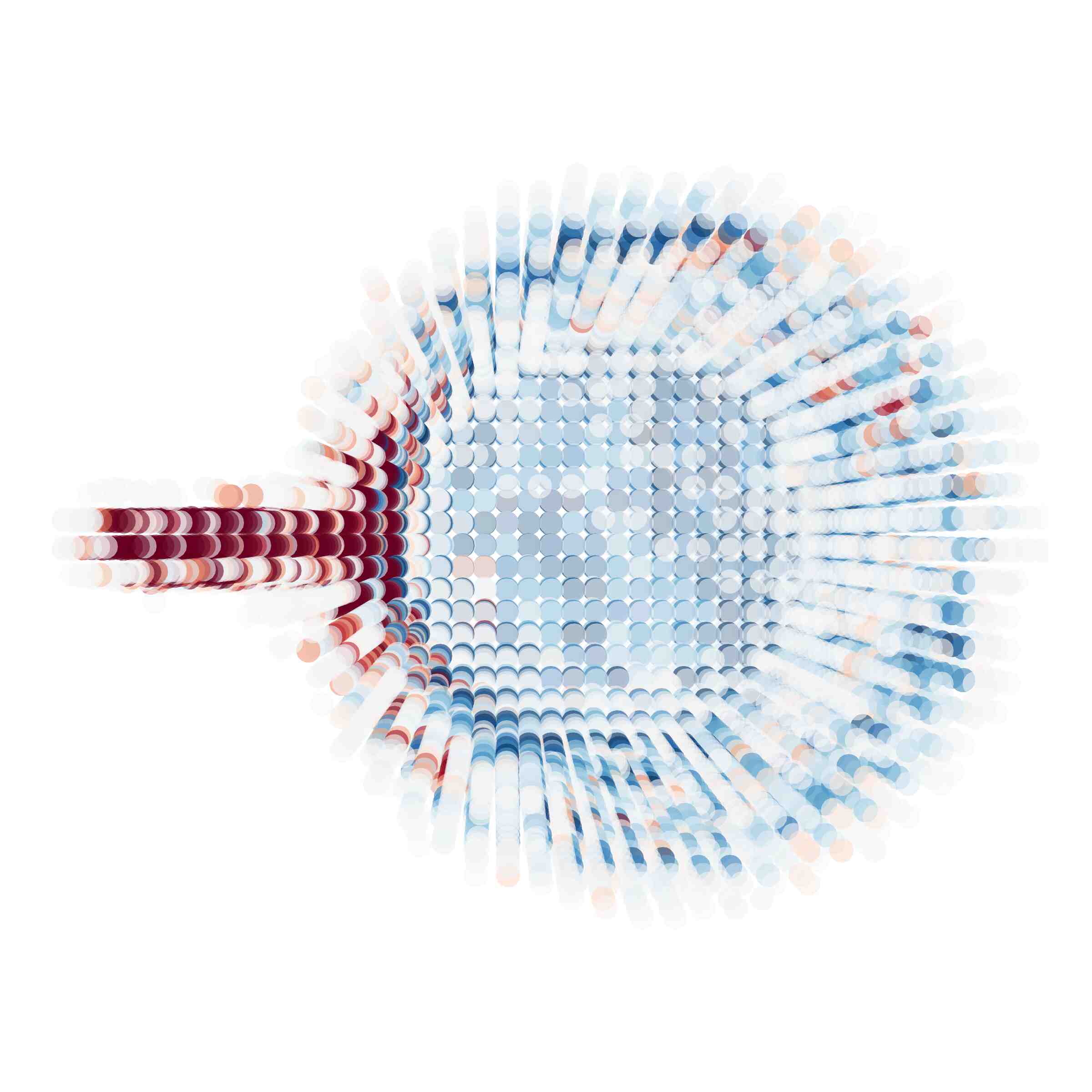}
    \includegraphics[height=17mm]{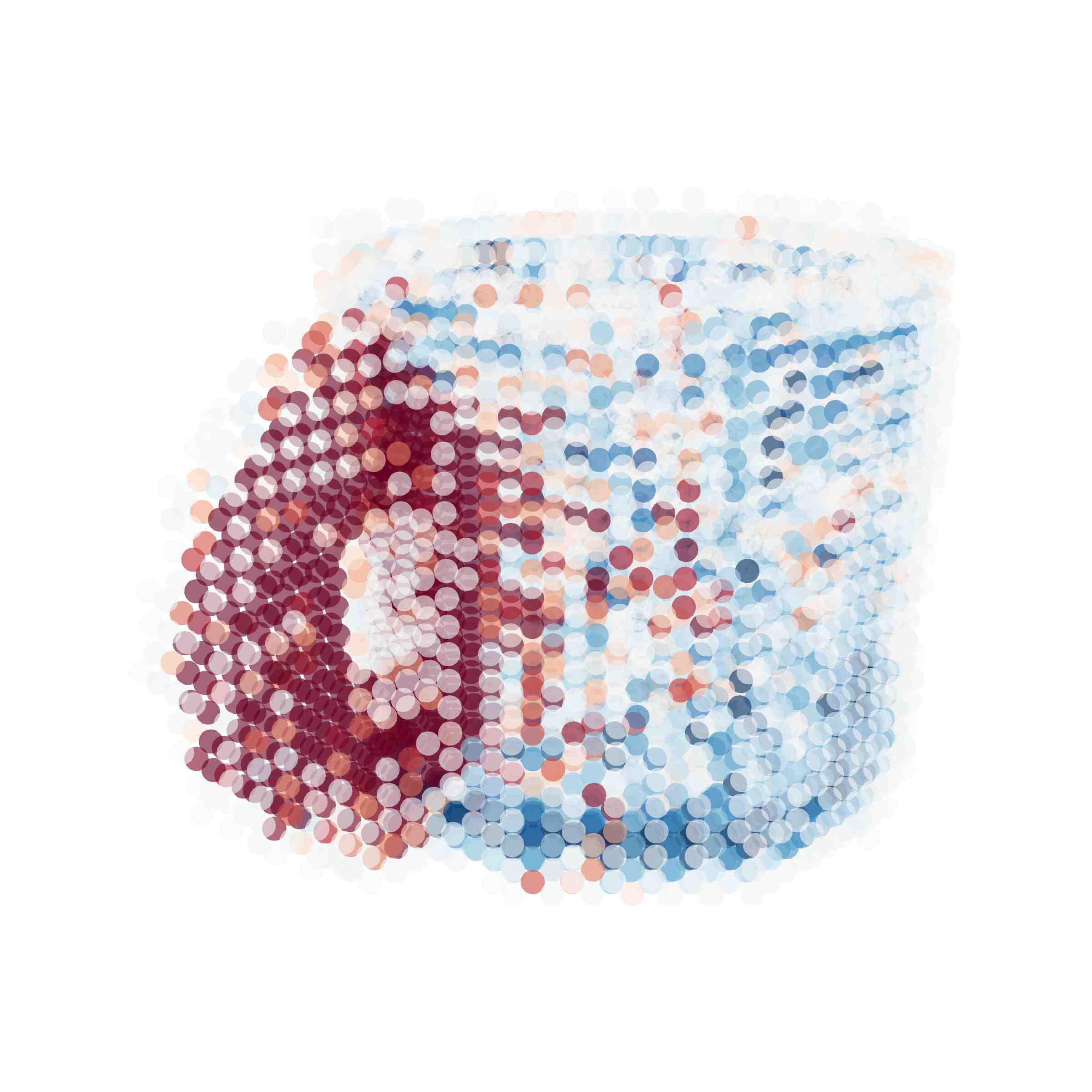}

    \includegraphics[height=17mm]{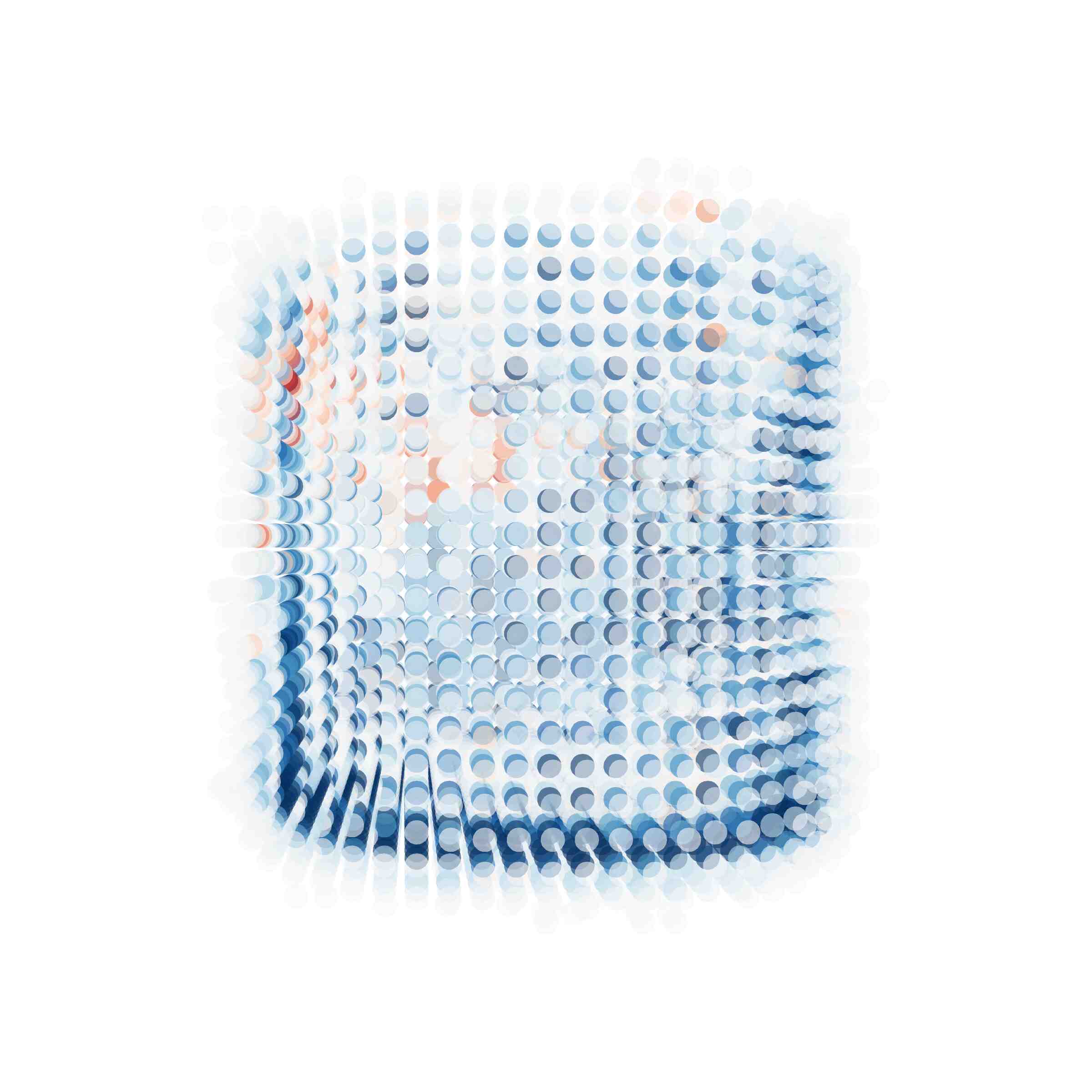}
    \includegraphics[height=17mm]{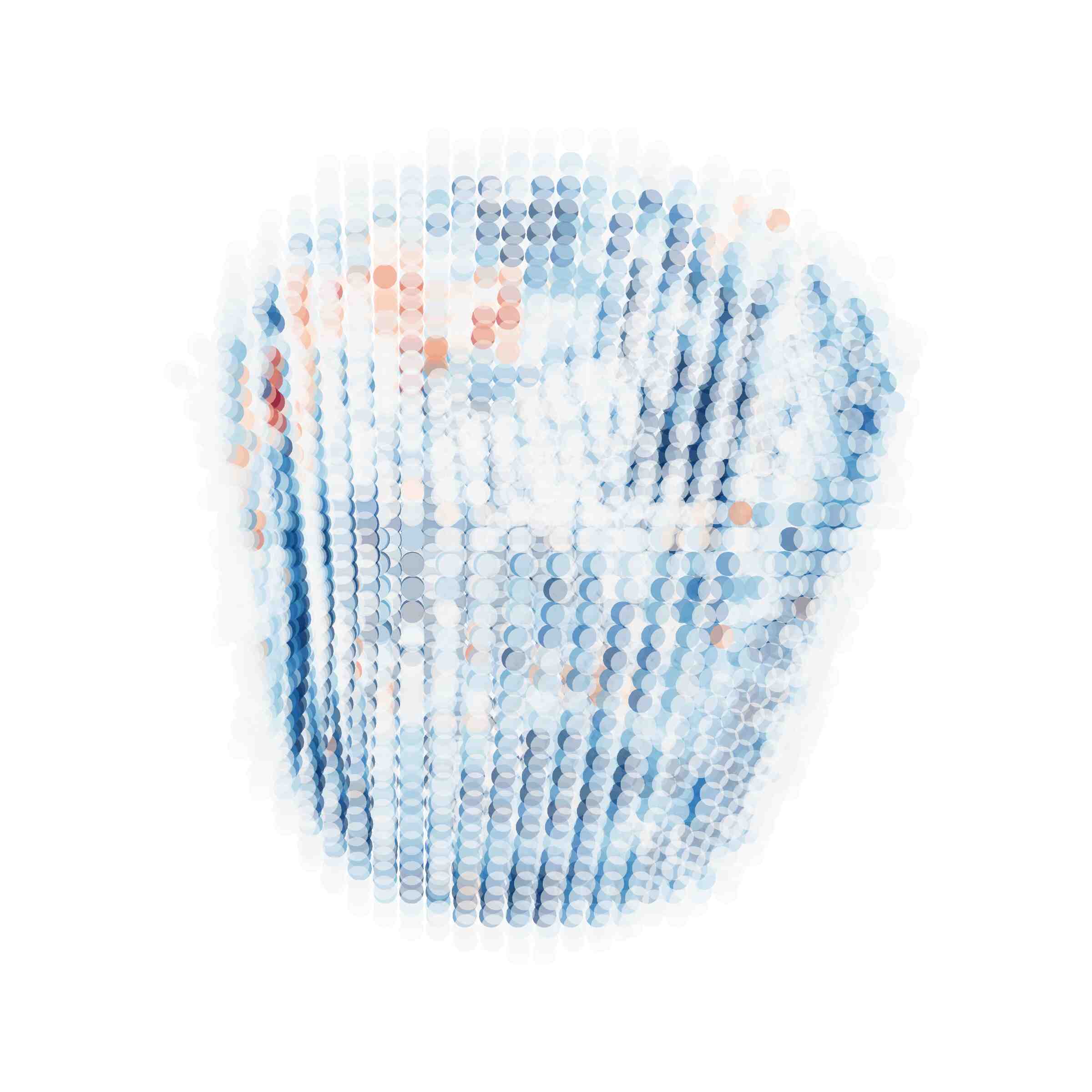}
    \includegraphics[height=17mm]{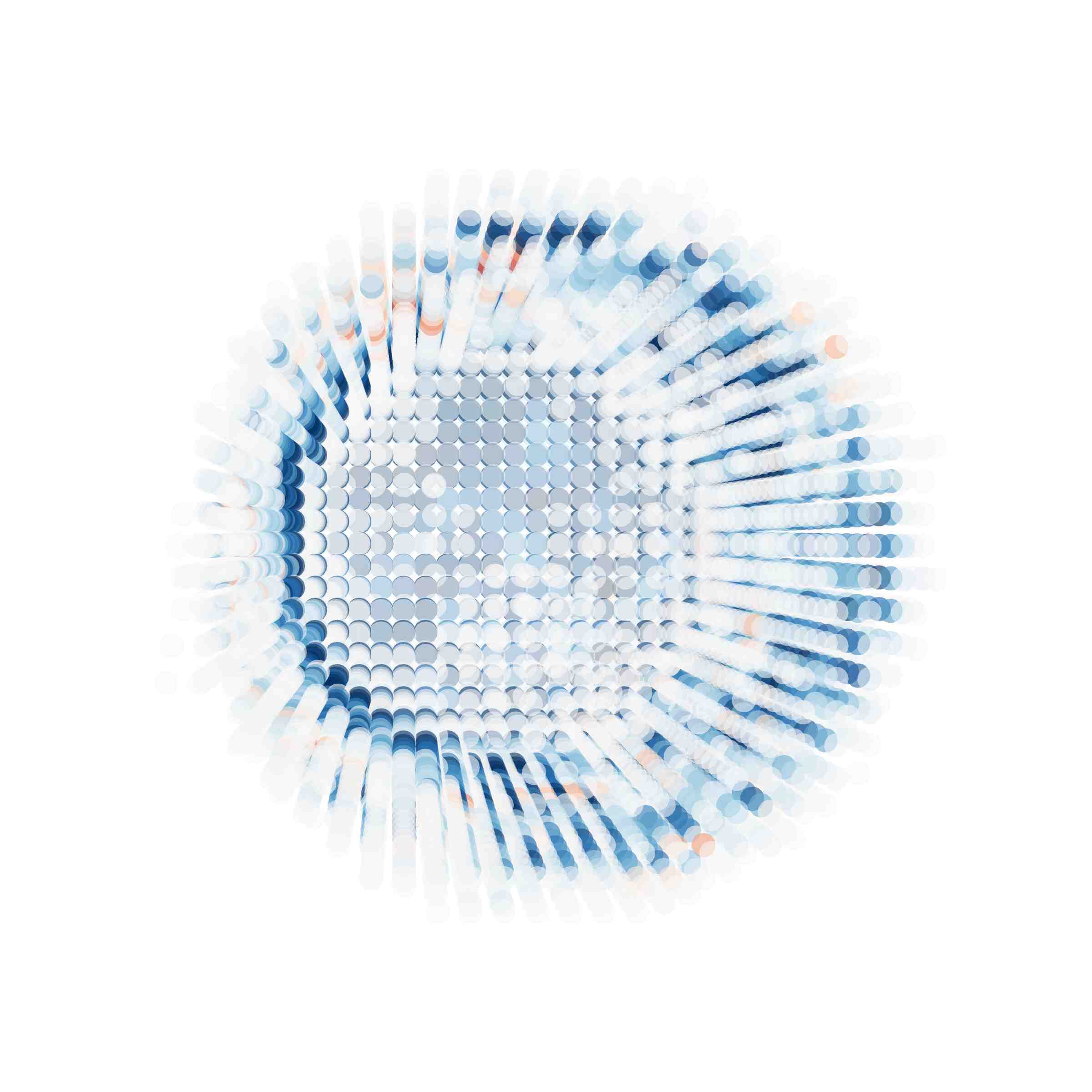}
    \includegraphics[height=17mm]{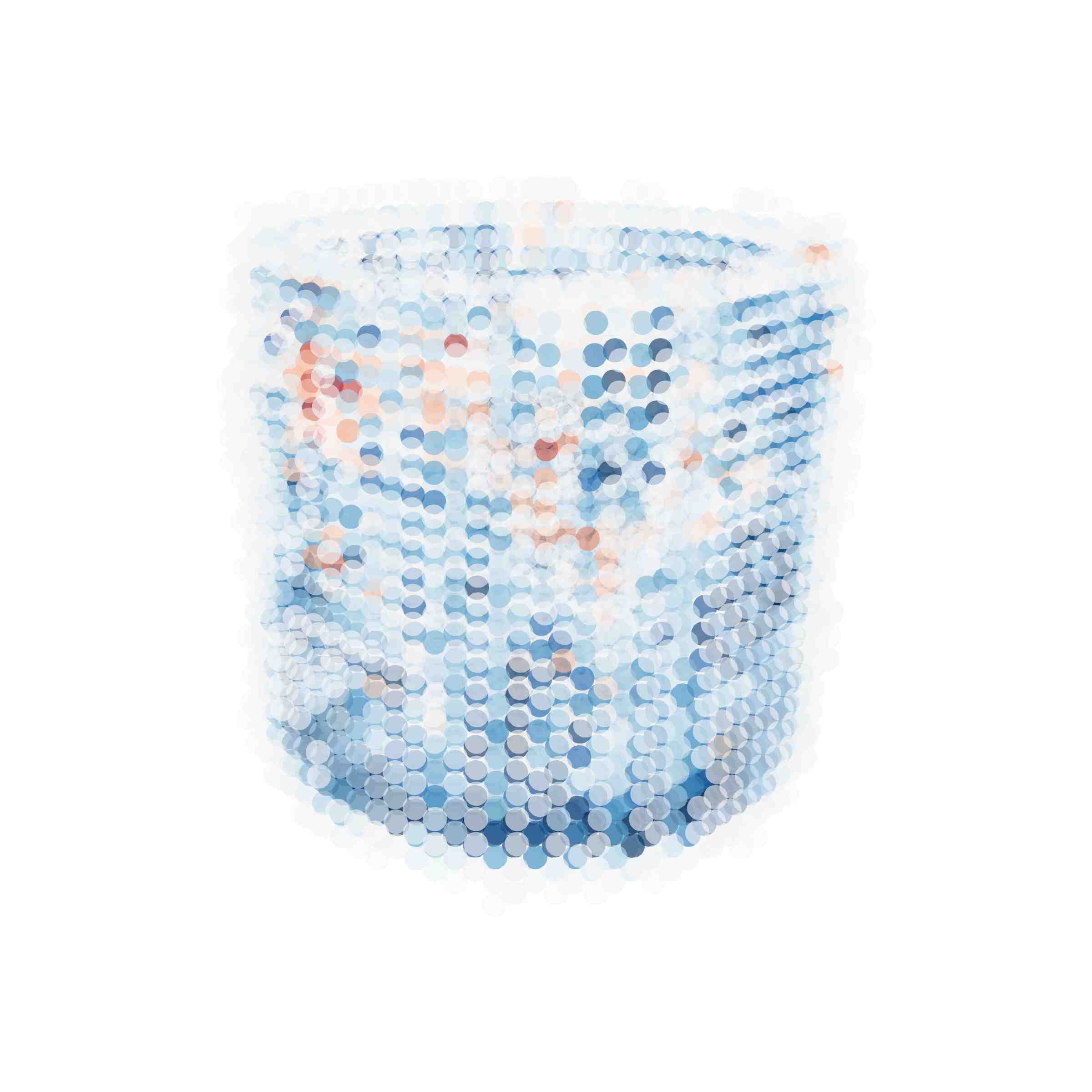}
    \includegraphics[height=17mm]{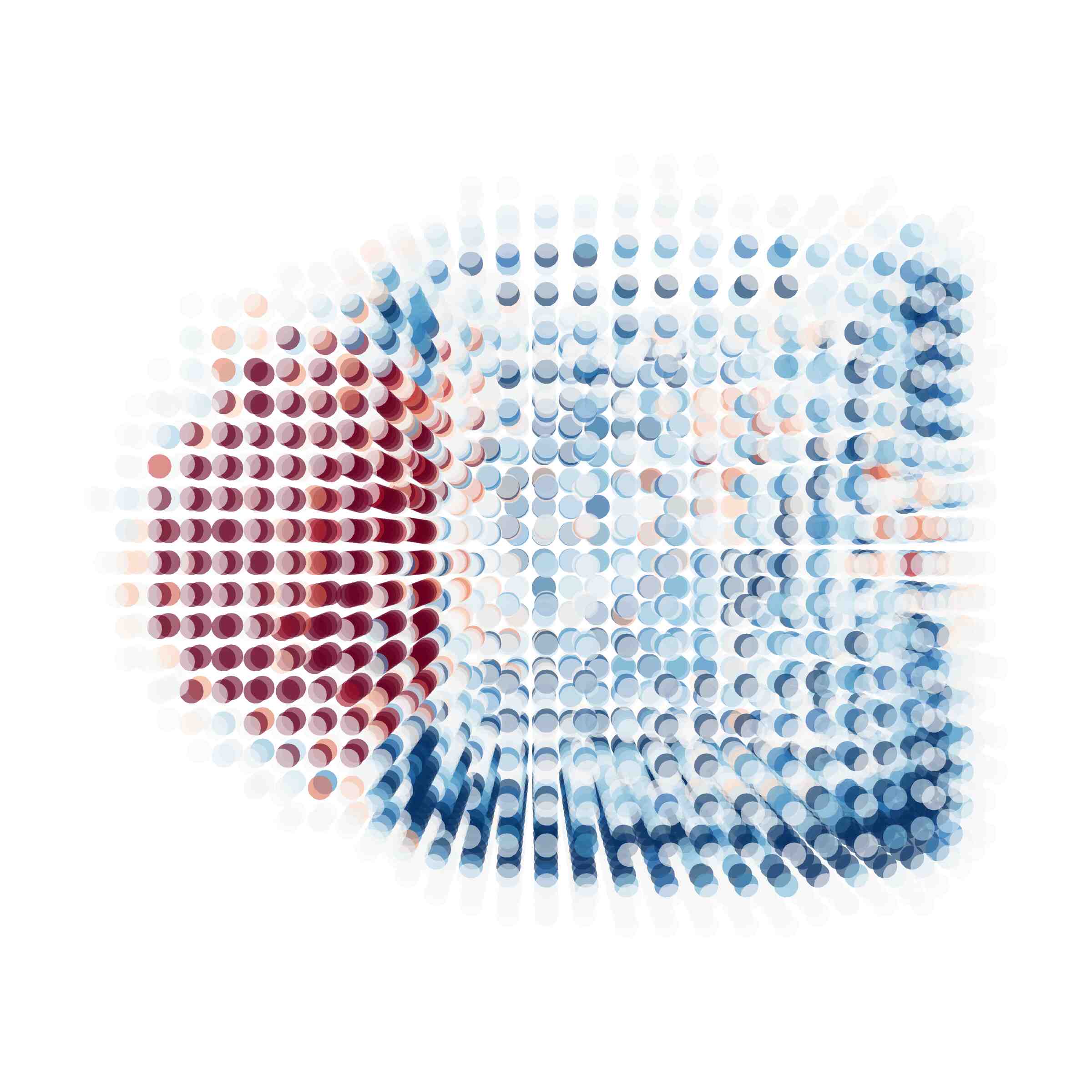}
    \includegraphics[height=17mm]{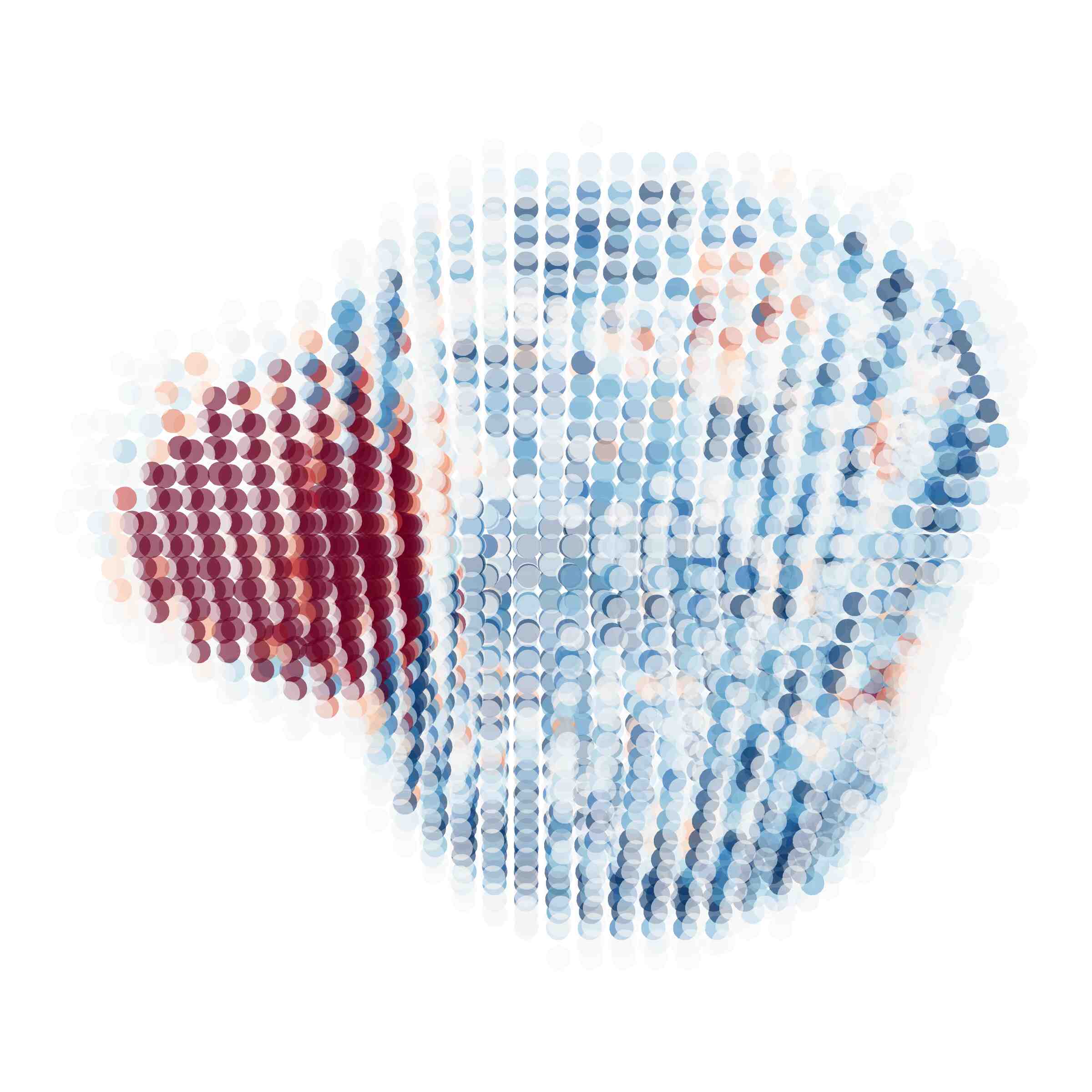}
    \includegraphics[height=17mm]{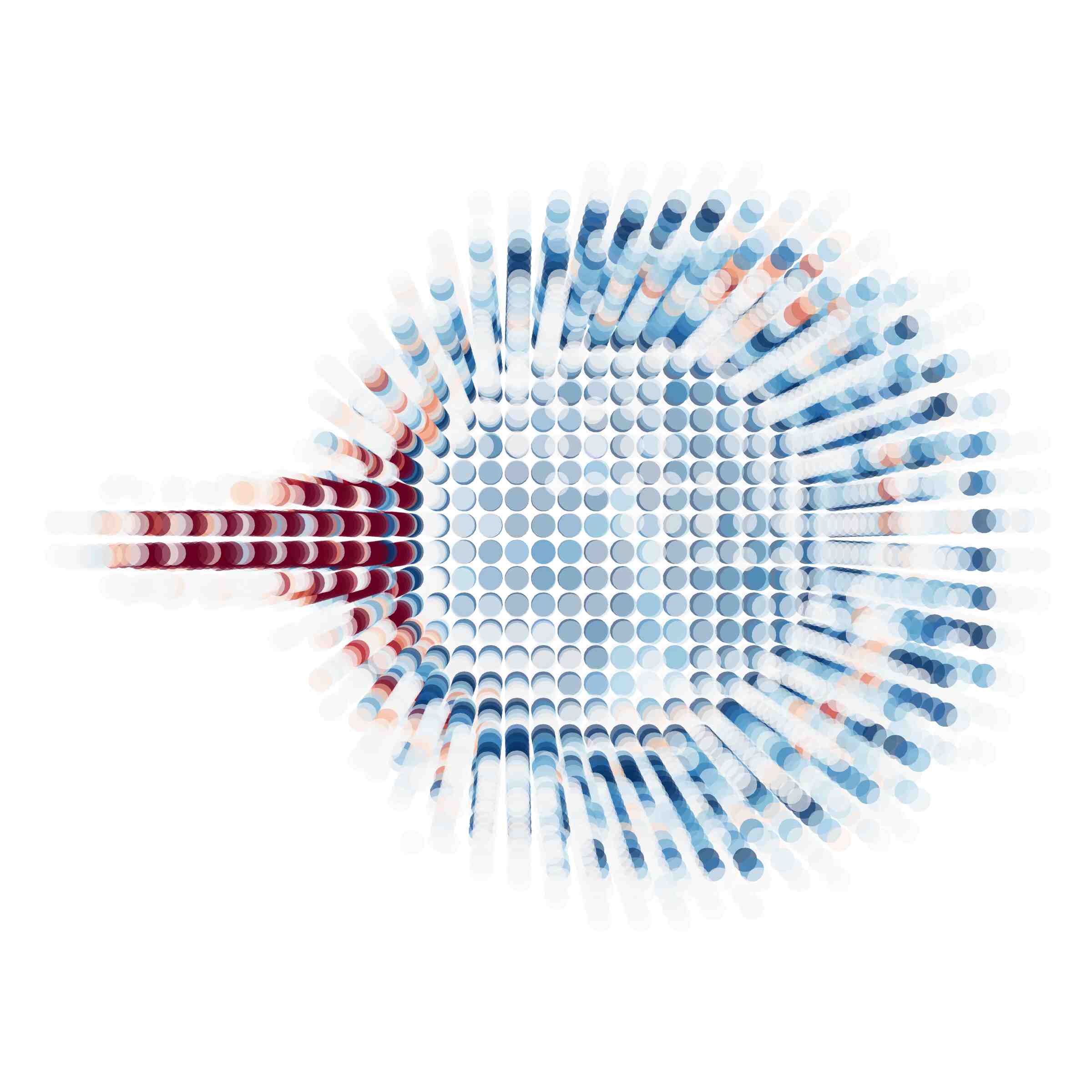}
    \includegraphics[height=17mm]{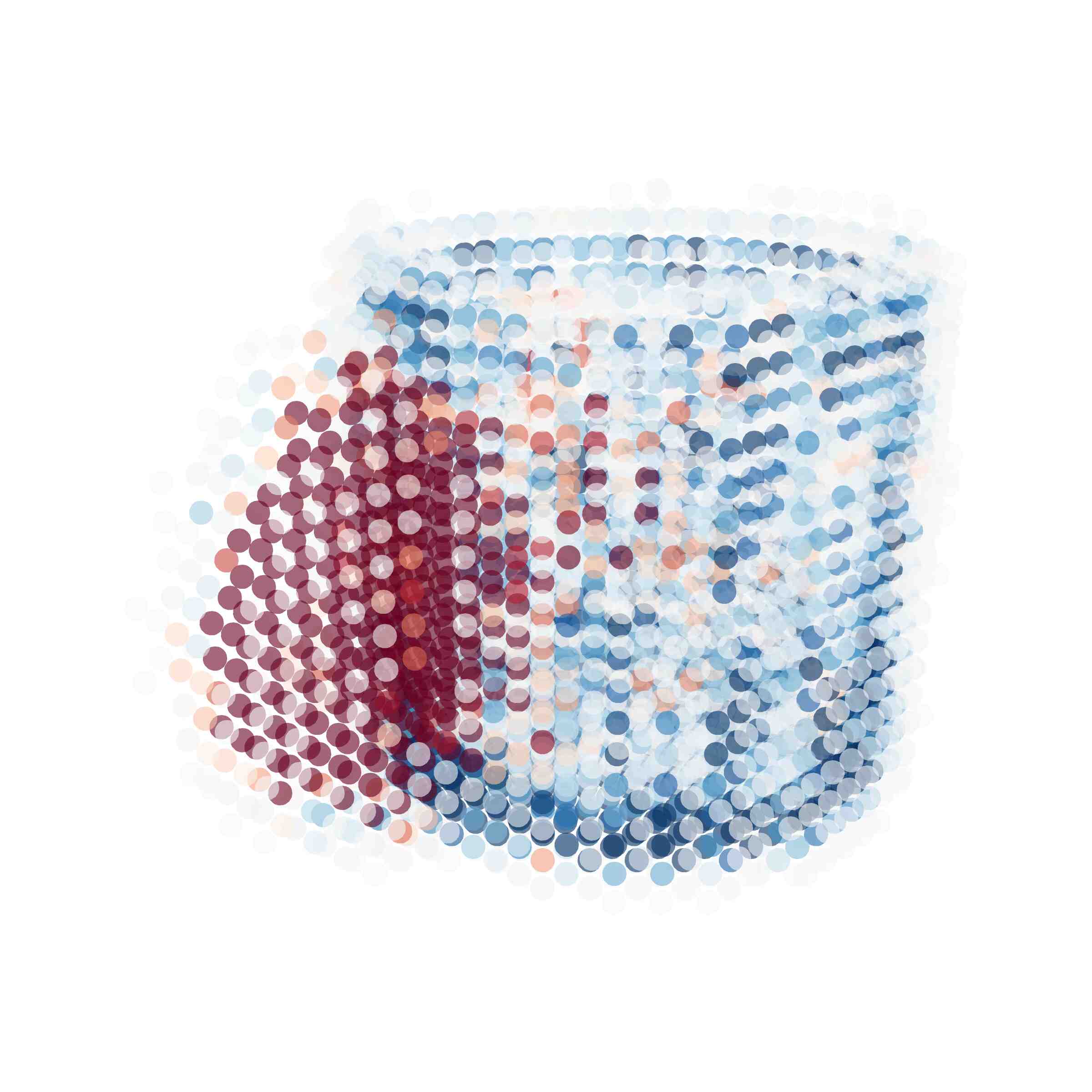}

    \caption{Visualization of the learned classifier that differentiates between
    the cup and the cup with a handle, using homology in degree~$1$ and labeled
    persistence landscapes. 
    Each row shows four viewpoints (side, top-side, top, angled) for the cup
    (left four) and the cup with a handle (right four).
    Rows correspond to birth edges (first row), death triangles (second row),
    representative cycles (third row), and bounding chains (fourth row).
    In each image the sum of the values of the voxels equals the empirical mean value of the classifier.
    }
    \label{fig:cup-handle-visualization-mean}
\end{figure}

We compute the Delaunay complex and the persistence diagram for homology in
degree~$1$ and the corresponding persistence landscape.
We train a linear SVM on training data consisting of $100$ samples from each class.
We then compute the empirical mean feature maps for $1000$ samples from each class.
See \cref{fig:cup-handle-visualization-mean}.

Note that the cup, which is common to both classes, is colored similarly in the two classes, whereas the handle is strongly colored by the color corresponding to the class of the cup with handle.



\subsection{Anomaly detection in a time series}
\label{sec:time-series-diff-phase}

In this example, 
we train a classifier to detect anomalies in a time series and mean feature maps to visualize the location of the anomaly.

We have two classes. The \emph{normal} class consists of time series 
$y_0,y_1,\ldots,y_{250}$, where 
\[
y_i = \sin\left(\frac{2 \pi i}{50}\right) + \varepsilon_i,
\]
with
$\varepsilon_i \sim N(0,\sigma^2)$, where $\sigma = 0.1$.
The \emph{anomalous} class is given by 
$z_0,z_1,\ldots,z_{250}$, where $N \in [0,\ldots,200]$, and
\[
z_i = y_i + 2 \sin \left( \frac{4 \pi (i - N)}{50} \right),
\]
if $i \in [N,N+50]$ and $z_i = y_i$ otherwise.
%

Each time series is transformed into a point cloud via a delay embedding with delay
$1$
and embedding dimension $d = 10$, yielding a point cloud of $242$ points in
$\mathbb{R}^{10}$.
We let $(K, w)$ be the corresponding Vietoris-Rips complex,
which is the complete simplicial complex on $242$ points in which the weight of a simplex is given by the maximum pairwise distance of its vertices.

We compute the 
persistence diagrams for homology in degree $1$
and the corresponding persistence landscapes.
We train a linear SVM classifier on $50$ samples per class and obtain a separating
hyperplane.

\begin{figure}[!htb]
    \centering
    \includegraphics[height=25mm]{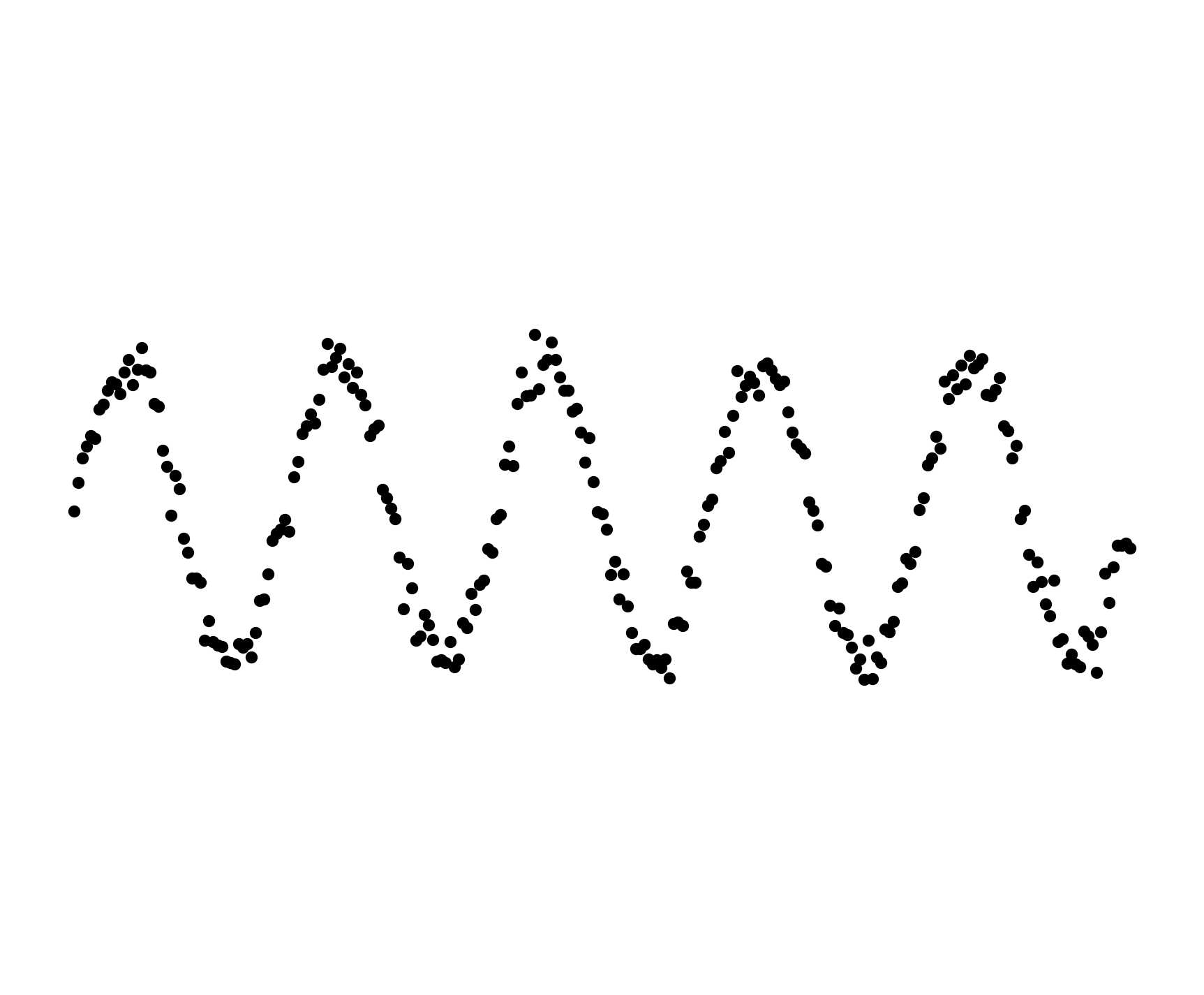}
    \includegraphics[height=25mm]{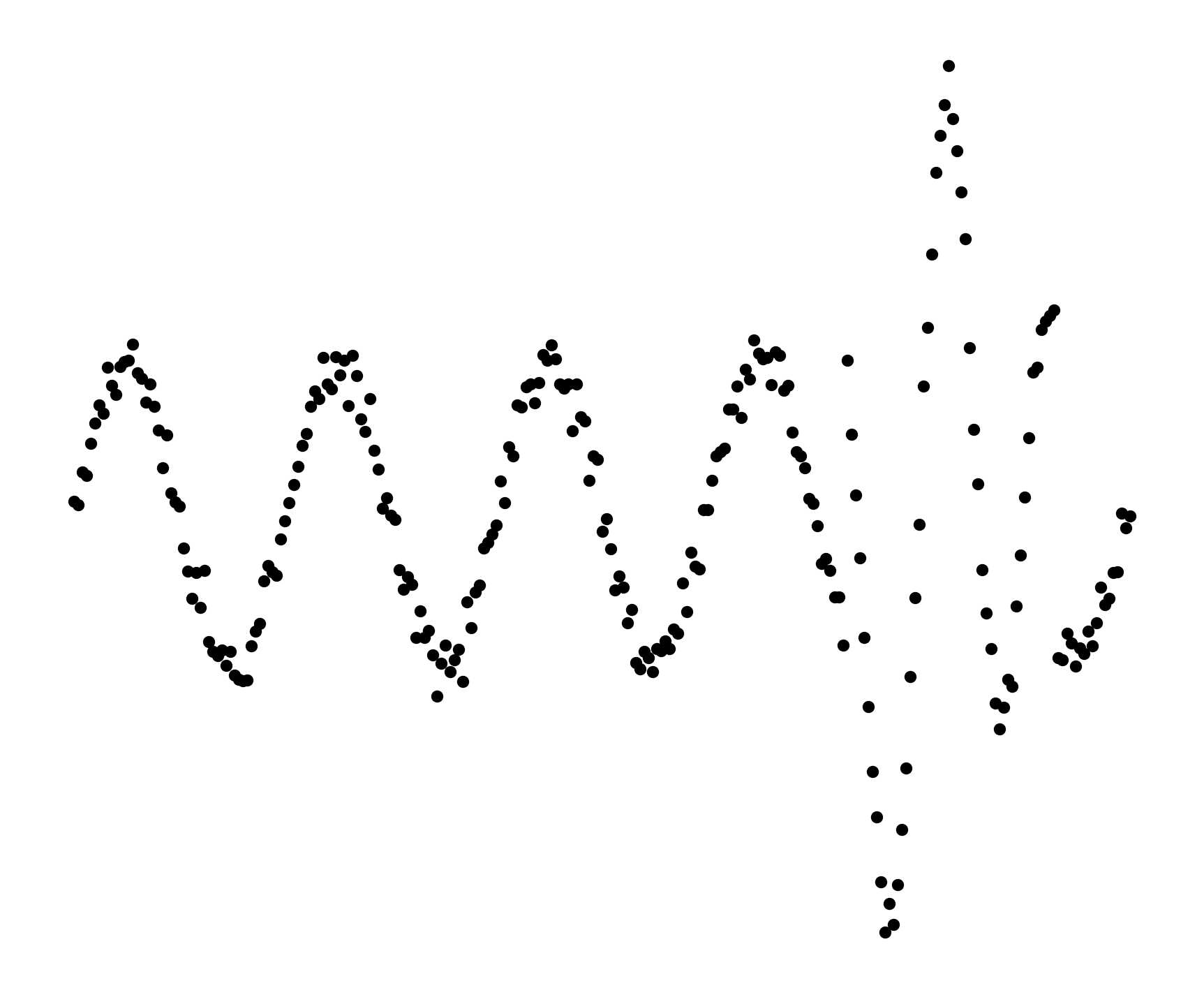} \quad 
    \includegraphics[height=25mm]{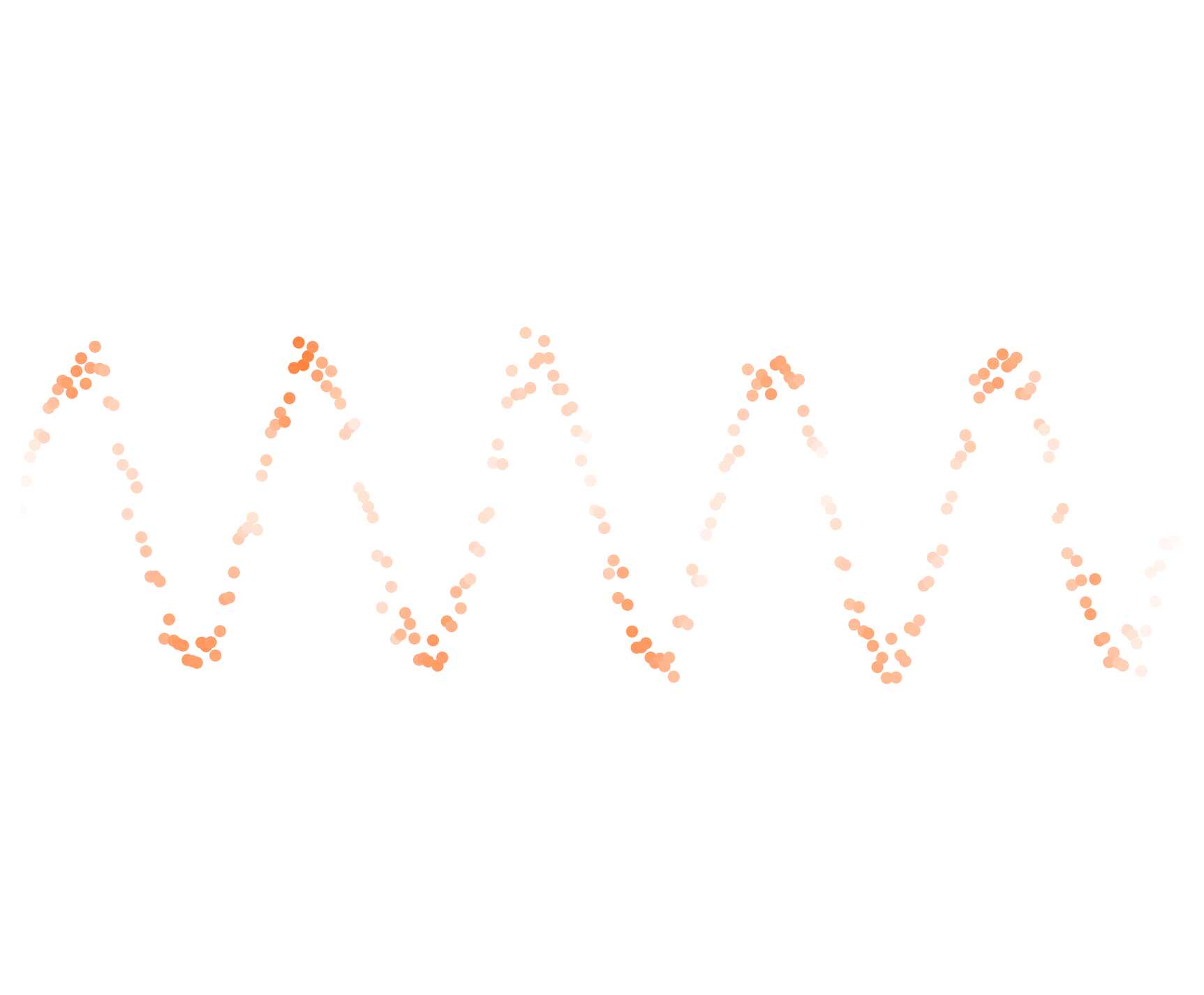}
    \includegraphics[height=25mm]{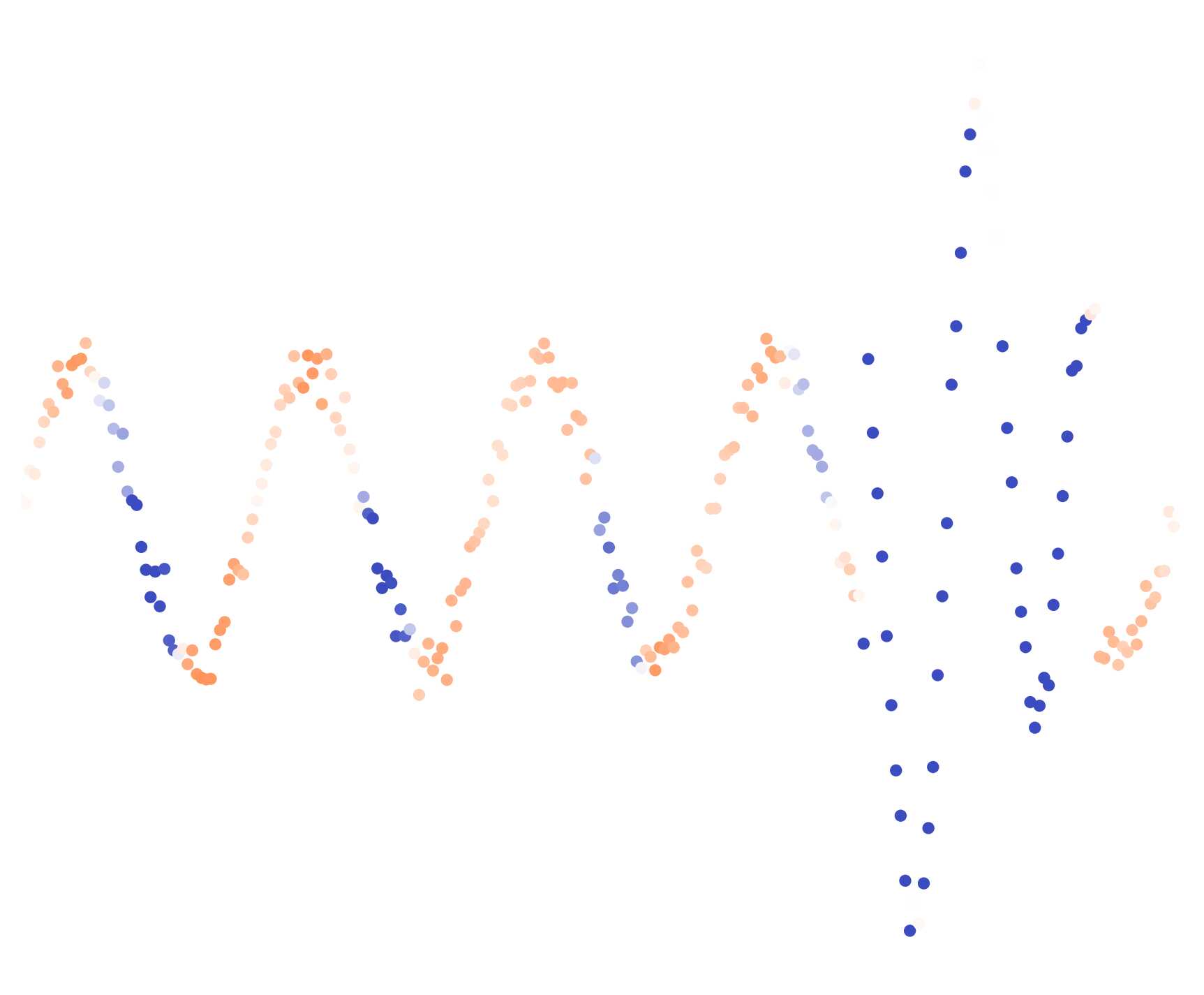}

    \caption{Visualization of the learned classifier that differentiates between normal and
    anomalous time series with a randomly located anomaly, using homology in degree $1$ and
    birth-simplex-labeled persistence landscapes.
    Left: a normal time series and an anomalous time series. 
    Right: Visualizations for the normal and anomalous time series.
    We use a blue-to-orange color gradient in which blue corresponds to the anomalous class and orange corresponds to the normal class.
    In each image the sum of the values of the points equals the empirical mean value of the classifier.
    }
    \label{fig:anomaly-time-series-diff-phase-deg1}
\end{figure}

We then take one sample from each class and compute the empirical mean feature
map for birth edges over $50$ perturbations.
Each perturbation is obtained
by adding independent Gaussian noise 
with mean $0$ and standard deviation $0.05$ to each element of the time series.
The values on the edges in the Vietoris-Rips complex are distributed evenly to the two boundary vertices.
These values and then further distributed evenly to the $10$ corresponding indices of the time series.
%
See \cref{fig:anomaly-time-series-diff-phase-deg1}.


Note that we are using TDA to both detect and locate the anomaly.




\subsection{Classification of discrete dynamical systems}
\label{sec:linked-twist-map}

In this example, we classify five cases of a one-parameter family of discrete dynamical systems.

For a choice of parameter $r$, we have the following discrete dynamical system on a flat torus given by a linked twist map in which each of the twists is given by the logistic map~\cite{Hertzsch:2007}.
\begin{equation} \label{eq:linked-twist-map}
    \begin{aligned}
        x_{n+1} &= x_n + r y_n(1-y_n) \ \mod 1\\
        y_{n+1} &= y_n + r x_{n+1}(1-x_{n+1}) \mod 1
    \end{aligned}
\end{equation}

\begin{figure}[!htb]
    \centering
    \includegraphics[height=25mm]{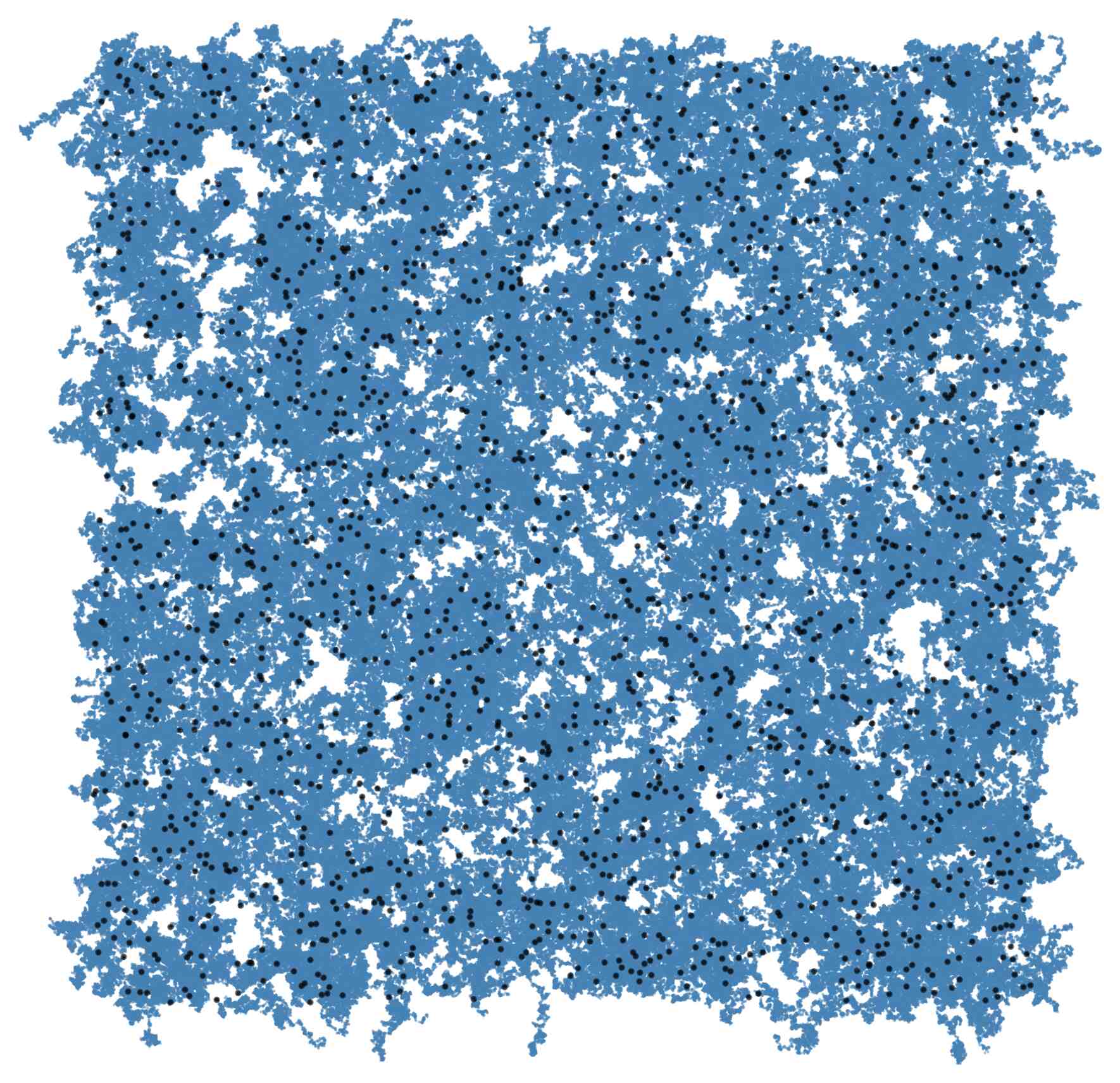}
    \centering
    \includegraphics[height=25mm]{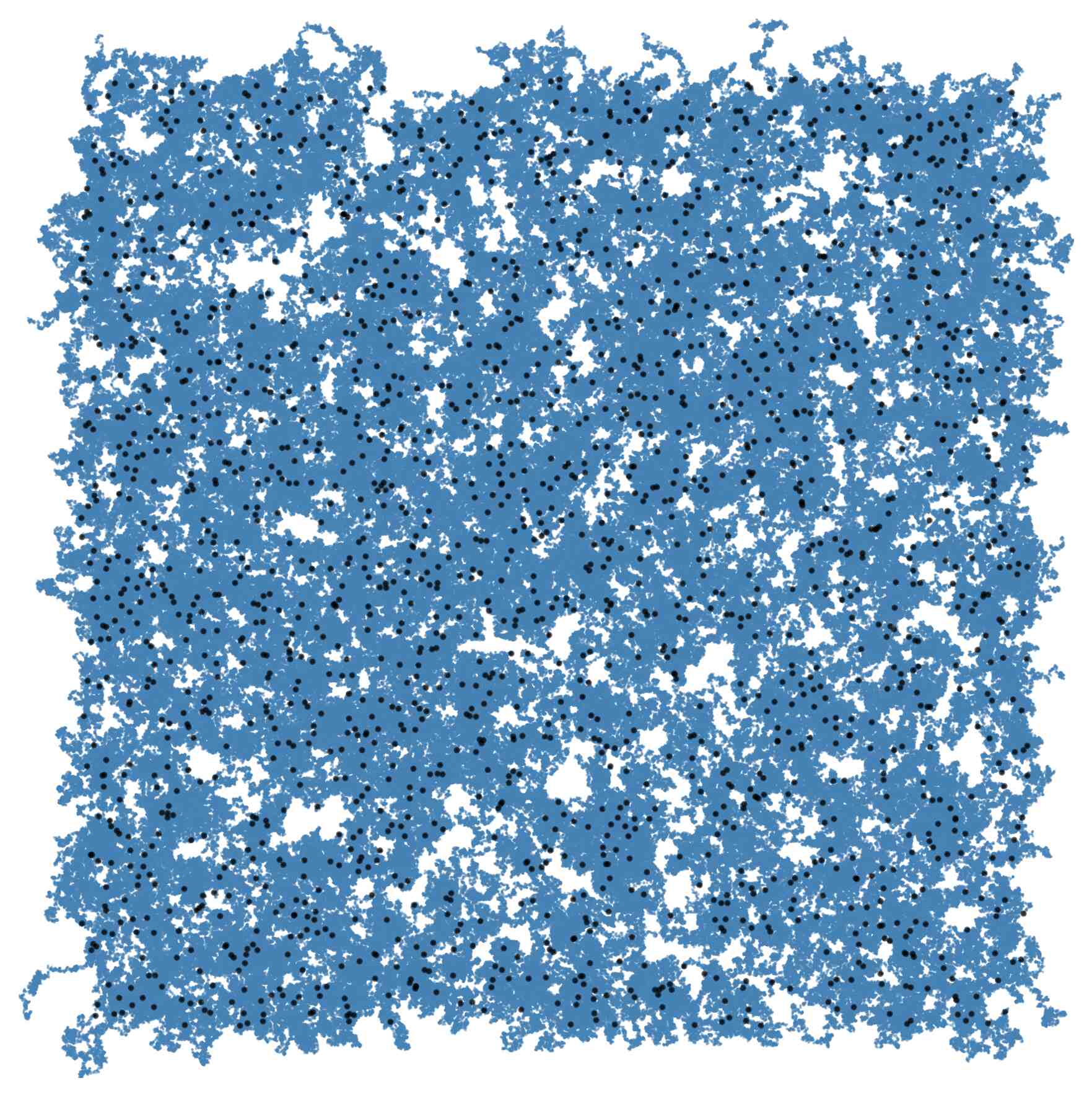}
    \centering
    \includegraphics[height=25mm]{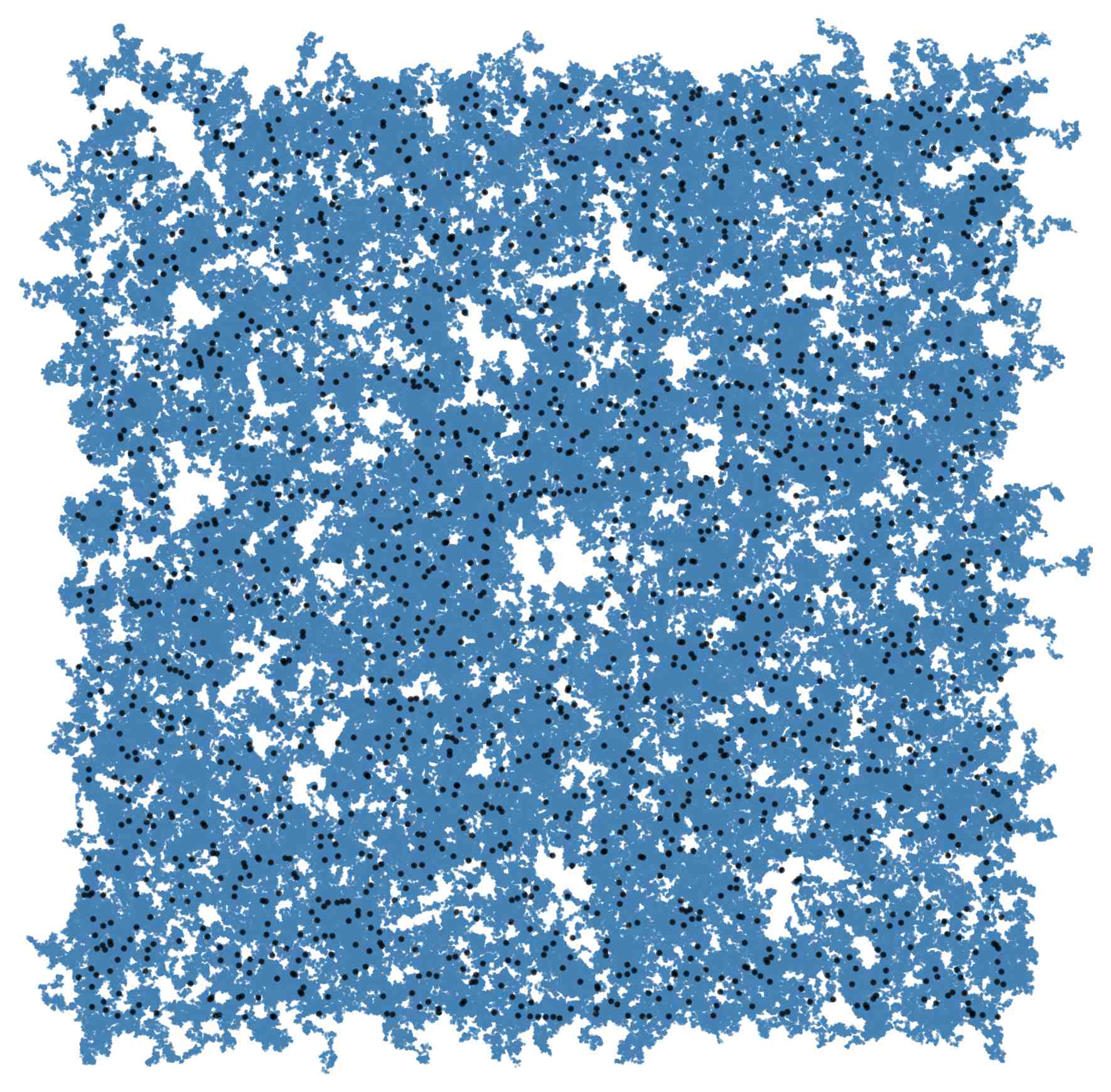}
    \centering
    \includegraphics[height=25mm]{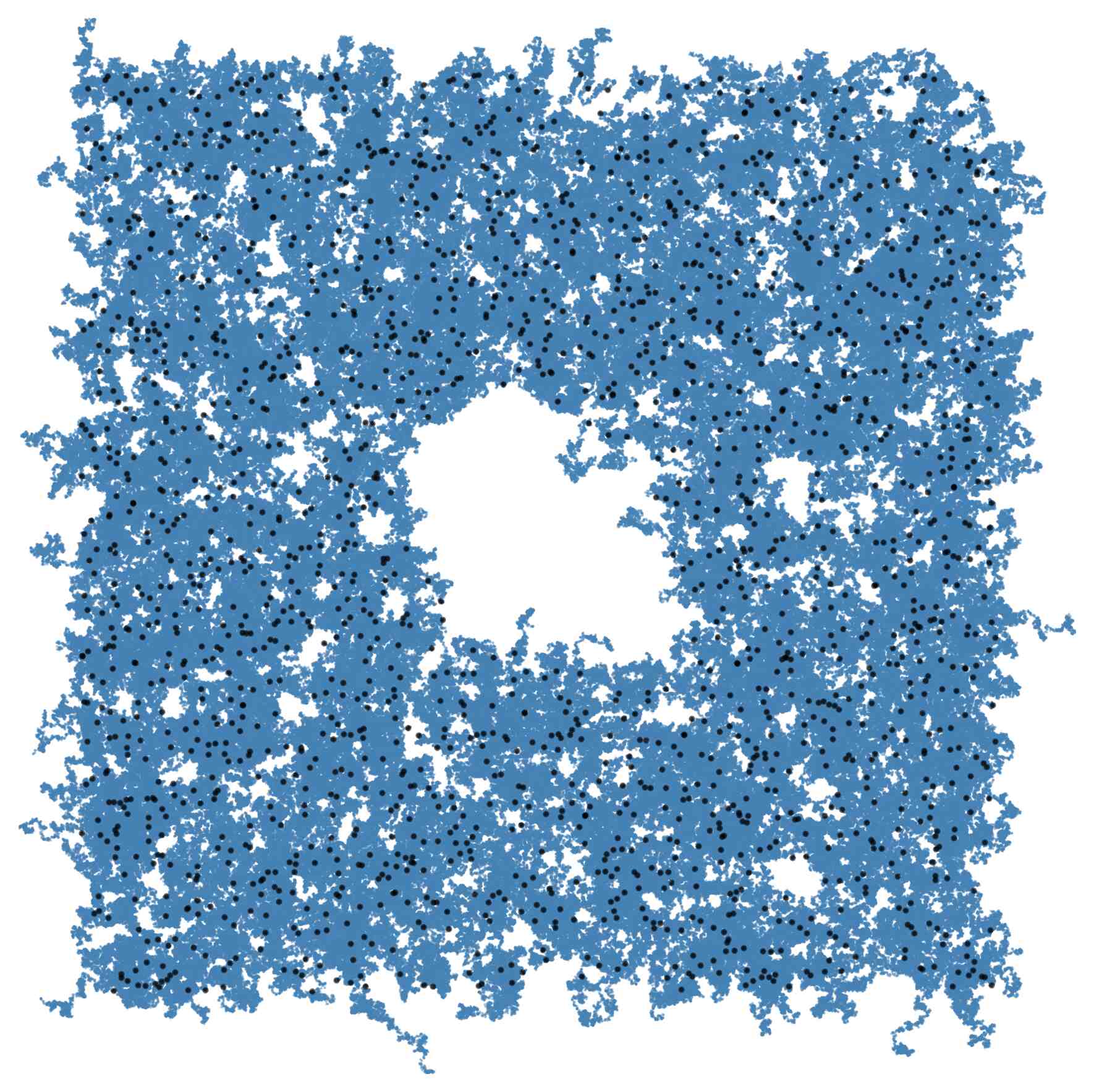}
    \centering
    \includegraphics[height=25mm]{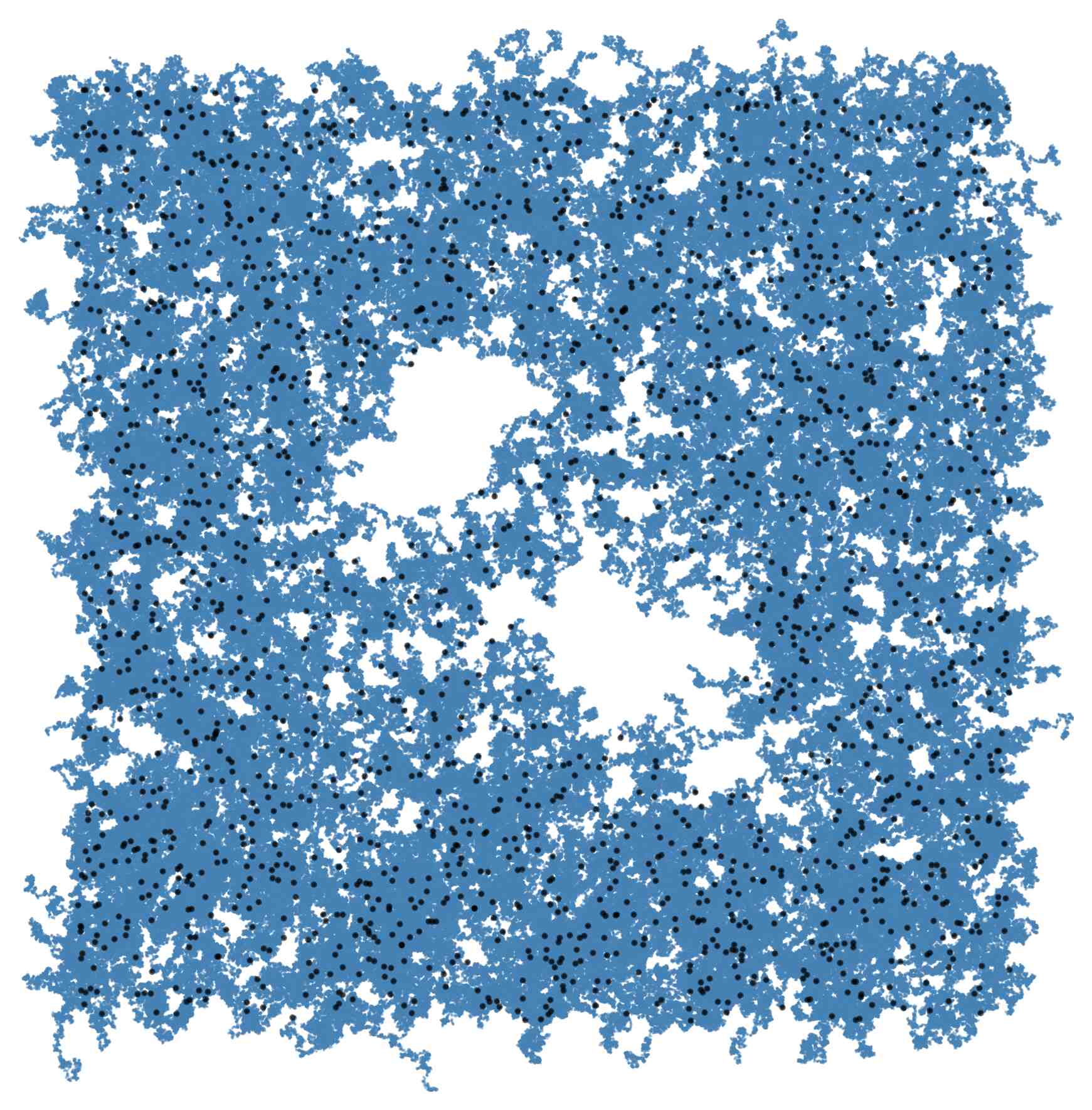}

    
    \centering
    \includegraphics[height=25mm]{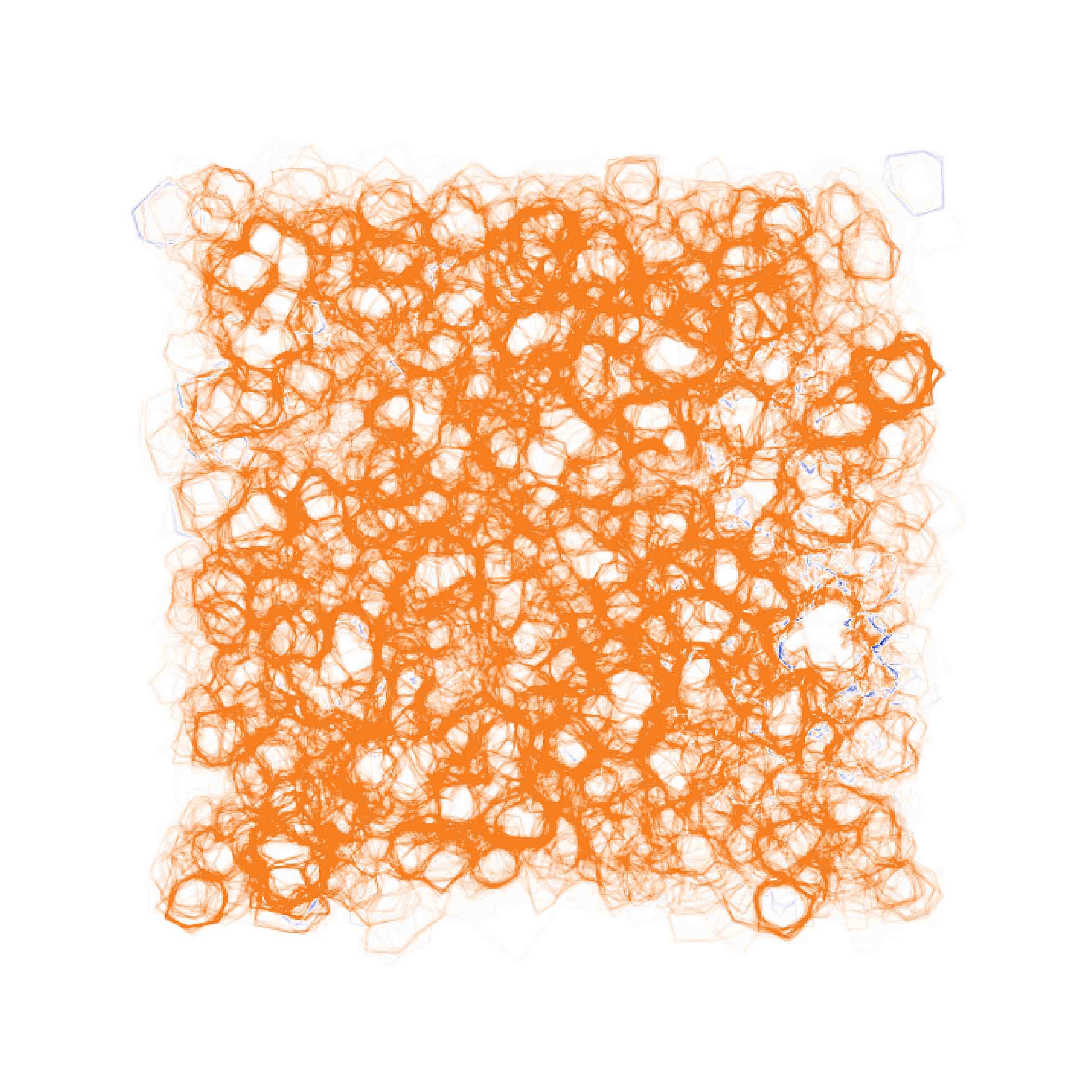}
    \centering
    \includegraphics[height=25mm]{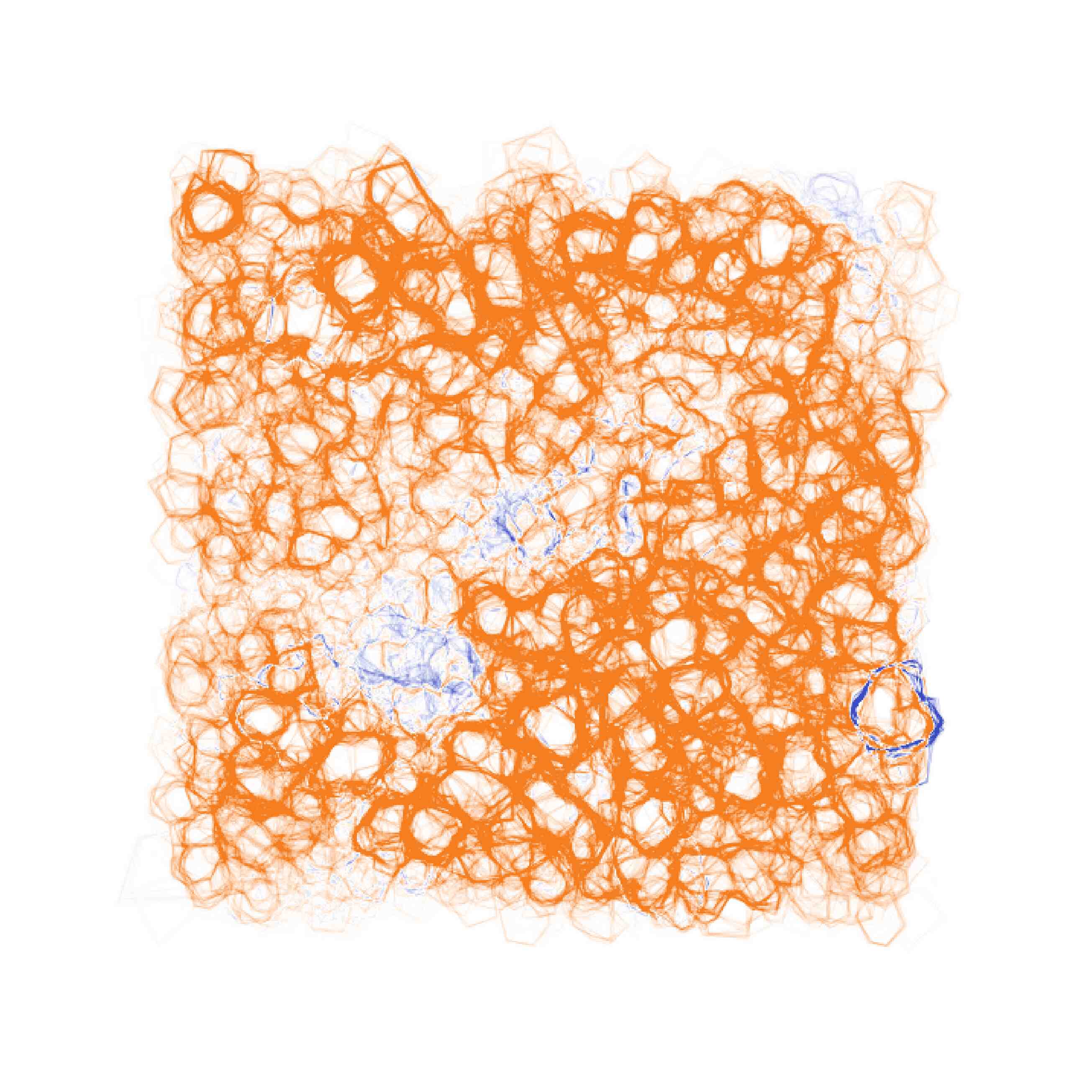}
    \centering
    \includegraphics[height=25mm]{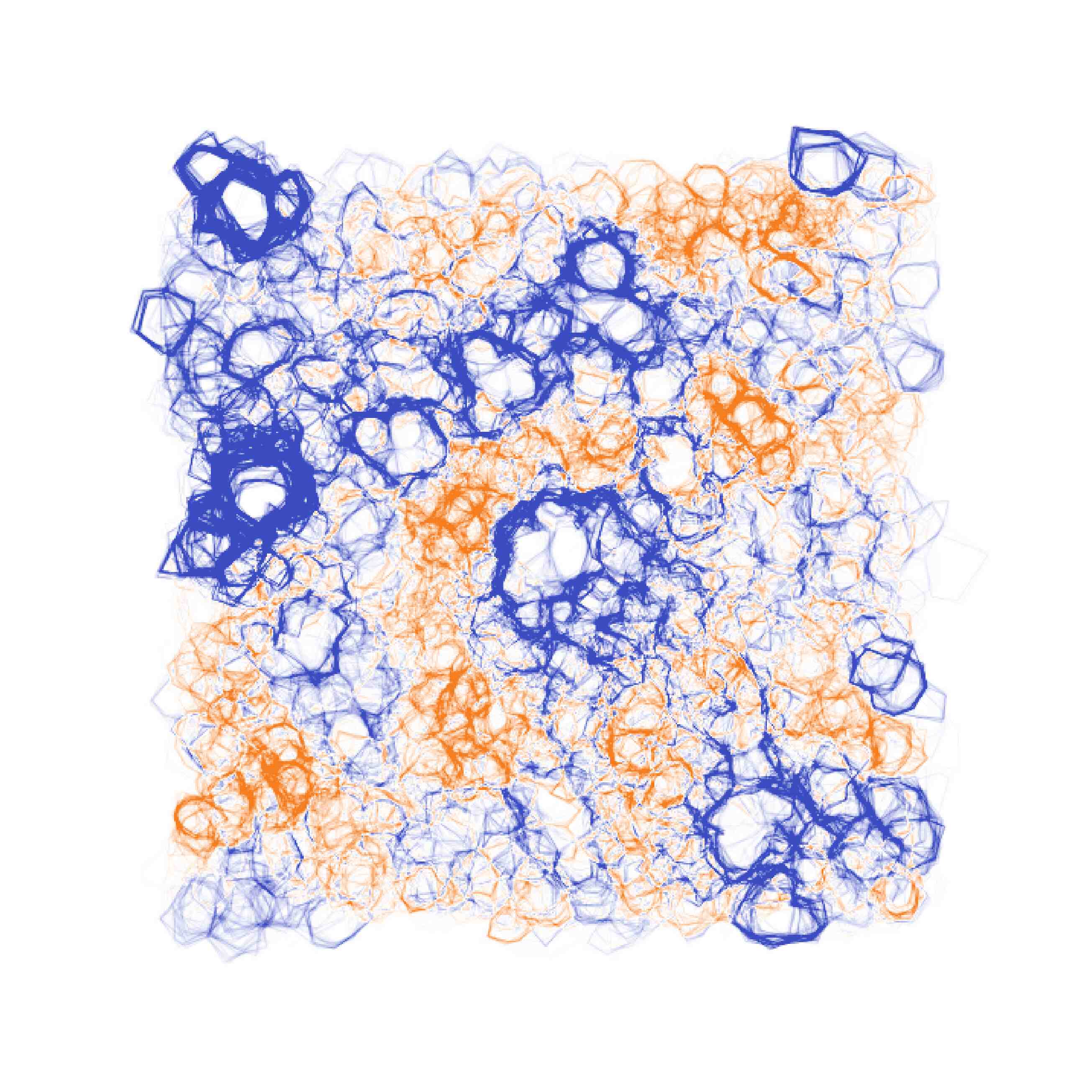}
    \centering
    \includegraphics[height=25mm]{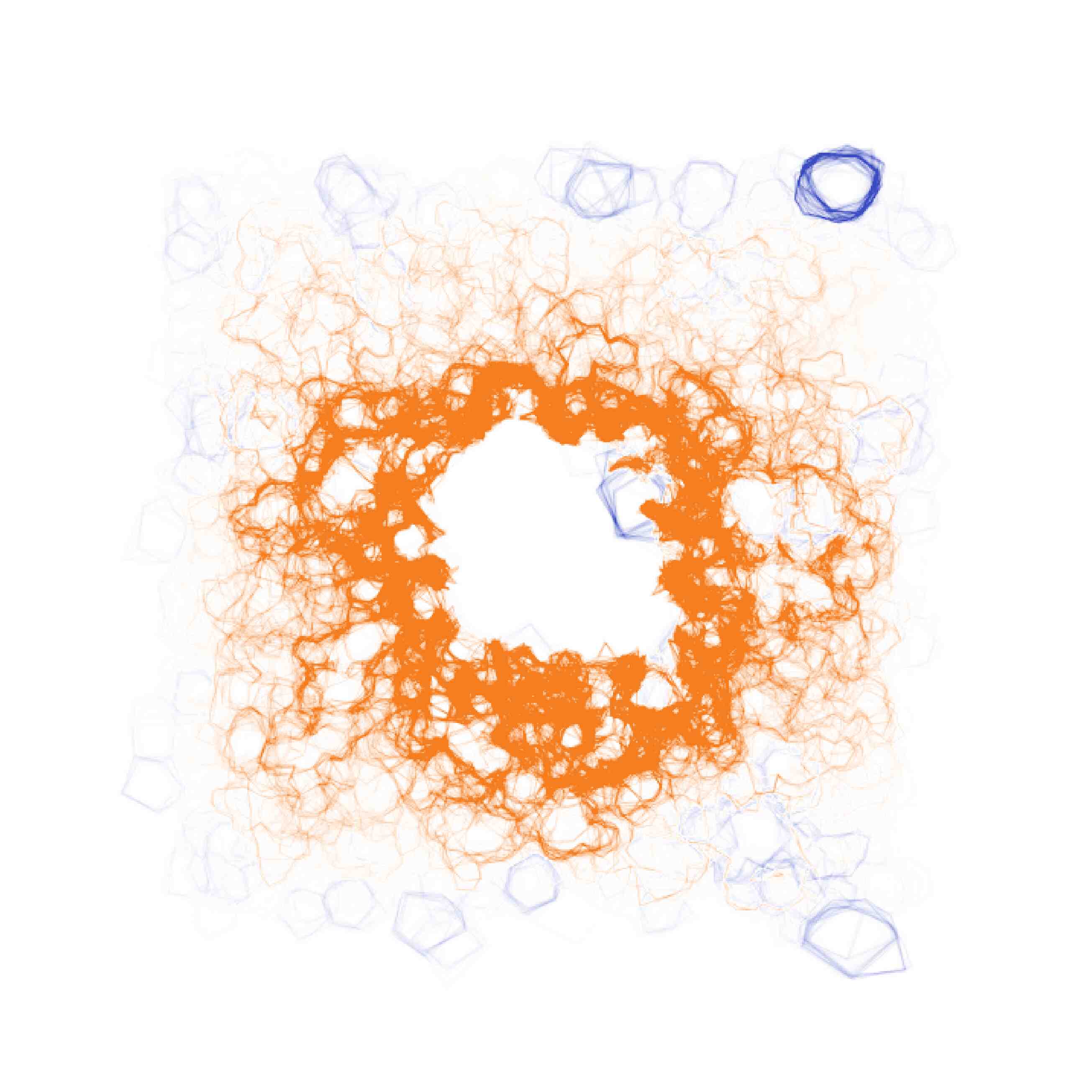}
    \centering
    \includegraphics[height=25mm]{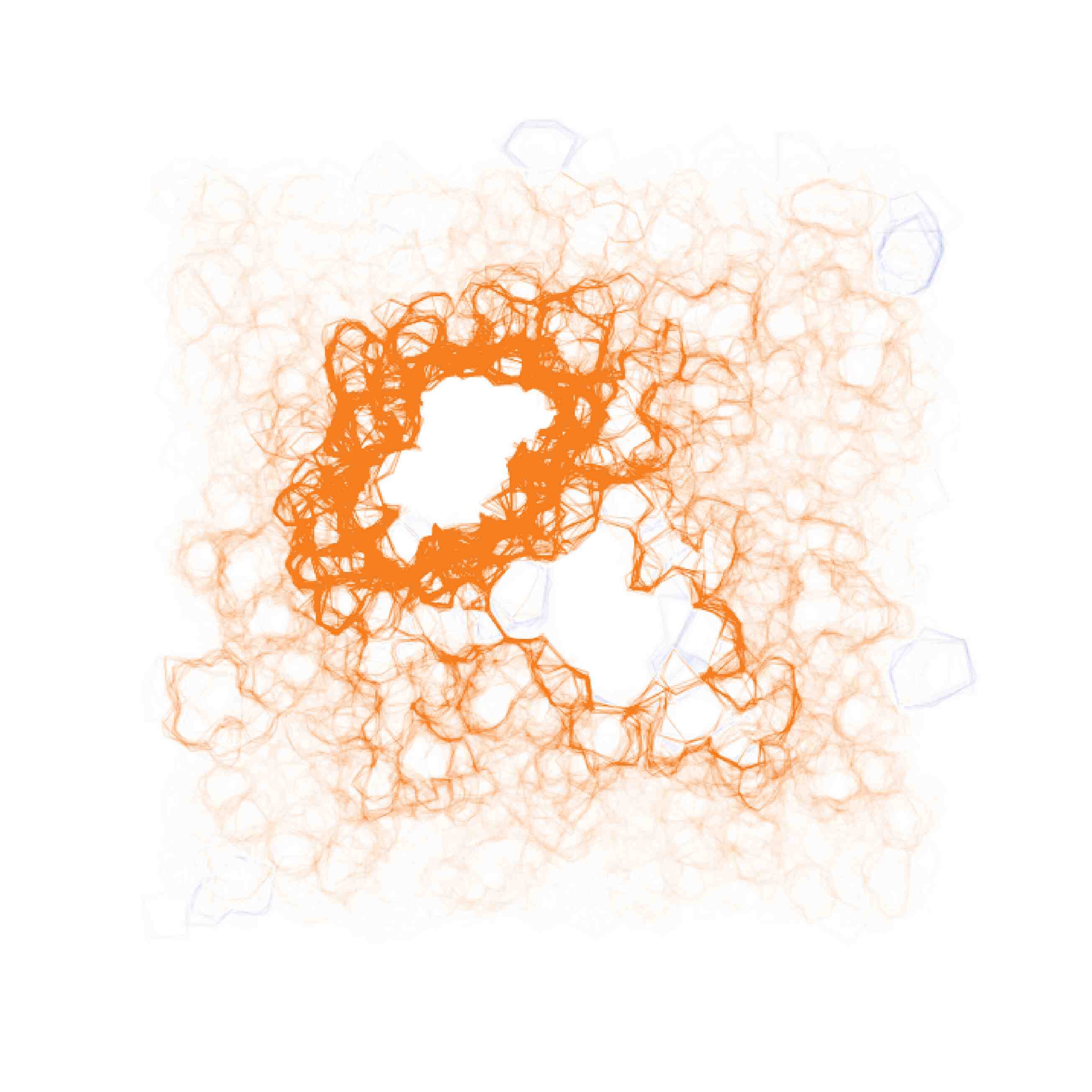}
    
    \caption{
    We consider the linked twist map for parameter values $r = 2.5, 3.5, 4.0, 4.1, 4.3$ (left to right).
    In the top row we plot the union of 1000 randomly initialized orbits (blue) overlaid with a single orbit (black).
    We compute persistence diagrams for homology in degree 1 and persistence landscapes, and use linear SVM to construct a classifier comparing one class to all others.
    The bottom row shows the empirical mean feature maps of $1000$ samples 
    using representative cycles, with a color gradient in which orange corresponds to the given class and blue corresponds to the other classes.
    In each image the sum of the pixel values equals the empirical mean value of the classifier.
    }
    \label{fig:linked-twist-visualization-deg1}
\end{figure}

For each $r \in \{2.5,3.5,4.0,4.1,4.3\}$, 
we generate $2,000$ points using \cref{eq:linked-twist-map}, with an initial condition drawn uniformly at random from the unit square $[0,1]^2$, repeating 100 times.
For simplicity we treat the points as a subset $X$ of the unit square.
As in the previous examples, we compute the Delaunay complex, the persistence diagram for homology in degree $1$ and the corresponding persistence landscape. 
For each class, we use linear SVM to produce a classification model by comparing that class to the union of the other classes.
We apply this model to 1000 samples from each class and use representative cycles 
and compute the empirical mean feature maps.
See \cref{fig:linked-twist-visualization-deg1}.

\bibliography{phm}

\end{document}